\documentclass[12pt]{article}
\usepackage{latexsym}
\usepackage{latexsym}
\usepackage{amsmath}
\usepackage{amsfonts}
\usepackage{amsbsy}
\usepackage{mathrsfs}
\usepackage[gray,pdflatex]{xcolor}
\usepackage{epsfig}
\usepackage[all]{xy}
\usepackage{url}
\DeclareMathOperator{\colim}{colim}
\DeclareMathOperator{\Co}{Co}
\DeclareMathOperator{\Ext}{Ext}
\newcommand{\C}{\mathcal{C}}
\begin{document}
\begin{titlepage}
\begin{center}
{\Large LOCALISATION AND COMPLETION}

\vspace{20pt}

with an\\
addendum on the use of\\
Brown-Peterson homology\\
in stable homotopy

\vspace{20pt}

by\\
J. F. Adams\\

University of Cambridge

\vspace{3.5in}

Lecture notes by Z. Fiedorowicz on a course\\
given at The University of Chicago in Spring 1973,\\
revised and supplemented by Z. Fiedorowicz, 2011



\end{center}
\end{titlepage}
\setcounter{page}{-2}
\pagestyle{empty}

\centerline{\bf FOREWORD}

\vspace{30pt}

In  spring 1973 Frank Adams gave a course at the University of Chicago on localisation and completion.  This
was in the very early days of the subject, which arose from disparate constructions of Quillen, Sullivan, Mimura, Nishida, Toda and others during the period 1969-1971.  In those days one usually assumed the spaces one was localising or completing were simply connected.  There were various proposals for extending the domain of definition of these constructions to more general spaces (e.g. \cite{BK}), but there was no clear consensus on how to proceed.

In his lectures Adams gave a lucid and compelling analysis of the properties one would want of such constructions.  He set up an elegant axiomatic treatment of localisation and completion in the framework of category theory and proposed a vast generalisation of the existing constructions. Unfortunately Adams'
program for constructing these localisation functors with respect to arbitrary generalised homology theories ran into a serious difficulty during the course of these lectures.  His proposal involved the use of the Brown Representability Theorem to construct his localisation functors, but he was unable to show that the relevant representable functors were set-valued rather than class-valued. Subsequent work by Bousfield established the existence of these generalised localisation functors, using more technical simplicial methods. These functors are now an essential tool in homotopy theory. 

At that time I was a graduate student at Chicago and was charged with the responsibility of taking notes for Adams' lectures. The resulting notes were briefly available in mimeographed form from the University of Chicago Mathematics Department. However these notes were never published in a more formal venue due to this apparent flaw in the proof.\footnote{An announcement of this work was published in \cite{A3}.} The notes also contain an addendum devoted to establishing that a certain element in the gamma family of the stable homotopy groups of spheres is nonzero, using Brown-Peterson (co)homology.  At that time this was a matter of controversy, as Oka and Toda claimed to have proved the contrary result.

I thought I had lost my only copy of these notes a long time ago, but I recently rediscovered them,
and I want to make them publicly available again. Besides being of historical interest, these notes give a very readable introduction to localisation and completion, with minimal prerequisites.  I was long aware that
the gap in Adams' proof was easily mendable, so I have supplemented these notes to include an epilogue
explaining this and have made a few other minor editorial changes. Thus it can now be seen in retrospect that Adams amazingly succeeded in his project of \textit{``constructing localizations and completions without doing a shred of work''} (cf. \cite{Ma}).

I would like to take this opportunity to thank Peter Landweber for his careful proofreading of an earlier
draft of this manuscript and Carles Casacuberta for providing me with the recent preprint \cite{BCMR}.

\vspace{10pt}

\hfill\parbox{155pt}{Zig Fiedorowicz\\ Columbus, Ohio\\ January, 2011\\fiedorow@math.ohio-state.edu}

\newpage

\centerline{CONTENTS}

\vspace{20pt}

1. Introduction to Localisation\dotfill\ \ \pageref{sec1}

2. Idempotent Functors\dotfill\ \pageref{sec2}

3. Axiomatic Characterisation of Classes S\dotfill\ \pageref{sec3}

4. A Further Axiom\dotfill\ \pageref{sec4}

5. Behaviour of Idempotent Functors with Respect to Fiberings;\\
$\begin{array}{c}\,\end{array}$
Construction of Localisation Using Postnikov Decomposition\dotfill\ \pageref{sec5}

6. Profinite Completion\dotfill\ \pageref{sec6}

7. Use of Brown-Peterson Homology in Stable Homotopy\dotfill\ \pageref{sec7}

8. Epilogue\dotfill\ \pageref{sec8}

$\begin{array}{c}\,\end{array}$
References\dotfill\ \pageref{refs}

\newpage

\pagestyle{plain}
\parindent15pt

\section{Introduction to Localisation}\label{sec1}

\bigskip

\hfill
\begin{minipage}{320pt}
\textit{``Il y a l\`a la possibilit\'e d'une \'etude locale (au sens arithm\'etique!)\\
des groupes d'homotopie \dots''}
\end{minipage}\\

\hfill
\begin{minipage}{90pt}
J-P. Serre \cite{Se}
\end{minipage}

\vspace{20pt}

In homotopy theory we have known for a long time that it is sufficient to attack problems one prime at a time. This insight goes back to the pioneering work of J-P. Serre \cite{Se}.

More recently we have gained a particularly convenient language and some particularly convenient machinery for exploiting this insight. This language and machinery was introduced following an analogy from commutative algebra. In commutative algebra we attack our problems one prime at a time by using the method of localisation. Thus we seek a comparable method in homotopy theory.

The earliest reference I have which develops such a method is Sullivan \cite{Su}. This was certainly very influential. At this point perhaps we should also mention Mimura-Nishida-Toda \cite{MNT}, Mimura-O'Neill-Toda \cite{MOT}, Mimura-Toda \cite{MT}, and Zabrodsky \cite{Z}. Another reference we might suggest is Quillen \cite{Q}.

In the first part of these lectures, I want to present a simple and uniform method of constructing all functors in homotopy theory which have formal properties similar to those of Sullivan's localisation functor. This opens the way to a study of such functors along axiomatic lines. I may also say something about Sullivan's completion functor. However it is clear from Sullivan's work that the completion functor enters it for a visibly good and sufficient reason which is particular to that piece of work. The localisation functor, however, is of very general use, and every graduate student of topology should learn about it.

I must begin by sketching some background, and I start with commutative algebra. Let $R$ be a commutative ring with 1 and
$S\subset R$ be a multiplicatively closed subset (i.e. a subset closed under finite products, such that $1\in S$). For example, if 
$R = \mathbf{Z}$ we may take
$$S=\{1, 2, 4, 8,\dots , 2^n, \dots\}$$
or
$$S = \{1, 3,5,\dots	2m+1, \dots \}.$$
Let $M$ be an $R$-module. We say $M$ is \textit{$S$-local} if the map $M\to M$ given by multiplication by $s$, i.e.,
$m \mapsto ms$, is an isomorphism for any $s\in S$.
To every $R$-module $M$ we can find a map $f:M\to M'$ so that
\begin{itemize}
\item[(i)] $M'$ is $S$-local
\item[(ii)] $f$ is universal with respect to (i). That is, if $g:M\to M''$ is another map such that $M''$ is $S$-local, then there is a
unique map $h: M'\to M''$ which makes the following diagram commute
$$\xymatrix{
&&M'\ar@{-->}[dd]^h\\
M\ar[urr]^f\ar[drr]_g\\
&&M''}$$
\end{itemize}
Such a map $f$ is called a localisation map; we say $f$ localises $M$ at $S$.

The usual construction of $M'$ is as a module of fractions. We first take pairs $(m,s)$, $m\in M$, $s\in S$.
We then define an equivalence relation on pairs:
$$(m,s)\sim(m ',s') \qquad\Longleftrightarrow\qquad\exists s''\in S\quad\ni\quad	ms's'' =m'ss''.$$
We define $S^{-1}M$ to be the set of all equivalence classes. The fraction
$\frac{m}{s}$ is the equivalence class containing $(m, s)$. We make $S^{-1}M$ into an
$R$-module in the obvious way. We define the map $f:M\to S^{-1}M$ by
$f(m) = \frac{m}{1}$ . We see that $S^{-1}M$ is $S$-local and the map $f$ is universal.

Since the ring $R$ is an $R$-module, $S^{-1}R$ is defined; we can make it into a ring so that the canonical map $R\to S^{-1}R$ is a map of rings:
$$\left(\frac{r}{s}\right)\left(\frac{r'}{s'}\right)=\frac{rr'}{ss'}.$$
Similarly $S^{-1}M$ becomes a module over $S^{-1}R$:
$$\left(\frac{m}{s'}\right)\left(\frac{r}{s}\right)=\frac{mr}{s's}.$$
Moreover we obtain a commutative diagram
$$\xymatrix{
M\otimes_R S^{-1}R\ar[rrd]^\cong\ar[dd]^\cong\\
&&S^{-1}M\\
S^{-1}M\otimes_R S^{-1}R\ar[urr]_\cong
}$$

We often use this fact just as a matter of notation when we have a convenient name for  $S^{-1}R$. For example, suppose $R=\mathbf{Z}$ and
$$S=\{1, 2, 4, 8,\dots , 2^n, \dots\}$$
so that $S^{-1}R=\mathbf{Z}\left[\frac{1}{2}\right]$. If $M$ is a $\mathbf{Z}$-module, we would usually
write $M\otimes\mathbf{Z}\left[\frac{1}{2}\right]$ for  $S^{-1}M$.

The most common example of a multiplicatively closed subset $S$ is the complement of a prime ideal $P$.
If $S =\mathcal{C}P$ we write $M_P$ for  $S^{-1}M$. 
For example, if $R=\mathbf{Z}$, $P=(2)$,
$$\mathcal{C}P=\{\pm 1,\pm 3,\pm 5,\dots,\pm(2m+1),\dots\},$$
then $\mathbf{Z}_{(2)}$ is the set of fractions
$$\left\{\frac{a}{2b+1}\right\}\subset\mathbf{Q}.$$

The construction of localisation in commutative algebra has many good properties.
Most of them need not delay us now; however, it is essential to know that
 localisation preserves exactness. If
$$L\stackrel{i}\longrightarrow M\stackrel{i}\longrightarrow N$$
is exact, then so is
$$S^{-1}L\stackrel{i}\longrightarrow S^{-1}M\stackrel{i}\longrightarrow S^{-1}N.$$

\noindent
{\bf Example (a)} Take the two localisation functors on $\mathbf{Z}$-modules
$$-\otimes \mathbf{Z}{\textstyle\left[\frac{1}{2}\right]}\mbox{ and }-\otimes\mathbf{Z}_{(2)}.$$
These two functors commute up to isomorphism and for any $\mathbf{Z}$-module $M$
we have a commutative diagram
$$\xymatrix{
M\ar[rrr]\ar[dd] &&&M\otimes\mathbf{Z}\left[\frac{1}{2}\right]\ar[d]\\
&&&\left(M\otimes\mathbf{Z}\left[\frac{1}{2}\right]\right)\otimes\mathbf{Z}_{(2)}\ar[d]^\cong\\
M\otimes\mathbf{Z}_{(2)}\ar[rrr]
&&&\left(M\otimes\mathbf{Z}_{2}\right)\otimes\mathbf{Z}\left[\frac{1}{2}\right]\cong M\otimes\mathbf{Q}
}$$
This diagram is both a pullback and pushout.

Conversely, if we are given : $M'$ a $\mathbf{Z}\left[\frac{1}{2}\right]$-module,
$M''$ a $\mathbf{Z}_{(2)}$-module ,
$M'''$ a $\mathbf{Q}$-module, and localising maps $f:M'\to M'''$, $g:M''\to M'''$,
then in the pullback diagram
$$\xymatrix{M\ar@{-->}[rr]^{h_1}\ar@{-->}[d]^{h_2}
&&M'\ar[d]^f\\
M''\ar[rr]^g &&M'''}$$
$h_1$ and $h_2$ are also localising maps.

\vspace{10pt}

Now I want to recall the basic theorem about Sullivan's localisation functor
and give one example of its use to show what it is meant for.

Before we do this, however, we must consider the category on which it is to be defined.
If we stick to simply-connected CW-complexes, everyone will feel happy and secure. Moreover,
it might be a matter of debate exactly how far we might wish to enlarge the domain of definition
of our functor; and if we take the domain too large, there might be more than one functor extending
the functor we all agree about for simply-connected spaces, and it might be a matter of debate
which extension is best.

So let us take $\C$ to be the category whose objects are 1-connected CW-complexes
with basepoint and whose morphisms are homotopy classes of maps, with both maps and homotopies
preserving basepoints.

Let $S\subset\mathbf{Z}$ be a multiplicatively closed subset. Then Sullivan showed there is a functor
$\C\to\C$ which at the level of homotopy and homology performs localisation at $S$.

\vspace{10pt}

\noindent
{\bf Theorem 1.1} (Sullivan). \textit{The following conditions on a map
$f:X\to Y$ in $\C$ are equivalent:
\begin{itemize}
\item[(i)] $f_{\#}:\pi_n(X)\to\pi_n(Y)$ localises at $S$ for each $n\ge1$.
\item[(ii)] $f_*:H_n(X)\to H_n(Y)$ localises at $S$ for each $n\ge1$.
\end{itemize}
Moreover there is a functor $E:\C\to\C$ and a natural transformation
${\eta:1\to E}$ so that for each $X$, $\eta_X:X\to EX$ satisfies both (i) and (ii).}

\vspace{10pt}

Before we go on, let me comment. Theorems like this usually tell us that there is something
with a stated property when it is not obvious that there is; but it is obvious that if there
is anything with that property, then the property characterises it. On the face of it this
is not a theorem of that form; the properties do characterise $EX$ and $\eta_X$, but it is
not obvious they do. Of course, these are comments on the way I've stated the theorem not
on the way Sullivan states it. Still we may make a note to look for other forms of the
statement.

Following the analogy from algebra we use the notation $X_S$ or $X\otimes\mathbf{Z}_S$
for $EX$. We defer discussing the properties of $E$ until we are forced to do so by examples.

As an application, I recall that at one time there was a conjecture of the following sort.

\vspace{10pt}

\noindent
{\bf Conjecture.} \textit{Any finite CW-complex which is an H-space is homotopy-equivalent
to a product of spaces from the classical list: $S^7$, $\mathbf{R}P^7$, compact Lie groups.}

\vspace{10pt}

This conjecture must always have looked optimistic, and it is now known to be false. The
first counterexample was due to Hilton and Roitberg. Let us see how we get one by
Zabrodsky's method of mixing homotopy groups, expressed in the language of localisation.

The classical list contains two entries $S^3\times S^7$ and $Sp(2)$.
These are different at the prime 2 and also at the prime 3 (e.g.
because $\pi_6(S^3\times S^7)=\mathbf{Z}/12\mathbf{Z}$, $\pi_6(Sp(2)) = 0$.)
However, if we apply the localisation functor $-\otimes\mathbf{Q}$ they become the same:
$$(S^3\times S^7)\otimes\mathbf{Q} = K(\mathbf{Q}, 3)\times K(\mathbf{Q},7) = Sp(2)\otimes\mathbf{Q}.$$
So consider for example
$$\xymatrix{
&&Sp(2)\otimes\mathbf{Z}\left[\frac{1}{2}\right]\ar[d]\\
&&Sp(2)\otimes\mathbf{Z}\left[\frac{1}{2}\right]\otimes\mathbf{Z}_{(2)}\ar[d]^\cong\\
(S^3\times S^7)\otimes\mathbf{Z}_{(2)}\ar[rr]
&&(S^3\times S^7)\otimes\mathbf{Z}_{(2)}\otimes\mathbf{Z}\left[\frac{1}{2}\right]
}$$
There should be a space $X$ such that
\begin{eqnarray*}
X\otimes\mathbf{Z}{\textstyle\left[\frac{1}{2}\right]} &\simeq 
&SP(2)\otimes\mathbf{Z}{\textstyle\left[\frac{1}{2}\right]}\\
X\otimes\mathbf{Z}_{(2)} &\simeq &Sp(2)\otimes\mathbf{Z}_{(2)}.
\end{eqnarray*}

It is easy to construct a candidate for $X$. Suppose we are given a diagram
$$\xymatrix{
&&X'\ar[d]^p\\
X''\ar[rr]^g &&X'''
}$$
then the weak pullback is the space of triples $(x',x'',\omega)$, where
$\omega:I\to X'''$ is a path from $p(x')$ to $g(x'')$. Strictly I shall
take a weakly equivalent CW-complex; anyway I get a diagram
$$\xymatrix{
X\ar[rr]\ar[d] &&X'\ar[d]^p\\
X''\ar[rr]^g &&X'''
}$$
and an exact homotopy sequence
$$\cdots\longrightarrow\pi_{n+1}(X''')\longrightarrow\pi_n(X)
\longrightarrow\pi_n(X')\oplus\pi_n(X'')\longrightarrow
\pi_n(X''')\longrightarrow\cdots$$
If we apply this to our case we get a diagram
$$\xymatrix{
X\ar[rr]\ar[d]&&Sp(2)\otimes\mathbf{Z}\left[\frac{1}{2}\right]\ar[d]\\
\left(S^3\times S^7\right)\otimes\mathbf{Z}_{(2)}\ar[rr]
&&K(\mathbf{Q}, 3)\times K(\mathbf{Q},7)
}$$
and we see that $X$ is 1-connected,
$$\pi_n(X)\longrightarrow \pi_n\left(Sp(2)\otimes\mathbf{Z}{\textstyle\left[\frac{1}{2}\right]}\right)$$
localises at $\mathbf{Z}\left[\frac{1}{2}\right]$ and
$$\pi_n(X)\longrightarrow\pi_n\left(\left(S^3\times S^7\right)\otimes\mathbf{Z}_{(2)}\right)$$
localises at $\mathbf{Z}_{(2)}$. In particular
$$\pi_6(X)=\mathbf{Z}/4\mathbf{Z}.$$
Therefore $X$ is not equivalent to $Sp(2)$ and it is not equivalent to $S^3\times S^7$;
it isn't equivalent to $Sp(2)$ at the prime 2 and it is not equivalent to $S^3\times S^7$
at the prime 3. Also it certainly isn't equivalent to anything else in the classical list.

Of course I should still convince you that X is an H-space and that it is equivalent
to a finite complex. To do the former, we need to note that localisation commutes with
products. More precisely, we have the projections
$$X\times Y\longrightarrow X,\qquad X\times Y\longrightarrow Y.$$
Localising we get
$$(X\times Y)_S\longrightarrow X_S,\qquad (X\times Y)_S\longrightarrow Y_S.$$
With these components we get
$$(X\times Y)_S\longrightarrow X_S\times Y_S.$$
This is an equivalence, because on homotopy it induces
$$\left(\pi_*(X)\oplus\pi_*(Y)\right)_S\stackrel{\cong}{\longrightarrow}\pi_*(X)_S\oplus\pi_*(Y)_S.$$
The product map
$$\left(S^3\times S^7\right)^2\longrightarrow S^3\times S^7$$
gives
$$\left(\left(S^3\times S^7\right)\otimes\mathbf{Z}_{(2)}\right)^2
\longrightarrow \left(S^3\times S^7\right)\otimes\mathbf{Z}_{(2)}$$
and
$$\left(\left(S^3\times S^7\right)\otimes\mathbf{Q}\right)^2
\longrightarrow \left(S^3\times S^7\right)\otimes\mathbf{Q}.$$
Similarly, the product map
$$\left(Sp(2)\right)^2\longrightarrow Sp(2)$$
gives
$$\left(Sp(2)\otimes\mathbf{Z}{\textstyle\left[\frac{1}{2}\right]}\right)^2
\longrightarrow Sp(2)\otimes\mathbf{Z}{\textstyle\left[\frac{1}{2}\right]}$$
and
$$\left(Sp(2)\otimes\mathbf{Q}\right)^2
\longrightarrow Sp(2)\otimes\mathbf{Q}.$$
The product maps on $Sp(2)\otimes\mathbf{Q}\simeq
\left(S^3\times S^7\right)\otimes\mathbf{Q}\simeq K(\mathbf{Q}, 3)\times K(\mathbf{Q},7)$
agree, since there is only one class of maps
$$\left(K(\mathbf{Q}, 3)\times K(\mathbf{Q},7)\right)^2\longrightarrow
K(\mathbf{Q}, 3)\times K(\mathbf{Q},7)$$
which has the basepoint as unit.

Hence we have the following diagram
$$\xymatrix @C=8pt{
X^2\ar[r]\ar@{-->}[drr]\ar[dd] 
&\left(Sp(2)\otimes\mathbf{Z}\left[\frac{1}{2}\right]\right)^2\ar[drr]\ar[dd]\\
&&X\ar[dd]\ar[r]&Sp(2)\otimes\mathbf{Z}\left[\frac{1}{2}\right]\ar[dd]\\
\left(\left(S^3\times S^7\right)\otimes\mathbf{Z}_{(2)}\right)^2\ar[r]\ar[drr]
&\left(K(\mathbf{Q}, 3)\times K(\mathbf{Q},7)\right)^2\ar[drr]\\
&&\left(S^3\times S^7\right)\otimes\mathbf{Z}_{(2)}\ar[r]
&K(\mathbf{Q}, 3)\times K(\mathbf{Q},7)
}$$
Since the right-hand square is a weak pullback, there is a map
$\mu:X^2\to X$ making the diagram commute. Moreover, consider the map
$$X\times\mbox{pt}\longrightarrow X\times X\stackrel{\mu}{\longrightarrow}X.$$
The induced map of homotopy groups becomes the identity if we apply
${-\otimes \mathbf{Z}\left[\frac{1}{2}\right]}$ or $-\otimes\mathbf{Z}_{(2)}$.
Therefore it is the identity. In particular this composite is an equivalence
by the theorem of J. H.C. Whitehead. Similarly
$$\mbox{pt}\times X\longrightarrow X\times X\stackrel{\mu}{\longrightarrow}X$$
is an equivalence. Hence we can alter $\mu$ so that these become the identity and
$X$ becomes an H-space.

Of course if we knew that $X$ was a strict pullback we wouldn't need that
last argument, but who cares.

Finally, we need to see that $X$ is equivalent to a finite complex.
Since $X$ is 1-connected it is sufficient to show that $\bigoplus_n H_n(X)$
is finitely-generated. From the fact that the maps
$$\xymatrix{
X\ar[rr]\ar[d]&&Sp(2)\otimes \mathbf{Z}\left[\frac{1}{2}\right]\\
\left(S^3\times S^7\right)\otimes\mathbf{Z}_{(2)}
}$$
induce localisation on homotopy, we infer they induce localisation on homology. Then we obtain
$$H_*(X)\cong H_*(Sp(2))\cong H_*(S^3\times S^7).$$

Before we go on, let me make one comment on the above example. If we want to construct
a fake Lie group then the original method of Hilton and Roitberg is simple and explicit,
and why should we use any other? This is certainly a reasonable objection. On the other hand,
if we want to state and prove that
$$(F/PL)\otimes\mathbf{Z}{\textstyle\left[\frac{1}{2}\right]}\simeq
BO\otimes\mathbf{Z}{\textstyle\left[\frac{1}{2}\right]}$$
then the method of Hilton and Roitberg can't do us much good. Since I actually want to talk about
localisation, I wanted to give a minimal example which would smuggle in a few of the things I wanted
to smuggle in. Thus I gave a minimal example in the direction of Zabrodsky mixing.

\newpage

\section{Idempotent Functors}\label{sec2}

I want to study functors in homotopy theory with the same formal properties
as localisation, so I'd better say what those properties are.

Suppose we are given a category $C$, a functor $E:C\to C$ and a natural
transformation $\eta: 1\to E$. On these I'm going to put two axioms:

\vspace{10pt}

\noindent
{\bf Axiom 2.1.}  $E\eta_X=\eta_{EX}:EX\to E^2X.$

\noindent
{\bf Axiom 2.2.} The common value of $E\eta_X$ and $\eta_{EX}$ 
is an equivalence from $EX$ to $E^2X$.

\vspace{10pt}

These axioms say that the functor $E$ is idempotent (in a particular way).
The categorists have already considered this axiom system, and they call
$(E,\eta)$ an idempotent triple or  idempotent monad.

\vspace{10pt}

\noindent
{\bf Example 2.3.} $C=R$-modules, $EM = S^{-1}M$ for a fixed
$S$, $\eta_M:M\to S^{-1}M$ is the canonical map $m\mapsto\frac{m}{1}$. Then
$$\eta_{EM}\mbox{ is the map } \frac{m}{s}\mapsto\frac{\frac{m}{s}}{1}$$
$$E\eta_{M}\mbox{ is the map } \frac{m}{s}\mapsto\frac{\frac{m}{1}}{s}.$$

\vspace{10pt}

\noindent
{\bf Example 2.4.} $C$ is the category in which the objects are metric spaces
and the maps are uniformly continuous functions. $EX$ is the completion $\hat{X}$
of $X$, constructed for example by taking equivalence classes of Cauchy sequences.
Also $\eta_X:X\to\hat{X}$ is the canonical map $x\mapsto\{x,x,\dots,x,\dots\}$.
$$\eta_{EX}\mbox{ is the map }
\{x_1,x_2,\dots\}\mapsto\{\{x_1,x_2,\dots\},\{x_1,x_2,\dots\},\dots\}$$
$$E\eta_{X}\mbox{ is the map }
\{x_1,x_2,\dots\}\mapsto\{\{x_1,x_1,\dots\},\{x_2,x_2,\dots\},\dots\}.$$

\vspace{10pt}

Since 2.1 and 2.2 are very simple formal properties, valid in Examples 2.3 and 2.4,
we may expect to see them hold in any case which presents a valid analogy with 2.3 or 2.4.
Equivalently suppose that for some $C$, $E$, $\eta$ either 2.1 or 2.2 is found to fail;
then that by itself would tend to discredit any analogy with 2.3 or 2.4 to the point
where we would hesitate to use the word ``localisation'' or ``completion'' for such a functor $E$.

Our programme is now as follows. First we must explore the consequences of our axioms. More
particularly we must understand how to characterise the map $\eta_X: X\to EX$
by universal properties. For this purpose we need to introduce two constructions. First one
needs to introduce a subset $D$ of the objects of $C$ and we must do it in such a way so that
for Example 2.4, $D$ becomes the subset of complete metric spaces (not just the spaces which
arise as $\hat{X}$ for some particular choice of completion). There are two equivalent
definitions.

\vspace{10pt}

\noindent
{\bf Definition 2.5.}\begin{itemize}
\item[(i)] $X\in D$ if $X\simeq EY$ for some $Y$ in $C$.
\item[(ii)] $X\in D$ if $\eta_X:X\to EX$ is an equivalence.
\end{itemize}

\vspace{10pt}

Clearly (ii) implies (i) (take $Y = X$). Also (i) implies (ii). For if
$X\stackrel{f}{\longrightarrow} EY$ is an equivalence, then $Ef$ is an equivalence and we have
$$\xymatrix{
EX\ar[rr]^{EF}_\simeq &&E^2 Y\\
X\ar[u]^{\eta_X}\ar[rr]^f_\simeq &&EY\ar[u]^\simeq_{\eta_{EY}}
}$$
so $\eta_X$ is also an equivalence.

We must also define a subclass $S$ of the morphisms of $C$.

\vspace{10pt}

\noindent
{\bf Definition 2.6.} The map $f:X\to Y$ lies in $S$ if and only if $Ef: EX\to EY$ is an equivalence.

\vspace{10pt}

The notation shows that this subclass of morphisms in a category is to be thought of as
analogous to a multiplicatively closed subset in a ring.

\vspace{10pt}

\noindent
{\bf Example 2.7.} Let $\C$ be the category in which the objects are all CW-complexes
with basepoint and the morphisms are homotopy classes of maps. When we construct Postnikov
systems, we may choose for each complex $X$ a complex $EX = X(1, 2, ... ,n)$ and a map
$\eta_X:X\to EX$ such that
\begin{itemize}
\item[(i)] $(\eta_X)_\#:\pi_i(X)\to\pi_i(EX)$ is an isomorphism for $i\le n$
\item[(ii)]$\pi_i(EX) = 0$ for $i>n$.
\end{itemize}
There is then a unique way to define $Ef$ for maps $f: X\to Y$ so that $E$ becomes a functor
$\C\to \C$ and $\eta$ becomes a natural transformation. Then
Axioms 2.1 and 2.2 are satisfied. We may call this $E$ the $n$-type functor. This shows that even in
homotopy theory there are functors $E$ satisfying 2.1 and 2.2 which are very different from localisation.
In this example it is particularly easy to describe the subsets $S$ and $D$. A map $f:X\to Y$ is in $S$ if and only if
$$f_\#:\pi_i(X)\to\pi_i(X)\mbox{ is an isomorphism for } i\le n.$$
An object $X$ is in $D$ if and only if
$$\pi_i(X)=0 \mbox{ for }i>n.$$

Our plan is now this. For those categories $C$ which concern us, we will characterise by axioms those subclasses $S$ which arise from pairs $(E,\eta)$ by Definition 2.6. In the applications we can then proceed as follows: first we write a definition which defines a subclass $S$; secondly we check that this subclass
$S$ satisfies the axioms; and finally we apply the main theorem to deduce that this subset $S$ arises from a pair $(E,\eta)$ satisfying 2.1 and 2.2.

Before we go on, let me make a few comments on why topologists should study idempotent functors on categories which interest them. First the general method of algebraic topology has the following pattern. We wish to study a category $C$ which is not properly known, like the homotopy category. So we apply functors like homotopy and homology which go from $C$ to categories which are known a bit better. Now if we are ever to get to knowing $C$, certainly we must begin by knowing some of those full subcategories $D$ which are a bit simpler. Once we do, it becomes fruitful
to consider functors on $C$ taking values in $D$ and apply our standard methods.

Secondly, by studying idempotent functors on $C$ we gain a great deal of understanding of the structure of
$C$. This is in good analogy with the study of idempotents in a ring or idempotents acting on a module.

Thirdly, the examples suggest that this study is reasonable.

Let us now define the notion of equivalence for idempotent functors. The appropriate equivalence relation is as follows:

\vspace{10pt}

\noindent
{\bf Definition 2.8.} If $(E,\eta)$, $(E',\eta')$ are two idempotent functors on
a category $C$ we say that $(E,\eta)\simeq(E',\eta')$, if there is a natural equivalence
$\epsilon: E\to E'$ such that $\eta'=\epsilon\eta$, i.e.
$$\xymatrix{
&&EX\ar[dd]^{\epsilon_X}\\
X\ar[urr]_{\eta_X}\ar[drr]^{\eta'_X}\\
&&E'X
}$$
If $(E,\eta)\simeq(E',\eta')$, then $S=S'$, $D=D'$; it is trivial to check it.

\vspace{10pt}

The categorists know some things about idempotent monads, and I must run through some of them.
The first thing you know about any monad is that the functor factors through a pair of adjoint functors.
In our case it is very easy to display the factorisation. On the one hand we can regard $E$ as taking values
in the full subcategory $D$
$$C\stackrel{E}{\longrightarrow}D\qquad D\stackrel{I}{\longrightarrow}C$$
and on the other hand we have the inclusion of $D$ in $C$. I must prove that these are adjoint functors.

\vspace{10pt}

\noindent
{\bf Lemma 2.9.}\textit{If $X\in C$, $Y\in D$ then
$$[X,Y]\stackrel{\eta_X^*}{\longleftarrow}[EX,Y]$$
is an isomorphism.}

\vspace{10pt}

\noindent
{\bf Proof.} (i) $\eta_X^*$ is epi: Suppose $f: X\to Y$. From the diagram
$$\xymatrix{
EX\ar[rr]^{EF} &&EY\\
X\ar[u]^{\eta_X}\ar[rr]^f &&Y\ar[u]_{\eta_Y}^\simeq
}$$
we see that $f$ factors through $\eta_X$.

\noindent
(ii) $\eta_X^*$ is mono: Take $f,g:EX\to Y$ and assume that $f\cdot\eta_X = g\cdot\eta_X$.
Then $Ef\cdot E\eta_X=Eg\cdot E\eta_X$. But since $\eta_{EX}=E\eta_X$, it follows from the diagram
$$\xymatrix{
E^2X\ar@<1ex>[rr]^{Ef}\ar@<-1ex>[rr]_{Eg} &&EY\\
\\
EX\ar[uu]^{\eta_{EX}}\ar@<1ex>[rr]^{f}\ar@<-1ex>[rr]_{g} &&Y\ar[uu]_{\eta_Y}
}$$
that $\eta_Y\cdot f=\eta_Y\cdot g$. Since $\eta_Y$ is an equivalence, $f=g$.\hfill$\Box$

\vspace{10pt}

The next proposition shows how $S$ and $D$ determine each other.

\vspace{10pt}

\noindent
{\bf Proposition 2.10.} \textit{(i) Suppose $f:X\to Y$ is a morphism in $C$. Then $f$
lies in $S$ if and only if
$$f^*:[Y,Z]\longrightarrow [X,Z]$$
is an isomorphism for all $Z$ in $D$.\\
(ii) If $Z$ is an object in $C$, then $Z$ is in $D$ if and only if 
$$f^*:[Y,Z]\longrightarrow [X,Z]$$
is an isomorphism for all morphisms $f:X\to Y$ in $S$. Indeed it is sufficient
to check that $f^*$ is epi for all $f$ in $S$.}

\vspace{10pt}

\noindent
{\bf Proof.} We must prove that if $f:X\to Y$ lies in $S$ and $Z$ lies in $D$, then
$$f^*:[Y,Z]\longrightarrow [X,Z]$$
is an isomorphism. This is immediate from the following commutative diagram
$$\xymatrix{
[Y,Z]\ar[rr]^{f^*} &&[X,Z]\\
[EY,Z]\ar[u]^{\eta_Y^*}_\simeq\ar[rr]^{(Ef)^*} &&[EX,Z]\ar[u]^{\eta_X^*}_\simeq
}$$
and since $\eta_Y^*$ and $\eta_X^*$ are isomorphisms 
by 2.9, while $(Ef)^*$ is an isomorphism because $f\in S$.

(i) Conversely suppose
$$f^*:[Y,Z]\longrightarrow [X,Z]$$
is an isomorphism for all $Z$ in $D$; we wish to show that $f\in S$.
By the same diagram as above, we see that
$$(Ef)^*:[EY,Z]\longrightarrow [EX,Z]$$
is an isomorphism for all $Z\in D$. Now we argue in the standard fashion.
Taking $Z = EX$ we see that there is a $g:EY\to EX$ such that 
$$g\cdot Ef = 1_{EX}.$$
Then
$$Ef\cdot g\cdot Ef = Ef\cdot 1_{EX}=1_{EY}\cdot Ef.$$
But taking $Z =EY$ and using the fact that $(Ef)^*$ is an isomorphism we get
$$Ef\cdot g=1_{EY}.$$
Hence $Ef$ is an equivalence and $f\in S$.

(ii) Suppose $f^*:[Y, Z]\to[X,Z]$ is epi for all $f$ in $S$; we wish to
show $Z\in D$. Take $f=\eta_Z:Z\to EZ$. This lies in $S$ because $Ef= E\eta_Z$ is
an equivalence by Axiom 2.2. Hence
$$\eta_Z^*:[EZ,Z]\to [Z,Z]$$
is epi and there is a map $\zeta:EZ\to Z$ such that
$$\zeta\cdot\eta_Z=1_Z.$$
Then
$$\eta_Z\cdot\zeta\eta_Z=\eta_Z\cdot 1_Z=1_Z\cdot\eta_Z.$$
But
$$\eta_Z^*:[EZ,EZ]\to [Z,EZ]$$
is an isomorphism by 2.9, since $EZ\in D$. Therefore
$$\eta_Z\cdot\zeta=1_{EZ},$$
$\eta_Z$ is an equivalence, and $Z\in D$.\hfill$\Box$

\vspace{30pt}

The next proposition shows that the pair $(E,\eta)$ is determined up to equivalence
by either $S$ or $D$.

\vspace{30pt}

\noindent
{\bf Proposition 2.11.} \textit{The following conditions on a map $f:X\to Y$ are equivalent.
\begin{itemize}
\item[(i)] There is a commutative diagram
$$
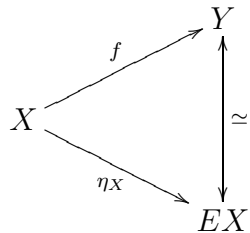
\xymatrix{
&&Y\ar@{<->}[dd]^\simeq\\
X\ar[urr]^f\ar[drr]_{\eta_X}\\
&&EX}$$
In other words, $f$ is a map $\eta_X$ up to equivalence.
\item[(ii)] $f$ is in $S$ and $Y$ is in $D$.
\item[(iii)]$f$ is in $S$ and is couniversal with respect to that property.
That is, given a map $s:X\to X'$ in $S$, there is a unique map $h:X'\to Y$
which makes the following diagram commute.
$$
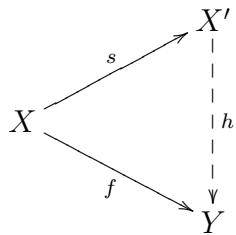
\xymatrix{
&&X'\ar@{-->}[dd]^h\\
X\ar[urr]^s\ar[drr]_f\\
&&Y}$$
\end{itemize}
}

\vspace{30pt}

\textit{
\begin{itemize}
\item[(iv)] $Y$ is in $D$ and is universal with respect to that property.
That is, if $g:X\to Z$ is a map with $Z\in D$, then there is a unique map
$h:Y\to Z$ which makes the following diagram commute.
$$
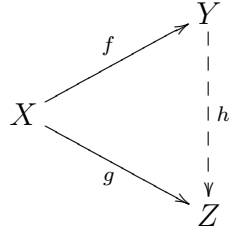
\xymatrix{
&&Y\ar@{-->}[dd]^h\\
X\ar[urr]^f\ar[drr]_g\\
&&Z}$$
\end{itemize}
}

\vspace{10pt}

\noindent
{\bf Proof.} (i) $\Longrightarrow$ (ii) Suppose we have a commutative diagram
$$
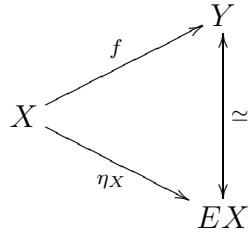
\xymatrix{
&&Y\ar@{<->}[dd]^\simeq\\
X\ar[urr]^f\ar[drr]_{\eta_X}\\
&&EX}$$
Applying E we obtain the commutative diagram
$$
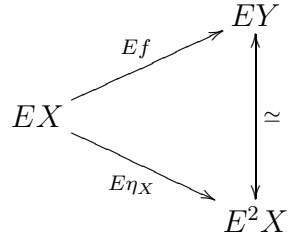
\xymatrix{
&&EY\ar@{<->}[dd]^\simeq\\
EX\ar[urr]^{Ef}\ar[drr]_{E\eta_X}\\
&&E^2X}$$
But $E\eta_X$ is an equivalence by 2.2. Therefore $Ef$ is an equivalence and $f\in S$. Also $Y\in D$
since $Y\simeq EX$.

(ii)$\Longrightarrow$ (iii) Suppose $f$ is in $S$ and $Y$ is in $D$. Suppose $s:X\to X'$ is in $S$. Then
$$s^*:[X',Y]\longrightarrow[X,Y]$$
is an isomorphism by 2.10, so $f$ has the couniversal property required.

(ii) $\Longrightarrow$ (iv) Suppose $f$ is in $S$ and $Y\in D$. Let $g:X\to Z$ with $Z\in D$ be given. Then
$$f^*: [Y, Z]\longrightarrow [X, Z]$$
is an isomorphism by 2.10, so $f$ has the universal property required.

(iii) $\Longrightarrow$ (i) Suppose $f$ is in $S$ and has the couniversal property. We know $\eta_X$ is in $S$
and we have shown it has the couniversal property, so $f$ and $\eta_X$ differ by a canonical equivalence:
$$\xymatrix{
&&Y\ar@{<->}[dd]^\simeq\\
X\ar[urr]^f\ar[drr]_{\eta_X}\\
&&EX}$$
Similarly for (iv) $\Longrightarrow$ (i).\hfill$\Box$

\vspace{10pt}

\noindent
{\bf Remark.} In Example 2.7, I gave a description of the $n$-type functor.
If you look back you will see that this is a description of the type 2.11(ii) 
-- $\eta_X$ is in $S$ and $EX$ is in $D$. In the situation we wish to study,
such a description is adequate to characterise $E$ and $\eta$.

\vspace{10pt}

You might think this was enough of \S 2, and so did I when I first wrote this lecture,
but then I found there was some more material I needed here. In fact, I must first recall
some of the standard material on the category of fractions.

Suppose we are given a category $C$ and a subset $S$ of the morphisms in $C$.
Then there is a category $S^{-1}C$ and a functor $Q:C\to S^{-1}C$ with the following properties:
\begin{itemize}
\item[(i)] If $s\in S$, $Qs$ is invertible in $S^{-1}C$;
\item[(ii)] $Q$ is universal with respect to property (i). That is, suppose we
are given a functor $T:C\to C'$ such that $Ts$ is invertible for any $s\in S$;
then there exists a unique functor $U:S^{-1}C\to C'$ such that $T = UQ$, i.e.,
$$\xymatrix{
&&S^{-1}C\ar@{-->}[dd]^U\\
C\ar[urr]^Q\ar[drr]_T\\
&&C'}$$
\end{itemize}

In fact, to construct these things, you take $S^{-1}C$ to have the same
objects as $C$ and you take $Q$ to be the identity on objects. To construct
the morphisms of $S^{-1}C$ you start with long words
$$\cdots(Qs_2)^{-1}(Qf_2)(Qs_1)^{-1}(Qf_1)$$
 which are grammatical, i.e.,
$$\xymatrix{
&Y_1 &&Y_2\\
X_1\ar[ur]^{f_1} &&X_2\ar[ul]^{s_1}\ar[ur]^{f_2} &&X_3\ar[ul]^{s_2} &\dots &.
}$$
You then divide these words into equivalence classes in an obvious way and check that everything works.

\vspace{10pt}

\noindent
{\bf N.B.} In practice, however, the collection of morphisms $S$ we wish to invert is not a set
but a class. In this case there are serious set-theoretic objections against the legitimacy of our
constructing the category $S^{-1}C$, e.g. for general classes $S$ and $X,Y\in C$,
$\mbox{Hom}_{S^{-1}C}(QX,QY)$ need not be a set. To overcome these objections we need a condition
analogous to the solution set condition in Freyd's adjoint functor theorem.

\vspace{10pt}

\noindent
{\bf Lemma 2.12.} \textit{Let $T_0, T_1:S^{-1}C\to\Gamma$ be functors and let
${\beta:T_0Q\to T_1Q}$ be a natural transformation. Then there is one and only
one natural transformation $\alpha:T_0\to T_1$ such that $\beta_X=\alpha_{QX}$
for each $X$ in $C$.}

\noindent
{\bf Proof.} (1) First proof. Since $Q$ is a one-to-one correspondence on objects,
the equation $\beta_X=\alpha_{QX}$ defines $\alpha$; we have only to check that it
is natural. We are given that the diagram
$$\xymatrix{
T_0QX\ar[rr]^{\alpha_{QX}}\ar[d]_{T_0g} &&T_1QX\ar[d]^{T_1g}\\
T_0QY\ar[rr]^{\alpha_{QX}} &&T_1QY
}$$
commutes when $g$ is of the form $Qf$, $f$ in C. Therefore it commutes
when $g$ is of the form $(Qs)^{-1}$. Therefore it commutes when g is of the
form
$$\cdots(Qs_2)^{-1}(Qf_2)(Qs_1)^{-1}(Qf_1).$$

(2) Second proof: We construct a category $\Gamma^I$ as follows: 
An object in $\Gamma^I$ is a map $f:X_0\to X_1$ in $\Gamma$. A morphism in
$\Gamma^I$ from $f:X_0\to X_1$ to $g:Y_0\to Y_1$ is a square diagram
of the following form
$$\xymatrix{
X_0\ar[rr]^f\ar[d]_{h_0} &&X_1\ar[d]^{h_1}\\
Y_0\ar[rr]^g &&Y_1
}$$

There are functors $\pi_i:\Gamma^I\to\Gamma$ for $i=0,1$; the value of $\pi_i$ on the
square diagram displayed above is $h_i:X_i\to Y_i$. The natural transformation
$\beta:T_0Q\to T_1Q$ may be interpreted as a functor $\beta:C\to\Gamma^I$
such that $\pi_i\beta=T_iQ$. Similarly,a natural transformation $\alpha:T_0\to T_1$
may be interpreted as a functor $\alpha:S^{-1}C\to\Gamma^I$ such that
$\pi_i\alpha=T_i$. It can be
checked that our problem is now to factorise the functor $\beta$ in the form
$\beta=\alpha Q$. This factorisation is possible and unique by the universal property
 of $Q$. More precisely, let $s\in S$; then $T_iQs$ is invertible in $\Gamma$.
That is, $\beta s$ is a square diagram whose two vertical arrows are invertible in $\Gamma$;
so $\beta s$ is invertible in $\Gamma^I$. Therefore $\beta$ factors uniquely through $Q$.
\hfill$\Box$

\vspace{10pt}

Conditions are known under which we can get a much better hold on $S^{-1}C$.
Among them are the following:

\noindent
{\bf Axiom 2.13.} The set $S$ is closed under finite compositions.

\noindent
{\bf Axiom 2.14.} Given any diagram
$$\xymatrix{
X\\
W\ar[u]^s\ar[r]^f &Y
}$$
with $s\in S$, there is a diagram
$$\xymatrix{
X\ar[rr]^g &&Z\\
W\ar[u]^s\ar[rr]^f &&Y\ar[u]_t
}$$
with $t\in S$ and $gs = tf$.

\vspace{10pt}

Axiom 2.14 allows us to rewrite $(Qf)(Qs)^{-1}$ as $(Qt)^{-1}(Qg)$; with the use of
Axiom 2.13 it allows us to reduce every long word
$$\cdots(Qs_2)^{-1}(Qf_2)(Qs_1)^{-1}(Qf_1)$$
to a short word
$$(Qs)^{-1}(Qf)$$
so that every morphism in $S^{-1}C$ can be represented as a short word.

\vspace{10pt}

\noindent
{\bf Axiom 2.15.} Given any diagram
$$\xymatrix{
W\ar[r]^s &X\ar@<1ex>[r]^f\ar@<-1ex>[r]_g &Y
}$$
with $s\in S$ and $fs = gs$, there exists a diagram
$$\xymatrix{
W\ar[r]^s &X\ar@<1ex>[r]^f\ar@<-1ex>[r]_g &Y\ar[r]^t &Z
}$$
with $t\in S$ and $tf = tg$.

\vspace{10pt}

Using Axiom 2.15 we can tell if two short words are equivalent.
Suppose the short words are $(Qs_1)^{-1}(Qf_1)$, $(Qs_2)^{-1}(Qf_2)$ so that
$$\xymatrix{
&&Y_1\ar@{-->}[drr]^{g_1}\\
X\ar[urr]^{f_1}\ar[drr]_{f_2} &&Y\ar[u]_{s_1}\ar[d]^{s_2}\ar@{-->}^{s_3}[rr] &&Y_3\\
&&Y_2\ar@{-->}[urr]_{g_2}
}$$
Then the condition is that there exists $s_3:Y\to Y_3$ in $S$ and maps $g_i:Y_i\to Y_3$ for $i=1,2$ such that
$$g_1f_1= g_2f_2\mbox{\quad and\quad }g_1s_1=g_3=g_2s_2.$$

For our purposes, it is sufficient to assume that $S$ satisfies 2.13, 2.14 and 2.15; hence we can construct the category of fractions $S^{-1}C$ as follows: Take as objects the objects of $C$. For morphisms from $X$ to $Y$, first take the diagrams
$$
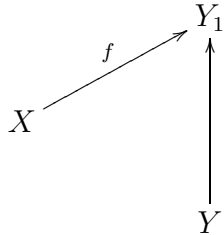
\xymatrix{
&&Y_1\\
X\ar[urr]^f\\
&&Y\ar[uu]
}$$
Define two diagrams
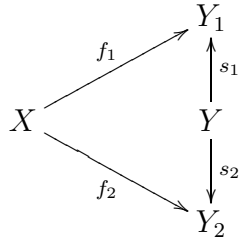
$$\xymatrix{
&&Y_1\\
X\ar[urr]^{f_1}\ar[drr]_{f_2} &&Y\ar[u]_{s_1}\ar[d]^{s_2}\\
&&Y_2
}$$
to be equivalent if there is a diagram
$$
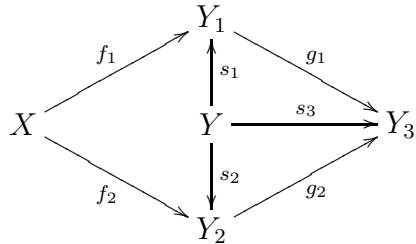
\xymatrix{
&&Y_1\ar[drr]^{g_1}\\
X\ar[urr]^{f_1}\ar[drr]_{f_2} &&Y\ar[u]_{s_1}\ar[d]^{s_2}\ar[rr]^{s_3} &&Y_3\\
&&Y_2\ar[urr]_{g_2}
}$$
such that $g_1f_1= g_2f_2$ and $g_1s_1=g_3=g_2s_2$. We check that this is an equivalence relation and take the morphisms from $X$ to $Y$ to be the equivalence classes. We define composition appropriately, check that we get a category. Finally we define $Q$ in the obvious way and check that it has the stated properties. See Gabriel and Zisman \cite{GZ}.

\newpage

\section{Axiomatic Characterisation of Classes $S$}\label{sec3}

Suppose we are given a subclass $S$ of the morphisms in a category $C$.
Then we inquire whether it satisfies the following six axioms, and I intend to prove
that they characterise the subsets $S$ which arise by Definition 2.6 at least for the categories which concern us.

\vspace{10pt}

\noindent
{\bf Axiom 3.1.} The class $S$ is closed under finite compositions.

\vspace{10pt}

\noindent
{\bf Axiom 3.2.} Given any diagram
$$\xymatrix{
X\\
W\ar[u]^s\ar[r]^f &Y
}$$
with $s$ in $S$, there is a diagram
$$\xymatrix{
X\ar[rr]^g &&Z\\
W\ar[u]^s\ar[rr]^f &&Y\ar[u]_t
}$$
with $t\in S$ and $gs = tf$.

\vspace{10pt}

\noindent
{\bf Axiom 3.3.} Given any diagram
$$\xymatrix{
W\ar[r]^s &X\ar@<1ex>[r]^f\ar@<-1ex>[r]_g &Y
}$$
with $s\in S$ and $fs = gs$, there exists a diagram
$$\xymatrix{
W\ar[r]^s &X\ar@<1ex>[r]^f\ar@<-1ex>[r]_g &Y\ar[r]^t &Z
}$$
with $t\in S$ and to $tf = tg$.

\vspace{10pt}

The next axiom is a solution set condition necessary for constructing the quotient category
$S^{-1}C$.

\noindent
{\bf Axiom 3.4.} To each object $Y$ in $C$ there is a set of arrows
$\left\{Y\stackrel{s_\alpha}{\longrightarrow}Z_\alpha\right\}$ in $S$ which are cofinal in $S$,
 i.e., given any $s:Y\to Z$, there is an arrow
$s_\alpha:Y\to Z_\alpha$ and a map $f:Z\to Z_\alpha$ such that $fs=s_\alpha$
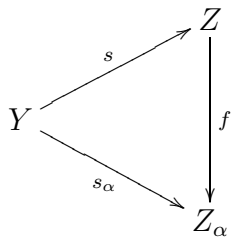
$$\xymatrix{
&&Z\ar[dd]^f\\
Y\ar[urr]^s\ar[drr]_{s_\alpha}\\
&&Z_\alpha
}$$

\vspace{10pt}

In algebra several subsets of a ring may give the same localisation. For instance if $R=\mathbf{Z}$
and
\begin{eqnarray*}
S_1 &= &\{1, 2, 4, 8,\dots , 2^n, \dots\}\\
S_2 &= &\{1, 4, 16, 64,\dots , 4^n, \dots\}
\end{eqnarray*}
then $S_1^{-1}R=S_2^{-1}R$. To simplify matters we usually consider the biggest such $S$.
The following axiom imposes a condition of a similar sort on our class $S$.

\vspace{10pt}

\noindent
{\bf Axiom 3.5.} If $f$ is a morphism in $C$ such that $Qf$ is invertible in $S^{-1}C$,
then $f\in S$.

\vspace{10pt}

\noindent
{\bf Remark.} Strictly speaking, Axiom 3.5 renders Axiom 3.1 redundant, but 3.1 is so elementary it seems foolish not to state it first.

Axiom 3.5 also admits the following equivalent formulation.

\vspace{10pt}

\noindent
{\bf Axiom $\mbox{3.5}^*$.} Given maps
$$W\stackrel{f}{\longrightarrow} X\stackrel{g}{\longrightarrow}Y\stackrel{h}{\longrightarrow}Z$$
with $gf$ and $hg$ in $S$, then $g\in S$. 

\vspace{10pt}

This version of 3.5 has the advantage that it makes perfect sense without assuming 3.4.

For the next axiom, I assume that my category $C$ has arbitrary coproducts; since I am interested
in the case of homotopy theory, I write the coproduct of factors $\{X_\alpha\}$ as 
$\bigvee_\alpha X_\alpha$.

\vspace{10pt}

\noindent
{\bf Axiom 3.6.} If $s_\alpha:X_\alpha\to Y_\alpha$ lies in $S$ for each $\alpha$, then
$$\bigvee_\alpha s_\alpha:\bigvee_\alpha X_\alpha\longrightarrow \bigvee_\alpha Y_\alpha$$
lies in $S$.

\vspace{10pt}

The easier half of our work is as follows.

\vspace{10pt}

\noindent
{\bf Proposition 3.7.} \textit{Let $C$ be a category with arbitrary coproducts,
$(E,\eta)$ a pair satisfying 2.1, 2.2, and $D$, $S$ as constructed in Definitions 2.5-2.6.
Then $S$ satisfies 3.1 to 3.6 inclusive. Moreover, $E$ factors as
$$C\stackrel{Q}{\longrightarrow} S^{-1}C\longrightarrow D\stackrel{I}{\longrightarrow}C$$
where $S^{-1}C\to D$ is an equivalence.}

\vspace{10pt}

The harder half is as follows.

\vspace{10pt}

\noindent
{\bf Theorem 3.8.} \textit{Let $\C$ be the category in which the objects are connected CW-complexes
with basepoint and the maps are homotopy classes. Let $S$ be a subclass of the morphisms of $\C$,
satisfying 3.1 to 3.6 inclusive. Then $S$ arises by 2.6 from a pair $(E,\eta)$ satisfying 2.1 and 2.2.}

\vspace{10pt}

\noindent
{\bf Note.} I really need only one thing about $\C$. I want to use Brown's Representability Theorem.
Therefore the same result is true if we take
$\C$ to be the category of 1-connected complexes or the category of spectra, etc.

\vspace{10pt}

In the categories which categorists usually consider one can probably get through even more easily by using standard adjoint functor theorems, but one presumably has to modify the axioms accordingly.

\vspace{10pt}

\noindent
{\bf Proof of Proposition 3.7.} To prove Axiom 3.1 suppose
$$X_1\stackrel{f_1}{\longrightarrow} X_2\stackrel{f_2}{\longrightarrow}X_3\longrightarrow\cdots
\longrightarrow X_n\stackrel{f_n}{\longrightarrow}X_{n+1}$$
all lie in $S$, i.e., $Ef_1$, $Ef_2$, \dots, $Ef_n$ are all equivalences. Then
$$E(f_nf_{n-1}\cdots f_2f_1)=(Ef_n)(Ef_{n-1})\cdots(Ef_2)(Ef_1)$$
is an equivalence, i.e., $f_nf_{n-1}\cdots f_2f_1$ lies in $S$.

(ii) Suppose we are given the diagram
$$\xymatrix{
X\\
W\ar[u]^s\ar[r]^f &Y
}$$
with $s$ in $S$. Then we can form the diagram
$$\xymatrix{
&&EX\\
X\ar[urr]^{\eta_X} &&EW\ar[u]^\simeq_{Es}\ar[rr]^{Ef} &&EY\\
W\ar[u]^s\ar[urr]^{\eta_W}\ar[rrr]^f &&&Y\ar[ur]^{\eta_Y}
}$$
Then $(Ef)(Es)^{-1}\eta_X:X\to EY$ is a map such that
$$\left[(Ef)(Es)^{-1}\eta_X\right]s=\eta_Yf,$$
and $\eta_Y\in S$, since $E\eta_Y$ is an equivalence.

(iii) To prove Axiom 3.3, suppose we are given
$$\xymatrix{
W\ar[r]^s &X\ar@<1ex>[r]^f\ar@<-1ex>[r]_g &Y
}$$
with $s\in S$ and $fs = gs$. Applying $E$ we obtain
$(Ef)(Es)=(Eg)(Es)$. Since $Es$ is an equivalence, $Ef=Eg$. But then
$(Ef)\eta_X=(Eg)\eta_X$ so by the diagrams
$$
\xymatrix{
X\ar[rr]^f\ar[d]_{\eta_X} &&Y\ar[d]^{\eta_Y}\\
EX\ar[rr]^{Ef} &&EY}
\qquad
\xymatrix{
X\ar[rr]^g\ar[d]_{\eta_X} &&Y\ar[d]^{\eta_Y}\\
EX\ar[rr]^{Eg} &&EY}
$$
we obtain $\eta_Yf=\eta_Yg$. Since $\eta_Y\in S$, this proves 3.3.

(iv) Take the set $\{s_\alpha:Y\to Z_\alpha\}$ to consist of the single map
$\eta_Y:Y\to EY$. Then given any $s:Y\to Z$ in $S$, we have the diagram
$$\xymatrix{
Z\ar[rr]^{\eta_Z} &&EZ\\
Y\ar[u]^s\ar[rr]^{\eta_Y} &&EY\ar[u]^\simeq_{Es}
}$$
whence $\eta_Y=(Es)^{-1}\eta_Zs$. Since $\eta_Y\in S$, $\eta_Y:Y\to EY$ is
cofinal in $S$ and $S$ satisfies Axiom 3.4.

(v) To prove Axiom 3.5, we first note that by definition, if $s\in S$, then
$Es$ is invertible. Hence $E$ factors through $S^{-1}C$:
$$\xymatrix{
C\ar[rr]^E\ar[dr]_Q &&C\\
&S^{-1}C\ar[ur]_R
}$$
Now suppose $f$ is a morphism in $C$ such that $Qf$ is invertible. Then
${RQf=Ef}$ is invertible, so $f\in S$.

(vi) To prove Axiom 3.6 suppose that $s_\alpha:X_\alpha\to Y_\alpha$ lies in $S$ for each $\alpha$.
Then
$$s_\alpha^*:[Y_\alpha, Z]\longrightarrow [X_\alpha, Z]$$
is an isomorphism for any $Z\in D$, by 2.10. Consider the following diagram
$$\xymatrix{
\left[\bigvee_\alpha Y_\alpha,Z\right]\ar[rr]^{(\bigvee_\alpha s_\alpha)^*}\ar[d]_\simeq
&&\left[\bigvee_\alpha X_\alpha,Z\right]\ar[d]^\simeq\\
\prod_\alpha[Y_\alpha,Z]\ar[rr]^{\prod_\alpha s_\alpha^*}_\simeq
&&\prod_\alpha[X_\alpha,Z]
}$$
We conclude that $(\bigvee_\alpha s_\alpha)^*$ is an isomorphism for any $Z\in D$.
Therefore $\bigvee_\alpha s_\alpha$ lies in $S$, by 2.10.

It is clear that the functor $E:C\to C$ factors through $D$ considered as a full subcategory
of $C$. Since $f\in S$ implies that $Ef$ is invertible the functor $C\to D$ factors through
$S^{-1}C$. It remains to prove that $S^{-1}C\to D$ is an equivalence of categories.

(i) Certainly every object of $D$ is equivalent to an object $EX$, i.e., an object in the image
of $C$ or $S^{-1}C$.

(ii) Take a map $f:EX\to EY$ in $D$. Then we have a diagram
$$\xymatrix{
E^2X\ar[rr]^{Ef} &&E^2Y\\
EX\ar[u]^{E\eta_X=E\eta_X}_\simeq\ar[rr]^f
&&EY\ar[u]_{E\eta_Y=E\eta_Y}^\simeq
}$$
Therefore
$$\xymatrix{
EX\ar[rr]^f &&EY\\
X\ar[u]^{\eta_X} &&Y\ar[u]_{\eta_Y}
}$$
represents an element of $S^{-1}C$ whose image in $D$ is the given element.

(iii) Take two morphisms of $S^{-1}C$ which have the same images in $D$, say represented by diagrams
$$
\xymatrix{
&&Y_1\\
X\ar[urr]^{f_1}&&Y\ar[u]_{s_1}
}\qquad\qquad\qquad
\xymatrix{
&&Y_2\\
X\ar[urr]^{f_2}&&Y\ar[u]_{s_2}
}
$$
I claim that without loss of generality we may suppose $s_1=s_2$. For by 3.2 we can construct
the following diagram
$$\xymatrix{
Y_1\ar[rr] &&Y_3\\
Y\ar[rr]^{s_2}\ar[urr]^s\ar[u]^{s_1} &&Y_2\ar[u]_{s_3}
}$$
with $s_3\in S$ so $s=s_3s_2\in S$. Then the two elements are equally well represented by
$$
\xymatrix{
&&Y_3\\
&Y_1\ar[ur]\\
X\ar[ur]^{f_1}&&Y\ar[uu]_s
}\qquad\qquad\qquad
\xymatrix{
&&Y_3\\
&Y_2\ar[ur]\\
X\ar[ur]^{f_2}&&Y\ar[uu]_s
}
$$
So suppose $s_1=s_2$ and revert to the notation $f_1$, $f_2$. Now the condition
that the two elements have the same images in $D$ is
$$(Es)^{-1}(Ef_1)=(Es)^{-1}(Ef_2)$$
whence 
$$Ef_1=Ef_2.$$
Now we can make the argument that
$$\eta_{Y_1}f_1=(Ef_1)\eta_X=(Ef_2)\eta_X=\eta_{Y_1}f_2$$
and $\eta_{Y_1}\in S$ whence the diagram
$$\xymatrix{
&Y_1\ar[drr]^{\eta_{Y_1}}\\
X\ar[ur]^{f_1}\ar[dr]_{f_2} &Y\ar[u]_s\ar[d]^s\ar[rr]^{(Es)\eta_Y} &&EY_1\\
&Y_1\ar[urr]_{\eta_{Y_1}}
}$$
shows that the two given elements are equal in $S^{-1}C$.\hfill$\Box$

\vspace{10pt}
Before we start the proof of Theorem 3.8 let's just observe one thing.
When we proved Axiom 3.2, we took $t$ to be the map $\eta_Y$ independent of
all the other data. But when we want to use Axiom 3.2, it only delivers
a map $t$ depending on all the other data. However, we can use Axiom 3.6
to mitigate this effect.

\vspace{10pt}

\noindent
{\bf Lemma 3.9.} \textit{Suppose there is given a set of diagrams 
$$\xymatrix{
X_\alpha\\
W_\alpha\ar[u]^{s_\alpha}\ar[r]^{f_\alpha} &Y
}$$
with a common Y and $s_\alpha$ in $S$. Then there exists
a single map $t:Y\to Z$ in $S$ so that we can complete all the diagrams
$$\xymatrix{
X_\alpha\ar[rr]^{g_\alpha} &&Z\\
W_\alpha\ar[u]^{s_\alpha}\ar[rr]^{f_\alpha} &&Y\ar[u]_t
}$$
}

\noindent
{\bf Proof.} Consider the diagram
$$\xymatrix{
\bigvee_\alpha X_\alpha\\
\bigvee_\alpha W_\alpha\ar[u]^{\bigvee_\alpha s_\alpha}\ar[r]^{\quad\{f_\alpha\}} &Y
}$$
By 3.6, $\bigvee_\alpha s_\alpha$ lies in $S$, so by 3.2 we can fill the diagram to get
$$\xymatrix{
\bigvee_\alpha X_\alpha\ar[rr]^{\quad g_\alpha} &&Z\\
\bigvee_\alpha W_\alpha\ar[u]^{\bigvee_\alpha s_\alpha}\ar[rr]^{\quad\{f_\alpha\}} &&Y\ar[u]^t
}$$
It is clear that $t$ has the required properties.\hfill$\Box$

\vspace{10pt}

\noindent
{\bf Corollary 3.10.} \textit{Suppose there is given a set of maps
$\{s_\alpha:Y\to X_\alpha\}$ in $S$. Then there exists a single map
$t:Y\to Z$ in $S$ which factors through each $s_\alpha$.}

\vspace{10pt}

\noindent
{\bf Proof.} This is a special case of 3.9 with $W_\alpha=Y$ and 
$f_\alpha=1_Y$ each $\alpha$.\hfill$\Box$

\vspace{10pt}

We operate similarly on Axiom 3.3.

\vspace{10pt}

\noindent
{\bf Lemma 3.11.} \textit{Suppose there is given a set of diagrams
$$\left\{\xymatrix{
X_\alpha\ar@<1ex>[r]^{f_\alpha}\ar@<-1ex>[r]_{g_\alpha} &Y
}\right\}$$
such that $Qf_\alpha=Qg_\alpha$ in $S^{-1}C$ for each $\alpha$. Then there exists
a single map $t:Y\to Z$ in $S$ such that $tf_\alpha=tg_\alpha$ for each $\alpha$.}

\vspace{10pt}

\noindent
{\bf Proof.} There is for each $\alpha$ a map $s_\alpha:Y\to Z_\alpha$ such that 
$s_\alpha f_\alpha=s_\alpha g_\alpha$. By Corollary 3.10 there is a single map $t:Y\to Z$ in $S$
which factors through each $s_\alpha$. Then $tf_\alpha=tg_\alpha$ for each $\alpha$.
\hfill$\Box$

\vspace{10pt}

\noindent
{\bf Lemma 3.12.} \textit{A coproduct $\bigvee_\alpha X_\alpha$ in $C$ remains a coproduct in
$S^{-1}C$.}

\vspace{10pt}

\noindent
{\bf Proof.} Let $i_\alpha:X_\alpha\to\bigvee_\alpha X_\alpha$ be the canonical injections;
then $Qi_\alpha:QX_\alpha\to Q\left(\bigvee_\alpha X_\alpha\right)$
induces a map
$$(Qi_\alpha)^*:\left[Q\left(\bigvee_\alpha X_\alpha\right),Y\right]_{S^{-1}C}
\longrightarrow [QX_\alpha,Y]_{S^{-1}C}$$
for any $Y$ in $S^{-1}C$. We have to show that
$$\left\{(Qi_\alpha)^*\right\}:\left[Q\left(\bigvee_\alpha X_\alpha\right),Y\right]_{S^{-1}C}
\longrightarrow \prod_\alpha [QX_\alpha,Y]_{S^{-1}C}$$
is an isomorphism.

First we prove that the map is epi. Suppose we are given diagrams
in $C$
$$\xymatrix{
&&Y_\alpha\\
X_\alpha\ar[urr]^{f_\alpha}&&Y\ar[u]_{s_\alpha}
}$$
representing elements of $[QX_\alpha,Y]_{S^{-1}C}$. By 3.10 there is a single map
${t:Y\to Z}$ in $S$ through which all the $s_\alpha$ factor. Then the elements in
$[QX_\alpha,Y]_{S^{-1}C}$ are equally well represented by diagrams
$$\xymatrix{
&&Z\\
X_\alpha\ar[urr]^{g_\alpha}&&Y\ar[u]_t
}$$
Then
$$\xymatrix{
&&Z\\
\bigvee_\alpha X_\alpha\ar[urr]^{\{g_\alpha\}}&&Y\ar[u]_t
}$$
represents an element in $\left[Q\left(\bigvee_\alpha X_\alpha\right),Y\right]_{S^{-1}C}$
which restricts to the given elements.

It remains to show that the function is mono. Suppose we are given
two elements in $\left[Q\left(\bigvee_\alpha X_\alpha\right),Y\right]_{S^{-1}C}$
with the same image in $\prod_\alpha [QX_\alpha,Y]_{S^{-1}C}$.
They may be represented as diagrams
$$\xymatrix{
&&Y_1\\
\bigvee_\alpha X_\alpha\ar[urr]^f&&Y\ar[u]_{s_1}
}\qquad\qquad\qquad
\xymatrix{
&&Y_2\\
\bigvee_\alpha X_\alpha\ar[urr]^g&&Y\ar[u]_{s_2}
}
$$
By the same argument as was used in the latter part of the proof of Proposition 3.7, we may assume 
$s_1=s_2$. Our data now says that for each $\alpha$, the components
$$\xymatrix{
X_\alpha\ar@<1ex>[r]^{f_\alpha}\ar@<-1ex>[r]_{g_\alpha} &Y
}$$
satisfy $(Qs_1)^{-1}(Qf_\alpha)=(Qs_1)^{-1}(Qg_\alpha)$, i.e. $Qf_\alpha=Qg_\alpha$.
Now Lemma 3.11 states that there is a single map $t:Y\to Z$ in $S$ such that
$tf_\alpha=tg_\alpha$ for each $\alpha$. Then we conclude that $tf=tg$ in $C$, and so
$Qf=Qg$. Hence
$$(Qs_1)^{-1}(Qf)=(Qs_1)^{-1}(Qg)$$
in $S^{-1}C$. This completes the proof of Lemma 3.12.\hfill$\Box$

\vspace{10pt}

The above results were derived for any category with arbitrary coproducts, with $S$
satisfying 3.1-3.6. Now we restrict our attention to the category $\C$ whose objects
are connected CW-complexes with basepoint and whose morphisms are homotopy classes.

\vspace{10pt}

\noindent
{\bf Lemma 3.13.} \textit{The diagram
$$\xymatrix{
U\ar@{^{(}->}[rr]^{i_1} &&U\cup V\\
U\cap V\ar@{^{(}->}[u]^{j_1}\ar@{^{(}->}[rr]^{j_2}
&&V\ar@{^{(}->}[u]_{i_2}
}$$
remains a weak pushout in $S^{-1}\C$.}

\vspace{10pt}

\noindent
{\bf Proof.} Let $Y$ be an object in $S^{-1}\C$ and let the diagrams
$$\xymatrix{
&&Y_1\\
U\ar[urr]^{f_1}&&Y\ar[u]_{s_1}
}\qquad\qquad\qquad
\xymatrix{
&&Y_1\\
V\ar[urr]^{f_2}&&Y\ar[u]_{s_2}
}
$$
with $s_1$, $s_2$ in $S$ represent morphisms in $S^{-1}\C$ such that
$$(Qs_1)^{-1}(Qf_1)(Qj_1)=(Qs_2)^{-1}(Qf_2)(Qj_2).$$
Now as noted previously we may assume $s_1=s_2$. Then $(Qf_1)(Qj_1)=(Qf_2)(Qj_2)$
Hence by 3.11 there is a $t:Y_1\to Z$ such that $tf_1j_1=tf_2j_2$. Since the
diagram is a weak pushout in $\C$ there is a map $f:U\cup V\to Z$ such that
$fi_1=tf_1$, $fi_2=tf_2$. Then
\begin{eqnarray*}
\left(\left[Q(ts_1)\right]^{-1}Qf\right)Qi_1
&=&\left[Q(ts_1)\right]^{-1}Q(fi_1)\\
&=&\left[Q(ts_1)\right]^{-1}Q(tf_1)\\
&=&(Qs_1)^{-1}Qf_1
\end{eqnarray*}
Similarly
$$\left(\left[Q(ts_1)\right]^{-1}Qf\right)Qi_2=(Qs_1)^{-1}Qf_2$$
\hfill$\Box$

\vspace{10pt}

\noindent
{\bf Lemma 3.14.} \textit{The functor $Q:\C\to S^{-1}\C$ has a right adjoint ${R:S^{-1}\C\to \C}$.}

\vspace{10pt}

\noindent
{\bf Proof.} Consider $[QX,Y]_{S^{-1}\C}$ where $X$ varies over $\C$ and $Y$
stays fixed in $S^{-1}\C$. We get a contravariant functor from $\C$ to sets.
I claim it satisfies the hypotheses of E. H. Brown's Representability Theorem.
In fact, I have just proved in Lemma 3.12 that it satisfies the Wedge Axiom,
and in Lemma 3.13 that it satisfies the Mayer-Vietoris
Axiom. We conclude that for fixed $Y$ in $S^{-1}\C$ there is an object $RY$ in $\C$
and an isomorphism
$$[QX,Y]_{S^{-1}\C}\longleftrightarrow[X,RY]_\C$$
natural for maps of $X$.

It is now standard that there is just one way to define $R$ on maps so that
$R$ becomes a functor and the above isomorphism becomes natural for maps of Y.
This proves 3.14.\hfill$\Box$

\vspace{10pt}

\noindent
{\bf Proof of Theorem 3.8.} We now define $E$ to be the composite
$$\C\stackrel{Q}{\longrightarrow}S^{-1}\C\stackrel{R}{\longrightarrow}\C.$$
We define $\eta$ as follows. Since there is an isomorphism
$$[QX,QX]_{S^{-1}\C}\longleftrightarrow[X,RQX]_\C,$$
let $\eta_X\in[X, RQX]$ be the element corresponding to $1_{QX}$. It is standard that $\eta$
is natural.

We have to check that $(E,\eta)$ satisfies 2.1 and 2.2, and that they yield (via 2.6)
the same class $S$ we started from.

Let us use the isomorphism
$$[QRY,Y]_{S^{-1}\C}\longleftrightarrow[RY,RY]_\C$$
to define $\xi_Y\in[QRY, Y]_{S^{-1}\C}$ as the map corresponding to $1_{RY}$.
Then it is standard that $\xi$ is natural and $\xi$, $\eta$ satisfy the following identities:\\
(3.15)\qquad\qquad\qquad\qquad\qquad $\xi_{QX}\cdot Q\eta_X = 1_{QX},$\\
(3.16)\qquad\qquad\qquad\qquad\qquad $R\xi_Y\cdot\eta_{RY}=1_{RY}.$\\
Now $Q\eta_X:QX\to QRQX$ is a natural transformation. By Lemma 2.12 there is one and only
one natural transformation
$$\zeta_Y:Y\longrightarrow QRY$$
natural for maps of $Y$ in $S^{-1}\C$ and such that
$$\zeta_{QX}=Q\eta_X.$$
Now (3.15) gives $\xi_{QX}\cdot\zeta_{QX}=1_{QX}$. Since $Q$ is a one-to-one correspondence
on objects, this is the same as
$$\xi_Y\cdot\zeta_Y=1_Y,\qquad Y\in S^{-1}\C.$$
By the naturality of $\zeta$ we have the following commutative diagram
$$\xymatrix{
QRY\ar[rr]^{\xi_Y}\ar[d]^{\zeta_{QRY}} &&Y\ar[d]^{\zeta_Y}\\
QRQRY\ar[rr]^{QR\xi_Y}&&QRY
}$$
But applying $Q$ to (3.16) we get
$$QR\xi_Y\cdot Q\eta_{RY}=1_{QRY}$$
whence
$$1_{QRY}=QR\xi_Y\cdot\zeta_{QRY}=\zeta_Y\cdot\xi_Y.$$
We thus conclude that $\xi$ and $\zeta$ are mutually inverse natural equivalences.

Applying $R$ to (3.15) we get
$$R\xi_{QX}\cdot\eta_{RQX}=1_{RQX};$$
substituting $Y=QX$ in (3.16) we get
$$R\xi_{QX}\cdot\eta_{RQX}=1_{RQX}.$$

Then $E\eta_X=RQ\eta_X$ and $\eta_{EX}=\eta_{RQX}$
 are both inverses to the natural equivalence $R\xi_{QX}$
so they must be equal and natural equivalences. This proves $(E,\eta)$ satisfy 2.1 and 2.2.

If $f\in S$, then $Qf$ is invertible so $Ef=RQf$ is invertible. Conversely, suppose
there is given $f:X\to Y$ in $\C$ such that $Ef=RQf$
is invertible. Then $QRQf$ is invertible. However, we have the following commutative diagram
$$\xymatrix{
QX\ar[d]^{Qf} &&QRX\ar[ll]_\simeq\ar[d]^{QRQf}_\simeq\\
QY &&QRQY\ar[ll]_\simeq
}$$
in which the horizontal arrows are given by $\xi$ or $\zeta$. Therefore $Qf$ is invertible.
Now Axiom 3.5 shows $f\in S$.

This proves that $S$ does arise from $(E,\eta)$, by 2.6. This completes the proof of Theorem 3.8.
\hfill$\Box$

\newpage

\section{A Further Axiom}\label{sec4}

In the applications, Axioms 3.1, 3.5 and 3.6 are very easy to verify, as we will see.
However, 3.2 and 3.3 (as well as 3.4) are less convenient because they are existence
statements and they leave us to construct
$$\xymatrix{
X\ar@{-->}[r]^g &{\color{gray}Z}\\
&Y\ar@{-->}[u]_t
}$$
for 3.2, or
$$\xymatrix{
X\ar@{-->}[r]^t &{\color{gray}Z}
}$$
for 3.3. It is much more convenient to be given a definite map and be told
to verify that it lies in $S$ (cf. 3.6). We therefore introduce the following axiom.
For it we assume that $U\cup V$ is a complex which is the union of subcomplexes $U$, $V$.

\vspace{10pt}

\noindent
{\bf Axiom 4.1.} If $i:U\cap V\longrightarrow U$ is in $S$, then 
$j:V\longrightarrow U\cup V$ is in $S$.

\vspace{10pt}

\noindent
{\bf Remarks.} This is a sort of excision axiom (cf. following page). To show
$j\in S$ it is sufficient to verify that
$$[U\cup V,Z]\stackrel{j^*}{\longrightarrow}[V,Z]$$
is a bijection for all $Z\in D$. We now show $j^*$ is automatically surjective.
Hence to check that $S$ satisfies 4.1, it suffices to verify that $j^*$ is injective.

For if $i\in S$ and $Z\in D$ the map
$$i^*:[U,Z]\longrightarrow[U\cap V,Z]$$
is a bijection. Hence given $g\in[V, Z]$ there is a map $f\in[U, Z]$ such that $gi'=fi$
where $i':U\cap V\longrightarrow U$. Since the following diagram is a weak pushout
$$\xymatrix{
&&&&Z\\
X\ar[urrrr]^g\ar[rr]_j &&U\cup V\ar@{-->}[urr]_{\exists h}\\
U\cap V\ar[u]_{i'}\ar[rr]^i &&U\ar[u]^{j'}\ar[uurr]_f
}$$
there is a map $h:U\cup V\longrightarrow Z$ such that $hj=g$. Hence $j^*$ is always surjective.

\vspace{10pt}

\noindent
{\bf Proposition 4.2.} \textit{If $S$ satisfies 3.1, 3.4, 3.5, 3.6 and 4.1, then $S$
satisfies 3.2 and 3.3.}

\vspace{10pt}

\noindent
{\bf Proposition 4.3.} \textit{There exist subclasses $S$ in $\C$ satisfying Axioms 3.1 to 3.6
inclusive but not 4.1.}

\vspace{10pt}

\noindent
{\bf Proof of 4.2.\,(i)}	We wish to prove Axiom 3.2. Suppose we are given a diagram
$$\xymatrix{
X\\
W\ar[u]^s\ar[r]^f &Y
}$$
with $s$ in $S$. Without loss of generality we may suppose $s$ and $f$ are injections of complexes. For we may represent $s$ and $f$ as cellular maps and then take the corresponding mapping cylinders $X'$, $Y'$.
The diagram
$$\xymatrix{
X\ar@{<->}[rr]^\simeq &&X'\\
&W\ar[ul]^s\ar@{^{(}->}[ur]
}$$
and Axiom 3.5 imply the inclusion $W\hookrightarrow X'$ also is in $S$.  (More generally,
3.5 allows us to change maps in $S$ by equivalences without leaving $S$.)

We now take $Z=X\cup_W Y$. Then we get a commutative diagram
$$\xymatrix{
X\ar@{^{(}->}[rr] &&Z\\
W\ar@{^{(}->}[u]^s\ar@{^{(}->}[rr]^f &&Y\ar@{^{(}->}[u]
}$$
Since $W\hookrightarrow X$ is in $S$, by Axiom 4.1 $Y\hookrightarrow Z$ is in $S$.

If you don't wish to replace $s$ and $f$ by inclusions, but merely to suppose
they are cellular maps, you can construct $Z$ as
$$X\cup_s\left(\frac{I\times W}{I\times\mbox{pt}}\right)\cup_f Y$$
and then of course you can divide this into two parts $U$ and $V$ at $\frac{1}{2}\times W$
and argue similarly.\hfill$\Box$

\vspace{10pt}

Before we go on with the proof of 4.2, we need the following lemma. We assume $S$
satisfies 3.1, 3.4, 3.5, 3.6, and 4.1.

\vspace{10pt}

\noindent
{\bf Lemma 4.4.} \textit{Suppose the inclusion $A\longrightarrow X$ is in $S$. Then the
inclusion
$$(0\times X)\cup\frac{I\times A}{I\times\mbox{\rm pt}}\cup(1\times X)\longrightarrow
\frac{I\times X}{I\times\mbox{\rm pt}}$$
is also in $S$.}

\vspace{10pt}

\noindent
{\bf Proof.} The inclusion
$$(0\times A)\vee(1\times A)\longrightarrow (0\times X)\vee(1\times X)$$
is in $S$ by Axiom 3.6. Let
\begin{eqnarray*}
U &= &(0\times X)\vee(1\times X),\\
V &= &\frac{I\times A}{I\times\mbox{pt}}.
\end{eqnarray*}
Then we have
\begin{eqnarray*}
U\cap V &= &(0\times A)\vee(1\times A)\\
U\cup V &= &(0\times X)\cup\frac{I\times A}{I\times\mbox{pt}}\cup(1\times X),
\end{eqnarray*}
so the inclusion $U\cap V\longrightarrow U$ is in $S$. By 4.1, the inclusion
$V\longrightarrow U\cup V$ is in $S$. Now since the inclusion $A\longrightarrow X$ is in $S$,
it follows from the diagram
$$\xymatrix{
A\ar@{<->}[rr]^\simeq\ar[d] &&\frac{I\times A}{I\times\mbox{\scriptsize pt}}=V\ar[d]\\
X\ar@{<->}[rr]^\simeq &&\frac{I\times X}{I\times\mbox{\scriptsize  pt}}
}$$
and 3.5 that the inclusion $V\longrightarrow\frac{I\times X}{I\times\mbox{\scriptsize pt}}$
\ is in $S$. Now in the following diagram
$$\xymatrix@R=+5pt{
V=\frac{I\times A}{I\times\mbox{\scriptsize pt}}\ar@{^{(}->}[rr]\ar@{^{(}->}[dr]
&&(0\times X)\cup\frac{I\times A}{I\times\mbox{\scriptsize pt}}\cup(1\times X)\ar@{_{(}->}[dl]\\
&\frac{I\times X}{I\times\mbox{\scriptsize pt}}
}$$
two of the three maps have been shown to be in $S$. By 3.5 it follows that the third is also in $S$.
\hfill$\Box$

\vspace{10pt}

\noindent
{\bf Proof of 4.2.\,(ii)} We wish to prove Axiom 3.3. Suppose we are given a diagram
$$\xymatrix{
W\ar[r]^s &X\ar@<1ex>[r]^{f_0}\ar@<-1ex>[r]_{f_1} &Y
}$$
of maps rather than homotopy classes with $f_0s\simeq f_1s$. Again we may assume $s$
is an inclusion of complexes. Now the homotopy $f_0s\simeq f_1s$ provides a map 
$h:(0\times X)\cup\frac{I\times W}{I\times\mbox{\scriptsize pt}}\cup(1\times X)
\longrightarrow Y$ and we get the diagram
$$\xymatrix{
\frac{I\times X}{I\times\mbox{\scriptsize pt}}\ar@{-->}[rr]
&&{\color{gray}Z}\\
(0\times X)\cup\frac{I\times A}{I\times\mbox{\scriptsize pt}}\cup(1\times X)
\ar[u]_i\ar[rr]^{\qquad\qquad h}
&&Y\ar@{-->}[u]_t
}$$
Now $i$ is in $S$ by 4.4. Hence by Axiom 3.2 which we have already proved,
we can fill in the above diagram. Then $tf_0\simeq tf_1$.\hfill$\Box$

\vspace{10pt}

\noindent
{\bf Proof of 4.3.} We need a counterexample. Let $ k$ be a field such that
the product map $ k\otimes_\mathbf{Z} k\longrightarrow k$
is an isomorphism; this holds precisely for the prime fields $ k=\mathbf{F}_p,\mathbf{Q}$.
Passing to direct sums, we see that for any vector space $V$ over $ k$
the product map $ k\otimes_\mathbf{Z}V\longrightarrow V$
is an isomorphism. So if $Y$ is an Eilenberg-Mac\,Lane space of type $(V,n)$ we have
$$H_n(Y; k)\cong V$$
as a $ k$-module. Now I recall that maps $f:X\to Y$ are in one-to-one
correspondence with $ k$-linear maps
$$H_n(X; k)\longrightarrow H_n(Y; k)$$
under $f\mapsto f_*$. So given $X$, choose $EX$ to be an Eilenberg-Mac\,Lane space of type $(V,n)$
such that
$$H_n(EX; k)\cong H_n(X; k)$$
and let the map $\eta_X:X\to EX$ realise the isomorphism. Clearly we
have an idempotent triple. A map $f:X\to Y$ is in $S$ if and only if 
$f_*:H_n(EX; k)\longrightarrow H_n(Y; k)$ is an isomorphism.

Now let $U$, $V$ denote the upper and lower hemispheres of the $n$-sphere, i.e.,
\begin{eqnarray*}
U &= &E^n_+\\
V &= &E^n_-\\
U\cap V &= &S^{n-1}\\
U\cup V &= &S^n
\end{eqnarray*}
Then $U\cap V\longrightarrow U$ lies in $S$, but $V\longrightarrow U\cup V$ does not lie in $S$.

\vspace{10pt}

\noindent
{\bf Example 4.5.} Let $K_*$ be a generalised homology theory satisfying the usual
axioms including
$$\colim_{\alpha}K_*(X_\alpha)\stackrel{\cong}{\longrightarrow}K_*(X)$$
where $X_\alpha$ runs over the finite subcomplexes of $X$. Define $S$ as follows:
a map $f:X\to Y$ is to lie in $S$ if and only if $f_*:K_*(X)\longrightarrow K_*(Y)$
is an isomorphism. I claim this satisfies all the axioms, except possibly 3.4.
We check 3.1. Suppose
$$X_1\stackrel{f_1}{\longrightarrow} X_2\stackrel{f_2}{\longrightarrow}X_3\longrightarrow\cdots
\longrightarrow X_n\stackrel{f_n}{\longrightarrow}X_{n+1}$$
are maps such that $(f_i)_*:K_*(X_i)\longrightarrow K_*(X_{i+1})$ are isomorphisms.
The composite
$$\xymatrix{
X_1\ar[rr]^{f_nf_{n-1}\cdots f_2f_1}&&X_{n+1}
}$$
has the same property.

We check 3.5. By construction, the functor $\C\stackrel{K_*}{\longrightarrow}\mbox{(graded groups)}$
has the property that $f\in S$ implies $K_*(f)$ is invertible. Therefore $K_*$ factors through
$S^{-1}\C$:
$$\xymatrix{
\C\ar[rr]^{K_*}\ar[dr]_Q &&\mbox{graded groups}\\
&S^{-1}\C\ar[ur]_T
}$$
Now suppose $f:X\to Y$ in $\C$ is such that $Qf$ is invertible. Then $TQf$ is
invertible, i.e., $K_*(f)$ is invertible so $f\in S$.

We check 3.6. Suppose $f_\alpha:X_\alpha\to Y_\alpha$ is in $S$ for each $\alpha$, i.e.
$${(f_\alpha)_*:K_*(X_\alpha)\to K_*(Y_\alpha)}$$
is an isomorphism for all $\alpha$.
Then we obtain the diagram
$$\xymatrix{
\widetilde{K}_*\left(\bigvee_\alpha X_\alpha\right)
\ar[rr]^{\left(\bigvee_\alpha f_\alpha\right)_*}
&&\widetilde{K}_*\left(\bigvee_\alpha Y_\alpha\right)\\
\bigoplus_\alpha\widetilde{K}_*\left(X_\alpha\right)
\ar[u]_\cong\ar[rr]^{\oplus_\alpha(f_\alpha)_*}_\cong
&&\bigoplus_\alpha\widetilde{K}_*\left(Y_\alpha\right)\ar[u]_\cong
}$$
whence $(\bigvee_\alpha f_\alpha)_*$ is an isomorphism and $\bigvee_\alpha f_\alpha$ lies in $S$.

We check 4.1. Suppose $i:U\cap V\longrightarrow U$ lies in $S$, i.e.,
$${i_*:K_*(U\cap V)\longrightarrow K_*(U)}$$
is an isomorphism. By exactness this is the same
as saying $K_*(U,U\cap V)=0$. By excision this is the same as saying $K_*(U\cup V, V)= 0$.
Then by exactness again
$$j_*: K_*(V)\longrightarrow K_*(U\cup V)$$
is an isomorphism, so $j\in S$.\hfill$\Box$

\vspace{10pt}

{\color{gray}
\noindent
{\bf Remark.} Unfortunately there seems to be no way in the above example to verify
that $S$ satisfies the set-theoretic condition expressed in Axiom 3.4. Until we get around
to verifying 3.4 under suitable extra assumptions or proving the existence of the desired
idempotent functors $(E,\eta)$ in some other way, the following two conjectures remain
pious hopes or indications of what we wish to prove.
}

\vspace{10pt}

\noindent
{\bf Editorial Note.} The pessimism expressed in the above remark proved to be unwarranted.
It turns out that there is a simple alternative to Axiom 3.4 which allows us to prove
the existence of idempotent functors $(E,\eta)$ for these classes $S$. See the epilogue
for details. The editor has thus taken the liberty of upgrading the following two statements
in the original manuscript from conjectures to theorem/corollary and made a few other revisions
in some subsequent statements referring to the conjectural status of these results.

\vspace{10pt}

\newpage

\noindent
{\bf Theorem 4.6.} \textit{Suppose there is given a generalised homology
theory $K_*$. Then for each $X\in \C$, there exists an object $EX$ and a map $\eta_X:X\to EX$
such that
\begin{itemize}
\item[(i)] $(\eta_X)_*:K_*(X)\longrightarrow K_*(EX)$ is an isomorphism,
\item[(ii)] $\eta_X$ is couniversal with respect to (i).
\end{itemize}}

\vspace{10pt}

This follows from the slight revision of the main theorem 3.8, discussed in the epilogue, and 
from 2.11. In this case the category of fractions $S^{-1}\C$ gives you as much of homotopy theory as you can see through the eyes of ${K_*\mbox{-theory}}$. Up to equivalence, it would be embedded in $\C$ as the full subcategory $D$.

\vspace{10pt}

\noindent
{\bf Corollary 4.7.} \textit{Take $K_*(-)$ to be $H_*(- ; A)$ where $A$ is a subring of $\mathbf{Q}$,
i.e., $A$ is obtained by localising $\mathbf{Z}$. In this way we get a localisation functor defined without restrictions on $\pi_*(X)$.}

\vspace{10pt}

For the remainder of this section, let $(E,\eta)$ and $S$ be as in Corollary 4.7, i.e., ${f:X\to Y}$ is in $S$ if
and only if $f_*: H_*(X; A)\longrightarrow H_*(Y; A)$ is an isomorphism. Let $D$ be the
corresponding subclass of $\C$, i.e., $Z\in D$ if and only if
$$f^*:[Y, Z]\longrightarrow [X, Z]$$
is an isomorphism for all $f\in S$.

\vspace{10pt}

\noindent
{\bf Lemma 4.8.} \textit{Let $Z\in D$ and let $\Gamma$ be a subgroup of $\pi_1(Z)$ such that 
$(\Gamma/[\Gamma,\Gamma])\otimes A=0.$ Then $\Gamma=0$.}

\vspace{10pt}

\noindent

\noindent
{\bf Proof.} Let $\{\gamma_\alpha\}$ be a set of generators for $\Gamma$. Since 
$[\gamma_\alpha]\otimes 1$ is zero in $(\Gamma/[\Gamma,\Gamma])\otimes A$, there
exists an integer $n$ which is invertible in $A$ such that
$$\gamma_\alpha^n=1\mbox{ in }\Gamma/[\Gamma,\Gamma].$$
Then $\gamma_\alpha^n$ is a product of commutators, say
$$\gamma_\alpha^n=\prod_\beta[\delta_{\alpha\beta},\epsilon_{\alpha\beta}],
\ \delta_{\alpha\beta}, \epsilon_{\alpha\beta}\in\Gamma.$$
Define a map $\bigvee_\alpha S^1\longrightarrow Z$ such that
$$\iota_\alpha\mapsto \gamma_\alpha,$$
where $\iota_\alpha$ is a generator of the $\alpha$-th free summand of 
$\pi_1\left(\bigvee_\alpha S^1\right)$. Then
$${\pi_1\left(\bigvee_\alpha S^1\right)\longrightarrow\Gamma\subseteq\pi_1(Z)}$$
is epi.
Choose elements, $\overline{\delta}_{\alpha\beta}$, $\overline{\epsilon}_{\alpha\beta}$
in $\pi_1\left(\bigvee_\alpha S^1\right)$ mapping onto $\delta_{\alpha\beta}$,
$\epsilon_{\alpha\beta}$. Form a space $X$ by attaching to $\bigvee_\alpha S^1$
a 2-cell $e_\alpha^2$, one for each $\alpha$, by a map in the class
$$\iota_\alpha^{-n}\prod_\beta[\overline{\delta}_{\alpha\beta},\overline{\epsilon}_{\alpha\beta}].$$
Then the map $\bigvee_\alpha S^1\longrightarrow Z$ extends to give a map $X\to Z$ and the attaching
maps induce isomorphisms on $H_*(-;A)$. Computing with the cellular chain complex of $X$ we obtain
$\widetilde{H}_*(X;A)=0$, and so the constant map $X\to\mbox{pt}$ lies in $S$. Since $Z\in D$, the
map $X\to Z$ factors to give
$$\xymatrix{
X\ar[rr]\ar[dr] &&Z\\
&\mbox{pt}\ar[ur]
}$$
Therefore the map $\pi_1(X)\to\pi_1(Z)$ is zero. But it mapped onto $\Gamma$,
so $\Gamma=0$.\\
${}_{}\qquad$\hfill$\Box$

\vspace{10pt}

\noindent
{\bf Corollary 4.9.} \textit{If $\pi_1(X)=0$, then $\pi_1(EX)=0$.}

\vspace{10pt}

\noindent
{\bf Proof.} If $\pi_1(X)=0$, then $H_1(X; A)=0$; therefore $H_1(EX;A)=0$
since $\eta_X\in S$. Now apply Lemma 4.8 with $Z=EX$, $\Gamma=\pi_1(Z)$.
\hfill$\Box$

\vspace{10pt}

Next we show that the homotopy groups of $EX$ are always local.

\vspace{10pt}

\noindent
{\bf Proposition 4.10.} \textit{If $Z\in D$ and $n$ is invertible in $A$, then division by $n$ is possible and unique in each homotopy group $\pi_i(Z)$ (including $i=1$).}

\vspace{10pt}

\noindent
{\bf Proof.} Let $n:S^i\to S^i$ be the map which induces multiplication by $n$ on the homotopy group
$\pi_i(S^i)$. Then $n_*:\widetilde{H}_*(S^i;A)\longrightarrow\widetilde{H}_*(S^i;A)$ also is multiplication by $n$ and
hence is an isomorphism, so $n\in S$. Hence
$$\pi_i(Z)=[S^i,Z]\longrightarrow[S^i,Z]=\pi_i(Z),$$
which is again multiplication by $n$, is an isomorphism.\hfill$\Box$

\vspace{10pt}

\noindent
{\bf Proposition 4.11.} \textit{Suppose $\pi_1(Y)=0$. Then there is a map $s:Y\to Y^{(\infty)}$ in $S$
with $\pi_1\left(Y^{(\infty)}\right)=0$ and with all $\pi_i\left(Y^{(\infty)}\right)$ being $A$-modules.}

\vspace{10pt}

\noindent
{\bf Proof.} We define inductively a sequence of spaces $Y^{(i)}$ and maps
$${s_i:Y^{(i)}\longrightarrow Y^{(i+1)}}$$ 
in $S$ such that $\pi_j\left(Y^{(i)}\right)$
is an $A$-module and 
$(s_i)_\#:\pi_j\left(Y^{(i)}\right)\longrightarrow\pi_j\left(Y^{(i+1)}\right)$
is an isomorphism for $j\le i$. We choose $Y^{(1)}=Y$.

Having chosen $Y^{(i)}$, let $W$ be a Moore space of type $\pi_{i+1}\left(Y^{(i)}\right)$
in dimension $i+1$, and let $f:W\longrightarrow Y^{(i)}$ be a map inducing an isomorphism
$\pi_{i+1}(W)\stackrel{\cong}{\longrightarrow}\pi_{i+1}\left(Y^{(i)}\right)$.
Now embed $W$ in a Moore space $X$ of type $\pi_{i+1}\left(Y^{(i)}\right)\otimes A$
in dimension $i+1$, so that the map $e:W\to X$ is in $S$. Then we have a
commutative diagram
$$\xymatrix{
W\ar[rr]^f\ar[d]^e &&Y^{(i)}\ar[d]^j\\
X\ar[rr]^k &&X\cup\frac{I\times W}{I\times\mbox{\scriptsize pt}}\cup Y^{(i)}
}$$
with $e,j\in S$. If we let
$Y^{(i+1)}=X\cup\frac{I\times W}{I\times\mbox{\scriptsize pt}}\cup Y^{(i)}$
and $s_i=j$, we see
that ${s_i:Y^{(i)}\longrightarrow Y^{(i+1)}}$ and $Y^{(i+1)}$
have all the required properties.

Having defined the $Y^{(i)}$'s, define $Y^{(\infty)}$ to be the colimit
of the $Y^{(i)}$'s. Clearly the induced map  $s:Y=Y^{(1)}\to Y^{(\infty)}$ is in $S$,
 $\pi_1\left(Y^{(\infty)}\right)=0$ and $\pi_i\left(Y^{(\infty)}\right)$ is
an $A$-module for all $i$.\hfill$\Box$

\vspace{10pt}

We note that the construction of Proposition 4.11 is precisely Sullivan's cellular construction of his
localisation functor for simply connected spaces. We now show that it is equivalent to our
localisation functor.

\vspace{10pt}

\noindent
{\bf Proposition 4.12.} \textit{If $\pi_1(Z)$ acts trivially on $\pi_i(Z)$ (which implies $\pi_1(Z)$
is abelian) and each $\pi_i(Z)$ is an A-module, then $Z\in D$.}

\vspace{10pt}

\noindent
{\bf Proof.} By 2.10 it suffices to show that whenever $f:X\to Y$ is in $S$, then
$$f^*: [Y,Z]\longrightarrow[X,Z]$$
is epi. We may assume that $f$ is an inclusion.

Now let $g:X\to Z$ be a map. Since $f$ is in $H_*(Y,X;A)=0$.
Hence $H^{i+1}(Y,X; \pi_i(Z))= 0)$. Therefore by obstruction theory, $g$ extends
to a map $\tilde{g}:Y\to Z$, i.e., $\tilde{g}f=g$. Hence $f^*$ is epi and $Z\in D$.
\hfill$\Box$

\vspace{10pt}

\noindent
{\bf Proposition 4.13.} \textit{For simply connected spaces $Y$ the localisation functor ${\eta_Y:Y\to EY}$
of Corollary 4.7 is equivalent to the construction of Proposition 4.11 (i.e. Sullivan's
localisation). Moreover $(\eta_Y)_\#:\pi_i(Y)\to \pi_i(EY)$ is a localising map for all $i$.}

\vspace{10pt}

\noindent
{\bf Proof.} By 4.11 $s:Y\to Y^{(\infty)}$ is in $S$, $\pi_1\left(Y^{(\infty)}\right)=0$ and 
$\pi_i\left(Y^{(\infty)}\right)$ are $A$-modules for all $i$. By Proposition 4.12 it follows
that $Y^{(\infty)}\in D$. By Proposition 2.11 $EY\simeq Y^{(\infty)}$ and we have a 
commutative diagram
$$\xymatrix{
&&Y^{(\infty)}\ar@{<->}[dd]^\simeq\\
Y\ar[urr]^s\ar[drr]_{\eta_Y}\\
&&EY}$$
Thus $EY$ is simply connected and the homotopy groups of $EY$ are $A$-local.

By definition of $S$
$$(\eta_Y)_*:H_*(Y;A)\longrightarrow H_*(EY;A)$$
is an isomorphism. But since $Y$ and $EY$ are simply connected,
then by Serre's $C$-theory this is equivalent to saying that
$$(\eta_Y)_\#\otimes 1:\pi_i(Y)\otimes A\longrightarrow \pi_i(EY)\otimes A$$
is an isomorphism for each $i$ (Spanier \cite[p. 512, Thm. 22]{Sp}).

Now consider the following diagram
$$\xymatrix{
\pi_i(Y)\ar[rr]^{(\eta_Y)_\#}\ar[d] &&\pi_i(EY)\ar[d]^\cong\\
\pi_i(Y)\otimes A\ar[rr]^{(\eta_Y)_\#\otimes 1}_\cong &&\pi_i(EY)\otimes A
}$$
We have just shown the bottom arrow is an isomorphism. The right-hand arrow is an
isomorphism since the homotopy groups of $EY$ are $A$-modules. The left-hand arrow
is a localising map. Therefore the top arrow is also a localising map.\hfill$\Box$

\newpage

\section{Behaviour of Idempotent Functors with Respect to Fiberings;
Construction of Localisation Using Postnikov Decomposition}\label{sec5}

Assume $\C$ is the homotopy category of connected CW-complexes, and $S$ is a class of morphisms which satisfies 3.1, 3.5, 3.6 and 4.1 but not necessarily 3.4. We define $D$ to be the class of all spaces $Z$ such that whenever $f:X\to Y$ is in $S$ then
$$f^*:[Y, Z]\longrightarrow[X, Z]$$
is an isomorphism.

If we knew that $S$ defined an idempotent functor $(E,\eta)$, then by 2.10 in order to show that $Z\in D$ it would suffice to check that $f^*$
is epi for all $f\in S$. Since we do not require $S$ to satisfy 3.4, we need the following lemma.

\vspace{10pt}

\noindent
{\bf Lemma 5.1.} \textit{A space $Z$ is in $D$ if and only if
$$f^*:[Y, Z]\longrightarrow[X, Z]$$
is epi for all $f\in S$.}

\vspace{10pt}

\noindent
{\bf Proof.} First note that if $i:X\to Y$ is an inclusion and if ${i^*:[Y, Z]\longrightarrow[X,Z]}$ is epi, then each map $g:X\to Z$ can be extended over $Y$. For if this is the case, some homotopic map $g'$ extends to $h':Y\to Z$ and we can use the homotopy extension property to extend $g$.

Now assume $f^*:[Y, Z]\longrightarrow[X, Z]$ is epi for all $f\in S$. Assume
$f:X\to Y$ is in $S$. We wish to show $f^*:[Y, Z]\longrightarrow[X, Z]$ is an isomorphism. Without loss of generality we may suppose $f$ is an inclusion.

We already know that $f^*$ is epi, so we only have to check that it is also mono. Suppose we are given 
$g_0, g_1: Y\to Z$ and a homotopy $g_0|X\simeq g_1|X$.
This defines a map $h:0\times Y\cup\frac{I\times X}{I\times\mbox{\scriptsize pt}}\cup 1\times Y\longrightarrow Z$.
Now by Lemma 4.4, the inclusion
$$0\times Y\cup\frac{I\times X}{I\times\mbox{pt}}\cup 1\times Y\longrightarrow 
\frac{I\times Y}{I\times\mbox{pt}}$$
lies in $S$. So by our assumption $h$ extends to a map $H$ over 
$\frac{I\times Y}{I\times\mbox{\scriptsize pt}}$
$$\xymatrix{
0\times Y\cup\frac{I\times X}{I\times\mbox{\scriptsize pt}}\cup 1\times Y
\ar[drr]^h\ar@{^{(}->}[dd]\\
&&Z\\
\frac{I\times Y}{I\times\mbox{\scriptsize pt}}\ar@{-->}[urr]^H
}$$
which gives a homotopy $g_0\simeq g_1$. Hence $f^*$ is mono and therefore an isomorphism.
\hfill$\Box$

\vspace{10pt}

\noindent
{\bf Proposition 5.2.} \textit{Suppose we are given a diagram
$$\xymatrix{
&&V\ar[d]^q\\
U\ar[rr]^p &&W
}$$
and form the homotopy pullback
$$\xymatrix{
F\ar[rr]\ar[d]&&V\ar[d]^q\\
U\ar[rr]^p &&W
}$$
If $U$, $V$, and $W$ are in $D$, then $F$ is (weakly) in $D$.}

\vspace{10pt}

\noindent
{\bf Note.} By the standard construction, a point in $F$ is a triple $(u,v,\omega)$
where $u\in U$, $v\in V$, and $\omega:I\to W$ is a path from $\omega(0)=p(u)$ to
$\omega(1)=q(v)$. So on the face of it, $F$ is not a CW-complex; we should simply
assert that any weakly equivalent CW-complex is in $D$.

\vspace{10pt}

\noindent
{\bf Proof of 5.2.} We want to prove  $f^*:[Y, F]\longrightarrow[X, F]$ is epi 
for all $f:X\to Y$ in $S$. Without loss of generality we can suppose $f$ is an inclusion.

Suppose we are given a map $X\to F$, that is, a pair of maps
$$g':X\longrightarrow U,\quad g'': X\longrightarrow V$$
and a specific homotopy
$$h:\frac{I\times X}{I\times\mbox{pt}}\longrightarrow W$$
between $pg'$ and $qg''$. Since $X\hookrightarrow Y$ 
is an inclusion in $S$ and
$U$, $V$ are in $D$, we can extend $g'$, $g''$ to $\tilde{g}'$, $\tilde{g}''$ over $Y$.
Now we have a map
$$\tilde{h}:0\times Y\cup\frac{I\times X}{I\times\mbox{pt}}\cup 1\times Y
\longrightarrow W$$
defined by $p\tilde{g}'$, $q\tilde{g}''$ and $h$. By 4.4, the inclusion
$$0\times Y\cup\frac{I\times X}{I\times\mbox{pt}}\cup 1\times Y\hookrightarrow
\frac{I\times Y}{I\times\mbox{pt}}$$
lies in $S$. Since $W$ is in $D$, we can extend $h$ to a map $H$ over
$\frac{I\times Y}{I\times\mbox{\scriptsize pt}}$.
Now the maps $\tilde{g}'$, $\tilde{g}''$ and $H$ give a map $Y\to F$ extending the given map
$X\to F$. This proves the proposition.\hfill$\Box$

\vspace{20pt}

\noindent
{\bf Remark.} The study of homotopy pullbacks includes that of fiberings as the special case:
$$\xymatrix{
F\ar[rr]\ar[d] &&E\ar[d]\\
\mbox{pt}\ar[rr] &&B
}$$
The more general study is desirable, e.g., for the study of Postnikov
decompositions in the non-simply connected case.

\vspace{20pt}

Now we introduce another axiom.

\newpage

\noindent
{\bf Axiom 5.3.} Suppose we are given a strictly commutative diagram of maps
$$\xymatrix{
&V\ar[rr]^v\ar'[d]^q[dd] &&V'\ar[dd]^{q'}\\
F\ar@{-->}[dd]\ar@{-->}[ur]\ar@{-->}[rr]^{\quad f} &{} &F'\ar@{-->}[dd]\ar@{-->}[ur]\\
&W\ar'[r]_w[rr] &{} &W'\\
U\ar[ur]^p\ar[rr]_u &&U'\ar[ur]_{p'}
}$$
and form the homotopy pullback. If $u$, $v$, $w$ are in $S$, then $f$ is in $S$.

\vspace{10pt}

\noindent
{\bf Corollary 5.4.} \textit{If $E$ exists and $S$ satisfies 4.1 and 5.3, then $E$ preserves homotopy
pullbacks. More precisely, if $EU$, $EV$ and $EW$ exist, then $EF$ exists and is their homotopy pullback.}

\vspace{10pt}

\noindent
{\bf Proof.}  Take a diagram
$$\xymatrix{
&V\ar[rr]^{\eta_V}\ar'[d]^q[dd] &&EV\ar[dd]^{Eq}\\
F\ar@{-->}[dd]\ar@{-->}[ur]\ar@{-->}[rr]^{\quad f} &{} &F'\ar@{-->}[dd]\ar@{-->}[ur]\\
&W\ar'[r]_{\eta_W}[rr] &{} &EW\\
U\ar[ur]^p\ar[rr]_{\eta_U} &&EU\ar[ur]_{Ep}
}$$
If you wish we can suppose the diagram is strictly commutative (replace $EU$, $EV$ by mapping cylinders). Complete the diagram of homotopy pullbacks. Then $f$ is in $S$ by Axiom 5.3 and $F'$ is in $D$
by Proposition 5.2, so up to equivalence $f$ is the map $\eta_F$.\hfill$\Box$

\vspace{10pt}

\noindent
{\bf Proposition 5.5.} \textit{The localisation functor of Corollary 4.7 preserves homotopy pullbacks with 1-connected base $W$. More precisely, if $\pi_1(W)=0$, then
$$\xymatrix{
EF\ar[rr]\ar[d] &&EV\ar[d]\\
EU\ar[rr] &&EW
}$$
is a homotopy pullback.}

\vspace{10pt}

\noindent
{\bf Proof.} If $W$ is 1-connected, then $EW$ is 1-connected by 4.9. So in order to apply the
argument of Corollary 5.4, it is sufficient to check Axiom 5.3 with $W$ and $W'$ 1-connected.

Consider first the diagram
$$\xymatrix{
G\ar@{-->}[rr]^g\ar@{-->}[d] &&G'\ar@{-->}[d]\\
V\ar[rr]^v\ar[d]^q &&V'\ar[d]^{q'}\\
W\ar[rr]^w &&W'
}$$
Let us supply homotopy fibres $G$, $G'$. In the Serre spectral sequence the operations
of $\pi_1(W)$, $\pi_1(W')$ on $H_*(G;A)$, $H_*(G';A)$ are trivial, because
$\pi_1(W)=0$, $\pi_1(W')=0$. We are given
$$v_*:H_*(V;A)\longrightarrow H_*(V';A),\quad w_*:H_*(W;A)\longrightarrow H_*(W';A)$$
are isomorphisms. Therefore the comparison theorem for spectral
sequences allows one to prove that
$$g_*:H_*(G;A)\longrightarrow H_*(G';A)$$
is an isomorphism. But now consider
$$\xymatrix{
G\ar@{-->}[rr]^g\ar@{-->}[d] &&G'\ar@{-->}[d]\\
F\ar[rr]^f\ar[d] &&F'\ar[d]\\
U\ar[rr]^u &&U'
}$$
Let us supply the homotopy fibres; they are $G$, $G'$ up to equivalence carrying the map we want into $g$. In the Serre spectral sequence the actions of $\pi_1(U)$, $\pi_1(U')$ on $H_*(G;A)$, $H_*(G';A)$ are trivial, since they factor through $\pi_1(W)$, $\pi_1(W')$. Now the Serre spectral sequence
shows that
$$f_*: H_*(F; A)\longrightarrow H_*(F'; A)$$
is an isomorphism. This concludes the proof.\hfill$\Box$

\vspace{10pt}

\noindent
{\bf Proposition 5.5.} \textit{Let $\pi$ be an abelian group. Then the localisation
 of $K(\pi ,n)$ is $K(\pi\otimes A ,n)$. More precisely, the following diagram is commutative
$$\xymatrix{
\tilde{H}_*\left(K(\pi ,n)\right)\ar[rr]\ar[d] &&\tilde{H}_*\left(K(\pi\otimes A ,n)\right)\ar[d]^\cong\\
\tilde{H}_*\left(K(\pi ,n);A\right)\ar[rr]^\cong&&\tilde{H}_*\left(K(\pi\otimes A ,n);A\right)
}$$}

\vspace{10pt}

\noindent
{\bf Proof.} First we tackle the case $n=1$. Suppose $\pi=\mathbf{Z}$. Then we get 0 except in dimension 1, and there
$$\xymatrix{
\mathbf{Z}\ar[rr]\ar[d] &&A\ar[d]\\
A\ar[rr] &&A\otimes_\mathbf{Z} A
}$$
with the obvious maps.

Suppose $\pi=\mathbf{Z}_{p^k}$, $p$ invertible in $A$. Then we have 0 in
even dimensions, while in odd ones
$$\xymatrix{
\mathbf{Z}_{p^k}\ar[rr]\ar[d] &&0\ar[d]\\
0\ar[rr] &&0
}$$

Suppose $\pi=\mathbf{Z}_{p^k}$, $p$ not invertible in $A$. Then we have 0 in
even dimensions, while in odd ones
$$\xymatrix{
\mathbf{Z}_{p^k}\ar[rr]\ar[d] &&\mathbf{Z}_{p^k}\ar[d]\\
\mathbf{Z}_{p^k}\ar[rr] &&\mathbf{Z}_{p^k}
}$$
all the maps being isomorphisms.

Suppose $\pi\cong\pi'\oplus\pi''$ and the result is true for $\pi'$, $\pi''$.
The map
$$K(\pi' ,1)\times K(\pi'' ,1)\longrightarrow K(\pi'\otimes A ,1)\times K(\pi''\otimes A ,1)$$
 induces an isomorphism of $\tilde{H}_*(-; A)$ by the K\"unneth theorem. Thus we have
$$\xymatrix{
\tilde{H}_*\left(K(\pi' ,1)\times K(\pi'' ,1)\right)\ar[rr]\ar[d]
&&\tilde{H}_*\left(K(\pi'\otimes A,1)\times K(\pi''\otimes A,1)\right)\ar[d]\\
\tilde{H}_*\left(K(\pi' ,1)\times K(\pi'' ,1);A\right)\ar[rr]^\cong
&&\tilde{H}_*\left(K(\pi'\otimes A,1)\times K(\pi''\otimes A,1);A\right)
}$$
For the right-hand vertical arrow, suppose $m$ is invertible in $A$. Then
it acts by an isomorphism on $\tilde{H}_i\left(K(\pi'\otimes A,1)\right)$ and
$\tilde{H}_j\left(K(\pi''\otimes A,1)\right)$
and therefore on the tensor and torsion products of these; so by the K\"unneth formula
$\tilde{H}_*\left(K(\pi' ,1)\times K(\pi'' ,1);A\right)$ is a module over $A$, and the right-hand vertical arrow is an isomorphism. So the result is true for $\pi\cong\pi'\oplus\pi''$.

At this point we have proved the result where $\pi$ is finitely-generated and $n=1$
by the structure theorem for finitely-generated abelian groups.

The result for general $\pi$ and $n=1$ follows by passing to colimits.

Now we prove the general result by induction over $n$. Assume it for $n$. Regard $K(\pi,n)$
as the space of loops on $K(\pi,n+1)$ and similarly with $\pi$ replaced by $\pi\otimes A$.
If $K(\pi,n)\longrightarrow K(\pi\otimes A,n)$ induces an isomorphism of $\tilde{H}_*(-;A)$,
then so does $K(\pi,n+1)\longrightarrow K(\pi\otimes A,n+1)$ by the Rothenberg-Steenrod
spectral sequence. Similarly if
$\tilde{H}_*\left(K(\pi\otimes A,n)\right)$ is a module over $A$, then so is
$\tilde{H}_*\left(K(\pi\otimes A,n+1)\right)$ by the Rothenberg-Steenrod spectral sequence.
\hfill$\Box$

\vspace{10pt}

\noindent
{\bf Theorem 5.7.} \textit{Let $(E,\eta)$ be the localisation of Corollary 4.7. Let $X$ be a space such that $\pi_1(X)$ is nilpotent and acts trivially on the higher homotopy groups. Then
\begin{itemize}
\item[(i)] the homology groups $\tilde{H}_*(EX)$ are $A$-modules (then of course\\ ${(\eta_X)_*:\tilde{H}_*(X)\longrightarrow\tilde{H}_*(EX)}$ localises).
\item[(ii)] $(\eta_X)_\#:\pi_i(X)\longrightarrow\pi_i(EX)$ is a localising map.	(If $i=1$ this means that the induced map on the subquotients of the lower central series is a localising map.)
\end{itemize}
}

\vspace{10pt}

\noindent
{\bf Proof.} First suppose that $X$ has all but a finite number of its homotopy groups zero. Then $X$ can be built up by a finite sequence of fiberings with 1-connected base which is of type $K(\pi,n+1)$
$$\xymatrix{
X_n\ar@{-->}[rr]^{\eta_{X_n}}\ar[d] &&F=EX_n\ar@{-->}[d]\\
X_{n-1}\ar[rr]^{\eta_{X_{n-1}}}\ar[d]^p &&EX_{n-1}\ar[d]^{p'}\\
B=K(\pi,n+1)\ar[rr] &&K(\pi\otimes A,n+1)
}$$
Now we assume as an induction hypothesis that (i) and (ii) hold for\\
${\eta_{X_{n-1}}:X_{n-1}\longrightarrow EX_{n-1}}$.
By 5.5, $EB=K(\pi\otimes A,n+1)$ and $p'=Ep$. By Corollary 5.4, $EX_n$ is
(up to equivalence) the homotopy fibre $F$ of $p'$.
Now the homology groups $\tilde{H}_*(EX_n)$ are
$A$-modules by a spectral sequence argument, and
$(\eta_{X_n})_\#:\pi_i(X_n)\longrightarrow\pi_i(EX_n)$ is a localising map.

It remains to pass to the general case.

Suppose we are given a diagram
$$\xymatrix{
Z_1 &Z_2\ar[l]^{p_1} &Z_3\ar[l]^{p_2} &Z_4\ar[l]^{p_3} &\ \dots\ar[l]^{p_4}
}$$
I define the homotopy inverse-limit $Z_\infty$ as follows: an element of $Z_\infty$
is a sequence of functions $\omega_n:[n,n+1]\longrightarrow Z_n$ such that
$$\omega_n(n+1)=p_n\omega_{n+1}(n).$$
To conclude the proof, I first prove:

\vspace{10pt}

\noindent
{\bf Lemma 5.8.} \textit{If all the objects $Z_n$ are in $D$, then $Z_\infty$ is in $D$.}

\vspace{10pt}

\noindent
{\bf Proof.} Suppose given $f:X\to Y$ in $S$. Without loss of generality we can assume $f$ is an inclusion. Suppose we are given a map $g:X\to Z_\infty$;	equivalently, it is a sequence of maps
$$g_n:\frac{[n,n+1]\times X}{[n,n+1]\times\mbox{pt}}\longrightarrow Z_n.$$
Since $X\hookrightarrow Y$ is in $S$, I can extend $g_n$ from $n\times X$ to
$n\times Y$. Of course,
the extension of $g_{n+1}$ forces the extension $p_ng_{n+1}$ of $g_n$ on 
$(n+1)\times Y$.  Now since the inclusion
$$n\times Y\cup\frac{[n,n+1]\times X}{[n,n+1]\times\mbox{pt}}\cup
(n+1)\times Y\longrightarrow
\frac{[n,n+1]\times Y}{[n,n+1]\times\mbox{pt}}$$
is in $S$ by 4.1, we can extend $g_n$ over 
$\frac{[n,n+1]\times Y}{[n,n+1]\times\mbox{\scriptsize pt}}$. This extends $g$
to a map $Y\to Z$. This concludes the proof of the lemma.\hfill$\Box$

\vspace{10pt}

\noindent
{\bf Conclusion of the Proof of 5.7.} To apply the lemma, consider the diagram
$$\xymatrix{
X\ar[d]\ar[rr] &&Z_\infty\ar[d]\\
\vdots\ar[d] &&\vdots\ar[d]\\
X_n\ar[rr]\ar[d] &&EX_n\ar[d]\\
X_{n-1}\ar[rr]\ar[d] &&EX_{n-1}\ar[d]\\
\vdots\ar[d] &&\vdots\ar[d]\\
X_1\ar[rr] &&EX_1
}$$
where the left-hand column is the Postnikov system of $X$. Take the
homotopy inverse limit of the right-hand side and call it $Z_\infty$. Since
the diagram is at least homotopy commutative we have a map
$X\to Z_\infty$	(not necessarily unique at this stage).

Now $Z_\infty$ lies in $D$ by 5.8. Since the $\pi_i(EX )$ are at least
groups, it makes sense to observe that there is an exact sequence
$$0\longrightarrow\lim_{\stackrel{\longleftarrow}{n}}^{\quad\ 1}\pi_{i+1}(EX_n)
\longrightarrow\pi_i(Z_\infty)\longrightarrow
\lim_{\stackrel{\longleftarrow}{n}}\pi_i(EX_n)\longrightarrow 0$$
But by construction $\pi_i(EX_n)$ is constant for $n\gg i$; so we get
${\pi_i(Z_\infty)\stackrel{\cong}{\longrightarrow}\pi_i(EX_n)}$ for $n\gg i$.
The homotopy diagram shows $\pi_i(X)\longrightarrow\pi_i(Z_\infty)$ localises.
Also for $n\gg i$ we get a diagram
$$\xymatrix{
H_i(X;A)\ar[rr]\ar[d]^\cong &&H_i(Z_\infty;A)\ar[d]^\cong\\
H_i(X_n;A)\ar[rr]^\cong &&H_i(EX_n;A)
}$$
Hence $X\to Z_\infty$ lies in $S$. Also 
$\tilde{H}_i(Z_\infty)\stackrel{\cong}{\longrightarrow}\tilde{H}_i(EX)$
for $n\gg i$
so $\tilde{H}_*(Z_\infty)$ is an $A$-module. This concludes the proof of 5.7.
\hfill$\Box$

\vspace{10pt}

\noindent
{\bf Corollary 5.9.} (after Serre). \textit{Let $X$ and $Y$ be spaces such
that $\pi_1(X)$, $\pi_1(Y)$ are nilpotent and act trivially on $\pi_n(X)$,
$\pi_n(Y)$ for $n>1$. Then $f\in S$, i.e.,
$$f_*:H_*(X;A)\longrightarrow H_*(Y;A)$$
is an isomorphism, if and only if
$$f_\#:\pi_i(X)\otimes A\longrightarrow\pi_i(Y)\otimes A$$
is an isomorphism. (For $i=1$ this means that the subquotients of the central series map isomorphically after tensoring with $A$.)}

\vspace{10pt}

\noindent
{\bf First proof.} (i) Suppose
$f_\#:\pi_i(X)\otimes A\longrightarrow\pi_i(Y)\otimes A$
is an isomorphism for all $i$. Then 5.5 shows that for all $i$
$$(*)\qquad\qquad\qquad\qquad\tilde{H}_*\left(K(\pi_i(X),i);A\right)\longrightarrow
\tilde{H}_*\left(K(\pi_i(U),i);A\right)\qquad\qquad\qquad\qquad$$
is an isomorphism. Now we work up the Postnikov system by induction.
For if we assume
$$H_*\left(X(1,\dots,n-1);A\right)\longrightarrow H_*\left(Y(1,\dots,n-1);A\right)$$
is an isomorphism, then applying $(*)$ with $i=n$ and the Serre spectral sequence, we conclude that
$$H_*\left(X(1,\dots,n);A\right)\longrightarrow H_*\left(Y(1,\dots,n);A\right)$$
is an isomorphism. Of course, to deal with $\pi_1$ we take as many steps as there are subquotients in the central series. By taking $n$ sufficiently
large we see that
$$H_m(X; A)\longrightarrow H_m(Y; A)$$
is an isomorphism.

(ii) Conversely, suppose $H_*(X; A)\longrightarrow H_*(Y; A)$ is an isomorphism. Suppose we have proved that
$$\pi_i(X)\otimes A\longrightarrow\pi_i(Y)\otimes A$$
is an isomorphism for $i<n$. Then as in part (i),
$$H_*\left(X(1,\dots,n-1);A\right)\longrightarrow H_*\left(Y(1,\dots,n-1);A\right)$$
is an isomorphism. By the five-lemma
$$H_*\left(X(1,\dots,n-1),X;A\right)\longrightarrow H_*\left(Y(1,\dots,n-1),Y;A\right)$$
is an isomorphism. But since $\pi_1$ operates trivially on $\pi_n$ we have
$$\xymatrix{
H_{n+1}\left(X(1,\dots,n-1),X;A\right)\ar[rr]^\cong
&&H_{n+1}\left(Y(1,\dots,n-1),Y;A\right)\\
\pi_{n+1}\left(X(1,\dots,n-1),X\right)\otimes A\ar[u]_\cong\ar[rr]
&&\pi_{n+1}\left(Y(1,\dots,n-1),Y\right)\otimes A\ar[u]_\cong\\
\pi_n(X)\otimes A\ar[rr]\ar[u]_\cong&&\pi_n(X)\otimes A\ar[u]_\cong
}$$
This completes the induction. Of course, for $n=1$ we need as many steps as there are in your central series, and the argument requires slight adaptation
(replace $\pi_n$ by the relevant subquotient).\hfill$\Box$

\vspace{10pt}

\noindent
{\bf Second proof.} We have the localisation functor already. So
$$f_*:H_*(X;A)\longrightarrow H_*(Y;A)$$
is an isomorphism if and only if $f\in S$, if and only if $Ef$ is an equivalence,
if and only if
$$(Ef)_\#:\pi_*(EX)\longrightarrow\pi_*(EY)$$
is an isomorphism, if and only if
$$f_\#:\pi_*(X)\otimes A\longrightarrow\pi_*(Y)\otimes A$$
is an isomorphism.\hfill$\Box$

\vspace{10pt}

\noindent
{\bf Corollary 5.10.} \textit{Let $Y$ be a space such that $\pi_1(Y)$ is nilpotent
and acts trivially on $\pi_n(Y)$ for $n>1$. Then $Y\in D$ if and only if 
$\tilde{H}_*(Y)$ is local, if and only if $\pi_*(Y)$ is local.}

\vspace{10pt}

\noindent
{\bf Proof.} Consider	
$$\xymatrix{
\tilde{H}_*(Y)\ar[rr]^b\ar[d]^a &&\tilde{H}_*(EY)\ar[d]^\cong
&\pi_*(Y)\ar[rr]^d\ar[d]^c &&\pi_*(EY)\ar[d]^\cong\\
\tilde{H}_*(Y;A)\ar[rr]^\cong &&\tilde{H}_*(EY;A)
&\pi_*(Y)\otimes A\ar[rr]^\cong &&\pi_*(EY)\otimes A
}$$
Then $\tilde{H}_*(Y)$
is local iff $a$ is an isomorphism, iff $b$ is an isomorphism, iff $c$ is an isomorphism (5.9 with $A=\mathbf{Z}$, iff $Y\in D$, iff $d$ is an isomorphism, iff
$\pi_*(Y)$ is local.\hfill$\Box$

\vspace{10pt}

\noindent
{\bf Corollary 5.11.} \textit{Let $X$, $Y$ be spaces such that $\pi_1$ is nilpotent and acts trivially on $\pi_n$ for $n>1$. Then the following conditions on
$f:X\to Y$ are equivalent.
\begin{itemize}
\item[(i)] There is a commutative diagram
$$\xymatrix{
&&Y\ar@{<->}[dd]^\cong\\
X\ar[urr]_f\ar[drr]_{\eta_X}\\
&&EX
}$$
\item[(iia)] $f_\#:\pi_*(X)\longrightarrow\pi_*(Y)$ is a localising map.
\item[(iib)] $f_*:\tilde{H}_*(X)\longrightarrow\tilde{H}_*(Y)$ is a localising map.
\item[(iiia)] $f_\#:\pi_*(X)\otimes A\longrightarrow\pi_*(Y)\otimes A$ is an isomorphism and is couniversal with that property (in the category of such $X$, $Y$).
\item[(iiib)] $f_*:H_*(X; A)\longrightarrow H_*(Y; A)$ is an isomorphism and is
couniversal with that property (in the category of such $X$, $Y$).
\item[(iva)] $\pi_*(Y)$ is local (i.e., an $A$-module) and $f$ is universal
with that property (in the category of such $X$, $Y$).
\item[(ivb)]  $\tilde{H}_*(Y)$ is local (i.e., an $A$-module) and $f$ is universal
with that property (in the category of such $X$, $Y$).
\end{itemize}}

\vspace{10pt}

This is precisely the statement of 2.11 with the alternative descriptions of $S$ and $D$ inserted from 5.9, 5.10.

\vspace{20pt}

\noindent
[A section suggesting possible directions for further work on idempotent functors is omitted here.]

\vspace{10pt}

\noindent
{\bf Editorial Note.} The results of this section can be strengthened with a little additional effort. Namely
one can weaken the hypothesis that the fundamental group is nilpotent and acts trivially on the higher
homotopy groups to requiring that the fundamental group act nilpotently on the higher homotopy groups.
This means that each higher homotopy group has a finite filtration such that the fundamental group acts trivially on the filtration quotients. Such spaces are called nilpotent spaces.  The proof is along
exactly the same lines as above, except that one refines the Postnikov tower of the space by a  larger tower
of fibrations corresponding to the filtration quotients in the higher homotopy groups.  See \cite{MP}
for details.

\newpage

\section{Profinite Completion}\label{sec6}

We start with the most classical example.

\vspace{10pt}

\noindent
{\bf Example 6.1.} Let $G$ be a topological group. Then we can consider all homomorphisms from $G$ to finite groups $F_\alpha$ (with the discrete topology). For our purposes we can replace $F_\alpha$ by the image of $G\to F_\alpha$; so it will be sufficient to consider all closed normal subgroups $N_\beta\subset G$ of finite index and consider the projections $G\to G/N_\beta$. Such normal subgroups are
ordered by inclusion, $N_\beta\subset N_\gamma$, so we get an inverse system of finite discrete groups, $G/N_\beta\to G/N_\gamma$. We define the profinite completion of $G$ to
be $\displaystyle\widehat{G}=\lim_{\stackrel{\longleftarrow}{\beta}}G/N_\beta$. It comes equipped with a compact Hausdorff topology,
the topology of the inverse limit; and there is a canonical map $G\to\widehat{G}$.

\vspace{10pt}

In this way we get an idempotent functor on the category of topological groups.

Of course, if you start with a discrete group $G$ and complete it, you get a topological group, say $\Co G$. You can then apply the functor $F$ which forgets the topology. But the composite functor $F\Co$ is not idempotent.

\vspace{10pt}

\noindent
{\bf Example 6.2.} Take $V$ to be a vector space over $\mathbf{F}_p$, with the
discrete topology. Then $\Co V=\widehat{V}\cong V^{**}$, the double dual of $V$, with the topology of the double dual. So $F\Co V\cong V^{**}$ with the discrete topology.
So $F\Co F\Co V\cong V^{****}$ which is not isomorphic to $V^{**}$ if $V$ is infinite-dimensional.

\vspace{10pt}

Now let me try to explain Sullivan's concerns. Suppose $X$ is say an affine algebraic variety $\subseteq\mathbf{C}^n$, defined by polynomial equations with integer coefficients. Then the automorphism $z\mapsto\overline{z}$ of $\mathbf{C}$ defines an automorphism of $\mathbf{C}^n$ and an automorphism of $X$. So it gives an automorphism of $H_*(X)$ or of any other topological invariant of $X$.

Now suppose we have an invariant $K(X)$ which can be defined by purely algebraic means starting from the equations which define $X$. Then,we reason that any algebraic automorphism of $\mathbf{C}$ which leaves $\mathbf{Z}$ fixed should induce an automorphism of $ K(X)$. In other words, the Galois group 
$\mbox{Gal}(\mathbf{C}/\mathbf{Q})$ should act on $K(X)$. Now, this is a monstrous big symmetry group. However, most of it should be irrelevant. Purely algebraic constructions should hardly require more of $\mathbf{C}$ than 
$\overline{\mathbf{Q}}$, the algebraic closure of $\mathbf{Q}$. So we might expect the action of $\mbox{Gal}(\mathbf{C}/\mathbf{Q})$ to factor through the quotient group
$\mbox{Gal}(\overline{\mathbf{Q}}/\mathbf{Q})$. But this is still a pretty large group, and such symmetry is very useful.

This still leaves open, of course the question of what invariants of $X$ can be defined by purely algebraic means; and in particular, how close we can come to the usual topological invariants. This question
has been much studied, and the original studies concerned the fundamental group. It is equivalent to study the covering spaces of $X$. It turns out that what you can do by purely algebraic means is to construct the finite covering spaces of $X$. What you cannot do by purely algebraic means is to
construct any infinite covering spaces of $X$. More precisely, let $p:\widetilde{X}\to X$ be a regular map between algebraic varieties; if $x_0\in X$, the counterimage $p^{-1}x_0$ has to be an algebraic subset of $X$, so if it is 0-dimensional it has to be finite. Therefore, what you can learn about the fundamental group $\pi_1(X)$ is exactly the finite quotients $\pi_1(X)/N_\beta$; so you can learn the profinite completion $\widehat{\pi_1(X)}$, which
indeed serves to summarize the information about the finite quotients.

The further developments of this story are the business of \'etale homotopy theory, and I don't claim to understand them. I am credibly informed that what you can do by purely algebraic means is the following: you can recover $[X, Z]$ whenever $Z$ is a complex all of whose homotopy groups are finite. For example, taking $Z$ to be an Eilenberg-Mac\,Lane space of type $(\pi,1)$, where $\pi$ is finite, you can recover the homomorphisms $\pi_1(X)\to\pi$.

What Sullivan now does is to construct a functor $\widehat{X}$ of $X$ which summarizes all the information you can obtain from $[X, Y]$ where $Y$ runs over complexes with finite homotopy groups.

You might think that we already have a plausible candidate for such a functor. That is, we define a subset $S$ as follows: a map
$f: X\to Y$ is to lie in $S$ if and only if $f^*:[Y, Z]\longrightarrow[X, Z]$ is an isomorphism whenever $Z$ has finite homotopy groups. We can replace this by two conditions:
\begin{itemize}
\item[(i)] $f^*:\mbox{Hom}\left(\pi_1(Y),F\right)\longrightarrow
\mbox{Hom}\left(\pi_1(X),F\right)$ is an isomorphism whenever $F$ is a finite group,
\item[(ii)] $f^*:H^*(Y;A)\longrightarrow H^*(X;A)$ is an isomorphism whenever $A$ is
a finite abelian group with operations from $\pi_1(Y)$.
\end{itemize}

The necessity of (i) is seen by taking $Z$ to be an Eilenberg-Mac\,Lane space of type
$(F,1)$. The necessity of (ii) is seen by taking $Z$ to be the obvious fibering
$$K(A,n)\longrightarrow Z\longrightarrow K(\mbox{Aut}(A),1).$$ 
(Since A is a finite group, $\mbox{Aut}(A)$ is a finite group.)

The sufficiency of (i) plus (ii) is seen by obstruction theory.

The subset S clearly satisfies 3.1, 3.5 and 3.6. We see it satisfies 4.1 by checking (i) and (ii); here (i) comes easily from
Van Kampen's Theorem and (ii) comes by the same argument we used in 4.5.

Modulo 3.4, we get a functor $EX$. It is plausible to conjecture
that it coincides with Sullivan's $\widehat{X}$ in the main cases of interest. In fact, Sullivan only proves that his functor $\widehat{X}$ has good properties under rather restrictive conditions ($\pi_1(X)=0$ and $\pi_n(X)$ finitely generated for each $n$; see Sullivan \cite[Theorem 3.9, p. 330]{Su}.) But since Sullivan indicates that these conditions can be weakened, let us disregard the conditions and turn to Sullivan's conclusions. The conclusion 3.9(iv) certainly implies that Sullivan's
map $X\to\widehat{X}$ lies in our class $S$; and
whenever this conclusion holds, then Sullivan's $\widehat{X}$ is equivalent to
our $EX$. It may even be the case that Sullivan's $\widehat{X}$ is equivalent to our $EX$ whenever $X$ is a finite complex. This would be completely satisfactory, because the main application arises when $X$ is an algebraic variety, and algebraic varieties tend to be equivalent to finite complexes. It is equally certain that I can't prove it because of genuine difficulties with the fundamental group. For example, one might like to know that if $G$ is a finitely presented group, such as
$\pi_1(X)$, then any homomorphism from $\widehat{G}$ to a finite group is necessarily continuous. So far as I know, this is not known.

There are reasons why $\widehat{X}$ and $EX$ should not be equivalent in general. The first is this
Sullivan 's approach assumes that you are given $X$ and know 
$[X, Z]$ whenever $Z$ has finite homotopy. The construction of $EX$ tends to assume that you know $[Y, Z]$ whenever Z has finite homotopy groups and $Y$ runs over a lot of other things as well as $X$. Clearly if you take $Y$ to be something which is not an algebraic variety you are feeding in information not accessible to algebraic geometry, so in this sense it is not obvious that algebraic geometry will suffice to construct $EX$. The second reason is this. Sullivan's approach does not give an idempotent functor defined on the category $\C$ we have been considering. If it is to have any chance to be idempotent it must be defined on a different category.

To explain this, I need the notion of a topologised object in a category. One way to define a group object in a category is as follows. It is an object $G$ such that the set $[X, G]$ is given the structure of a group for each $X$ so that the induced functions $f^*:[Y, G]\longrightarrow[X, G]$ are group homomorphisms. Similarly, we say that $T$ is a topologised object in a category $C$ if for each $X$ in $C$,
$[X, T]$ is given the structure of a topological space, and for each $f:X\to Y$
in $C$, the induced function $f^*:[Y, T]\longrightarrow[X, T]$ is continuous. Let $T$ and $U$ be topologised objects in $C$; we say that a map $g:T\to U$ in $C$ is a map of
topologised objects if $g^*:[X, T]\longrightarrow[X, U]$ is continuous for each $X$.

\vspace{10pt}

\noindent
{\bf Example 6.3.} Let $\pi$ be a compact (Hausdorff) abelian group, for example,
$S^1$ or the $p$-adic integers. Then the Eilenberg-Mac\,Lane space $K(\pi,n)$ is a topologised object in the homotopy category. In fact, for any $X$ we topologise
$C^q(X;\pi)$ with the product topology; we take the obvious topology on cocycles, coboundaries and $H^n(X;\pi)$. This gives a topology in $[X, K(\pi,n)]$ for each $n$ (indeed making it a compact group), and the induced functions are continuous. A continuous homomorphism $\theta:\pi_1\to\pi_2$
induces a map of the topologised objects $K(\pi_1,n)\longrightarrow K(\pi_2,n)$.

\vspace{10pt}

I will say that $T$ is a compact Hausdorff object in $C$ if it is topologised and the topology on $[X, T]$ is compact and Hausdorff for each $X$. Example 5.3 is of this kind. However, Sullivan speaks of a
compact representable functor; that is, he identifies $T$ with the representable functor $[X, T]$.

In particular, let $\C$ be the homotopy category, and let $\C^T$ be the category of topologised objects in $\C$. Then our ignorance of $\C^T$ is really deplorable; in fact, it is one of the objects of this lecture to persuade some of you to give a little thought to $\C^T$. Anyway, Sullivan's profinite completion functor is a functor $\Co$ from $\C$ to $\C^T$. Of course, we can compose it with the functor $F:\C^T\to \C$, which forgets the topology; but the composite $F\Co$ is not idempotent. (Take $X$ to be an Eilenberg-Mac\,Lane space of type $(V, n)$ where $V$ is as in Ex. 5.2.) Therefore, the only chance of making Sullivan's profinite completion idempotent is to set it up as a functor from $\C^T$ to $\C^T$ (compare 5.1). I do not know whether Sullivan's profinite completion can be set up as a functor from $\C^T$ to $\C^T$;
a fortiori I do not know if it can be set up as an idempotent functor from $\C^T$ to $\C^T$.

At this point, I should start to do the work. What follows is my interpretation of a letter from Sullivan.

\vspace{10pt}

\noindent
{\bf Proposition 6.4.} \textit{Let $Z\in \C$ be an object with finite homotopy groups. Then there is one and only one way to give $Z$ the structure of a compact Hausdorff object. Any map $Z\to Z'$ between such objects is a map of topologised objects.}

\vspace{10pt}

\noindent
{\bf Proof.} First let $X$ be a finite complex. Then $[X, Z]$ is a finite set, and there is one and only one way to make it a compact Hausdorff space.

Next let $X$ be a general complex, and define 
$\displaystyle H(X)=\lim_{\stackrel{\longleftarrow}{\alpha}}[X_\alpha, Z],$
where $X_\alpha$
runs over all the finite subcomplexes of $X$. We topologise 
$\displaystyle\lim_{\stackrel{\longleftarrow}{\alpha}}[X_\alpha, Z]$
with the inverse-limit topology. Clearly $H$ is a contravariant
functor from the category $\C$ to the category of compact Hausdorff spaces. I claim $H$ is representable. To check the wedge axiom is automatic. To check the Mayer-Vietoris axiom is straightforward, but one uses in an essential way the fact that the inverse limit of a system of nonempty finite sets is nonempty. We conclude that 
$\displaystyle H(X)\stackrel{\cong}{\longleftrightarrow}[X,\overline{Z}]$, where $\overline{Z}$ is a representing object. Now we clearly have a natural transformation
$$[X,Z]\longrightarrow \lim_{\stackrel{\longleftarrow}{\alpha}}[X_\alpha, Z]
=H(X)\stackrel{\cong}{\longleftrightarrow}[X,\overline{Z}]$$
This must be induced by a map $Z\to\overline{Z}$. Now if $X$ is a sphere then
clearly $\displaystyle[X,Z]\longrightarrow \lim_{\stackrel{\longleftarrow}{\alpha}}[X_\alpha, Z]$
is an isomorphism because the limit is a attained. Therefore $Z\to\overline{Z}$
induces isomorphisms of homotopy groups, and is an equivalence. So for any $X$, the obvious map
$\displaystyle[X,Z]\longrightarrow\lim_{\stackrel{\longleftarrow}{\alpha}}[X_\alpha, Z]$ is an isomorphism. We already topologised
$\displaystyle\lim_{\stackrel{\longleftarrow}{\alpha}}[X_\alpha, Z]$, so this puts a compact Hausdorff topology on $[X, Z]$. We
check that we have made $Z$ a compact Hausdorff object, and such a structure on $Z$ is unique. For any map $Z\to Z'$ the map
$[X_\alpha, Z]\longrightarrow[X_\alpha, Z']$ is continuous, so
 $\displaystyle\lim_{\stackrel{\longleftarrow}{\alpha}}[X_\alpha, Z]\longrightarrow
\lim_{\stackrel{\longleftarrow}{\alpha}}[X_\alpha, Z']$
is continuous.\hfill$\Box$

\vspace{20pt}

Now we turn to the construction of Sullivan's profinite completion $Y\to\widehat{Y}$.

\newpage

Let $Y$ be a fixed object in $\C$. We consider a category:
the objects are maps $g_\alpha:Y\to Z_\alpha$ such that $Z_\alpha$ has finite
homotopy groups; the arrows are diagrams
$$
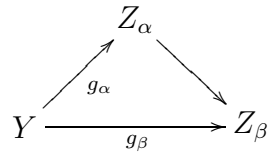
\xymatrix{
&Z_\alpha\ar[dr]\\
Y\ar[ur]_{g_\alpha}\ar[rr]_{g_\beta}&&Z_\beta
}$$
We check that this is a directed category (directed away from infinity).
Given two objects, I have to find one further from infinity than each. This is easy: given
$$
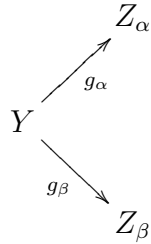
\xymatrix{
&Z_\alpha\\
Y\ar[ur]_{g_\alpha}\ar[dr]_{g_\beta}\\
&Z_\beta
}$$
I construct
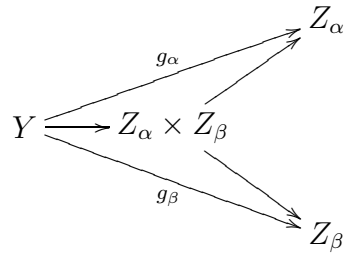
$$\xymatrix{
&&Z_\alpha\\
Y\ar[urr]^{g_\alpha}\ar[drr]_{g_\beta}\ar[r]
&Z_\alpha\times Z_\beta\ar[ur]\ar[dr]\\
&&Z_\beta
}$$
and $Z_\alpha\times Z_\beta$ has finite homotopy groups.
Secondly, suppose we are given
a diagram
$$\xymatrix{\bullet\ar@<1ex>[r]\ar@<-1ex>[r]&\bullet};$$
I have to construct
$$\xymatrix{{}\ar@{-->}[r]&\bullet\ar@<1ex>[r]\ar@<-1ex>[r]&\bullet}.$$
So suppose given
$$
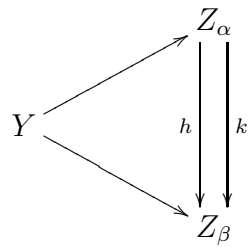
\xymatrix{
&&Z_\alpha\ar@<1ex>[dd]^k\ar@<-1ex>[dd]_h\\
Y\ar[urr]\ar[drr]\\
&&Z_\beta
}$$
Form
$$\xymatrix{
&&Z_\alpha\ar[d]^{(1,k)}\\
Z_\alpha\ar[rr]^{(1,h)}&&Z_\alpha\times Z_\beta
}$$
and take the homotopy pullback
$$\xymatrix{
F\ar[rr]\ar[d]&&Z_\alpha\ar[d]^{(1,k)}\\
Z_\alpha\ar[rr]^{(1,h)}&&Z_\alpha\times Z_\beta
}$$
(If I want to stay inside the category of connected spaces, I take the component of the basepoint.) Then $F$ has finite homotopy groups, and it is easy to show that there is a map $Y\to F$ which does as required.

It is very tempting to modify this construction, taking $Y$ to be
a fixed object in $\C^T$ and considering maps $g_\alpha:Y\to Z_\alpha$ in $\C^T$.
The difficulty comes at the last step above: how can we know that we can choose the map $Y\to F$ to be a map of topologised objects?

Now that we have a directed category, we form 
$\displaystyle H(X)=\lim_{\stackrel{\longleftarrow}{\alpha}}[X,Z_\alpha]$.
Each set $[X, Z_\alpha]$ is a compact Hausdorff space by 5.4, so we give
$\displaystyle\lim_{\stackrel{\longleftarrow}{\alpha}}[X,Z_\alpha]$ the inverse-limit topology, and it becomes a compact
Hausdorff space. Thus $H$ becomes a contravariant functor from $\C$ to
compact Hausdorff spaces.

\vspace{10pt}

\noindent
{\bf Lemma 6.5.} \textit{H is representable.}

\vspace{10pt}

\noindent
{\bf Proof.} The proof of the wedge axiom is automatic. The proof of the Mayer-Vietoris axiom is straightforward, but one makes essential use of the fact that an inverse limit of nonempty compact Hausdorff spaces is nonempty.\hfill$\Box$

\vspace{10pt}

So we have 
$\displaystyle\lim_{\stackrel{\longleftarrow}{\alpha}}[X,Z_\alpha]\stackrel{\cong}{\longleftrightarrow}
[X,\widehat{Y}]$ for an object $\widehat{Y}$ in $\C$. Since
we have given $\displaystyle\lim_{\stackrel{\longleftarrow}{\alpha}}[X,Z_\alpha]$ a compact Hausdorff topology, $\widehat{Y}$ is a compact Hausdorff object in $\C$.
We have an obvious natural transformation.
$$[X, Y]\longrightarrow\lim_{\stackrel{\longleftarrow}{\alpha}}[X,Z_\alpha]\stackrel{\cong}{\longleftrightarrow}[X,\widehat{Y}];$$
this must be induced by a map $Y\to\widehat{Y}$.

It is rather clear that the homotopy groups of $\widehat{Y}$ are profinite:
we have 
$\displaystyle[S^n,\widehat{Y}]\stackrel{\cong}{\longleftrightarrow}
\lim_{\stackrel{\longleftarrow}{\alpha}}[S^n,Z_\alpha]$
where $[S^n,Z_\alpha]$ is finite.

\vspace{10pt}

The following are among the main properties of Sullivan's profinite completion.

\vspace{10pt}

\noindent
{\bf Theorem 6.6.} (Sullivan \cite[Theorem 3.9, p. 330]{Su})
\textit{Suppose $Y$ is simply connected and its homotopy groups
$\pi_n(Y)$ are finitely-generated for each $n$. Then
\begin{itemize}
\item[(i)] The homotopy groups of $\widehat{Y}$ are profinite and
universal with this property: if
$$
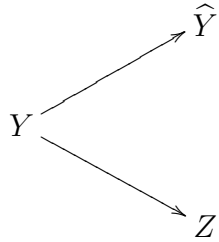
\xymatrix{
&&\widehat{Y}\\
Y\ar[urr]\ar[drr]\\
&&Z}$$
and the homotopy groups of $Z$ are profinite, then
$$
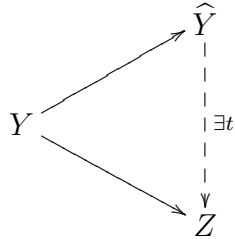
\xymatrix{
&&\widehat{Y}\ar@{-->}[dd]^{\exists t}\\
Y\ar[urr]\ar[drr]\\
&&Z}$$
\item[(ii)] More precisely, 
$\widehat{\pi_i(Y)}\stackrel{\cong}{\longrightarrow}\pi_i(\widehat{Y})$.
\item[(iii)] The topological structure of $\widehat{Y}$ is determined by
the underlying object. In fact
$$\Co F\Co Y\cong \Co Y.$$
\end{itemize}
}

\vspace{10pt}

The proof is to see what profinite completion does to the Postnikov system.

\vspace{10pt}

\noindent
{\bf Editorial Note.} It turns out that Sullivan's profinite completion functor, for spaces satisfying the
hypothesis of Theorem 6.6, is a special case of the localisation functor of Theorem 4.6. Namely if one
takes as $S$ the class of morphisms which induce isomorphisms in $H_*(-;A)$, where $A$ is a torsion
abelian group which has $p$-torsion for all primes $p$ (e.g. $A=\bigoplus_p\mathbf{Z}/p\mathbf{Z}$), then
the corresponding idempotent functor $EY$ is naturally equivalent to $\widehat{Y}$.  Again a good
reference is \cite{MP}.

\newpage

\section{Use of Brown-Peterson Homology in Stable Homotopy}\label{sec7}

You may have heard that we have a contradiction in homotopy theory these days.
There is supposed to be a certain map $\gamma_1$; and
Zahler and Thomas claim to have proved that $\gamma_1\ne0$, while Oka and
Toda claim to have proved that $\gamma_1=0$ (\cite[Remark on p. 147]{O}). 
The situation reminds me of these words of Andr\'e Weil \cite{W} about the need for rigour in mathematics: \textit{``It is possible for the advancing army to outrun its services of supply
and incur disaster unless it waits for the quartermaster to perform his
inglorious but indispensable tasks.''} In the present situation the duty of
the quartermaster is clearly to call in both parties to the dispute and audit their books.

I have tried to repeat the work of Zahler and Thomas. What I would like you to do is to watch with the most critical attention while I try to perform it in public. I don't think I've got anything up my sleeve, but someone, somewhere has got a joker hidden away; probably not because he put it there, but more likely because he accepted someone else's work in too trusting and uncritical a fashion.

I must begin by explaining the element $\gamma_1$ and for that I need the spaces $V(n)$
of Toda and others. Let $p$ be an odd prime. Let us recall that according to Milnor, the mod $p$ homology of the Eilenberg-Mac\,Lane spectrum $H\mathbf{Z}_p$ is the tensor product of an exterior algebra and a polynomial algebra
$$(H\mathbf{Z}_p)_*(H\mathbf{Z}_p)=\Lambda[\tau_0,\tau_1,\dots]\otimes
\mathbf{Z}_p[\xi_1,\xi_2,\dots].$$
The Brown-Peterson spectrum $BP$ admits a map $f:BP\to H\mathbf{Z}_p$ so that
$$f_*:(H\mathbf{Z}_p)_*(BP)\to (H\mathbf{Z}_p)_*(H\mathbf{Z}_p)$$
is mono, and its image is $\mathbf{Z}_p[\xi_1,\xi_2,\dots]$. Equivalently
$(H\mathbf{Z}_p)_*(BP)$ is a cyclic module over the mod $p$ Steenrod algebra $A^*$ generated by one generator $f$ and is isomorphic to the quotient $A^*/A^*\beta_p A^*$, where $\beta_p$ denotes the Bockstein operation.

Similarly, we would admit $X$ as a spectrum of type $V(\infty)$ if it had a map
$f:X\to H\mathbf{Z}_p$ such that
$$f_*:(H\mathbf{Z}_p)_*(X)\to (H\mathbf{Z}_p)_*(H\mathbf{Z}_p)$$
is mono and its image is $\Lambda[\tau_0,\tau_1,\dots]$. It would be well to specify also that $X$ should be connected and $(H\mathbf{Z}_q)_*(X)=0$ for $q$ prime to $p$. If such a spectrum existed, we could form
$$V(\infty)\wedge BP\longrightarrow H\mathbf{Z}_p\wedge H\mathbf{Z}_p\longrightarrow
H\mathbf{Z}_p,$$
and this would induce an isomorphism of $(H\mathbf{Z}_p)_*$ so that 
$V(\infty)\wedge BP$ would be an Eilenberg-Mac\,Lane spectrum $H\mathbf{Z}_p$.	However, no such spectrum $V(\infty)$ exists, so we turn to finite approximations to it.

We admit $X$ is a spectrum of type $V(n)$ if it has a map 
$f:X\to H\mathbf{Z}_p$ such that
$$f_*:(H\mathbf{Z}_p)_*(X)\to (H\mathbf{Z}_p)_*(H\mathbf{Z}_p)$$
is mono and its image is $\Lambda[\tau_0,\tau_1,\dots,\tau_n]$.
Essentially we are prescribing $(H\mathbf{Z}_p)_*(X)$ as a module over the mod $p$ Steenrod algebra. Again we specify that $X$ should be connected and
$(H\mathbf{Z}_q)_*(X)=0$ for $q$ prime to $p$.

\vspace{10pt}

\noindent
{\bf Example.} $V(0)$ exists; we can take $V(0)=S^0\cup_p e^1$.

\vspace{10pt}

\noindent
{\bf Example.} $V(1)$ exists. More precisely, there is a map 
$\alpha:S^{2p-2}\longrightarrow S^1$
such that the operation $P^1$ is non-zero in 
$(H\mathbf{Z}_p)_*(S^1\cup_\alpha e^{2p-1}$. Since
$p\alpha=0$ and $\{p,\alpha,p\}=0$, we can find a map $A$ to fill the diagram below
$$\xymatrix{
S^{2p-1}\cup_p e^{2p-1}\ar@{-->}[rr]^A &&S^0\cup_p e^1\ar[d]\\
S^{2p-2}\ar[u]\ar[rr]^\alpha &&S^1
}$$
The mapping cone of $A$ is a complex
$$V(1) = V(0)\cup_A CS^{2p-2}V(0)=S^0\cup_p e^1\cup_\alpha e^{2p-1}\cup_p e^{2p}$$
in which the operation $\beta P^1\beta$ is non-zero.

\vspace{10pt}

It was shown by Larry Smith \cite{Sm} that $V(2)$ exists for $p\ge 5$ and Toda shows that
$V(3)$ exists for $p\ge 7$ ( \cite[p. 53, Thm. 1.1]{T}).\footnote{L. S. Nave \cite{N} has shown that for any $p$,
$V(n)$ does not exist for $n\ge\frac{p+1}{2}$.}
Toda also proves that these $V(n)$ are unique up to homotopy equivalence (Theorem 4.1, p. 57).
 This is not actually essential for our purposes, but it is convenient.

Now we want to see how the spaces $V(n)$ give rise to systematic families of elements in homotopy theory. We have seen that there is a
map	
$$A:S^{2(p-1)}V(0)\longrightarrow V(0).$$
By iterating it we can construct
$$\xymatrix{
S^{2r(p-1)}V(0)\ar[r]&\dots\ar[r] &S^{4(p-1)}V(0)\ar[rr]^{S^{2(p-1)}A}
&&S^{2(p-1)}V(0)\ar[r]^{\quad\  A} &V(0)\ar[d]\\
S^{2r(p-1)}\ar[u]\ar[rrrrr]^{\alpha_r}&&&&&S^1
}$$
So we construct $\alpha_r\in\pi^s_{2r(p-1)-1}$.
It is usual to write $q=2(p-1)$; 
then $\alpha_r$ is of degree $rq-1$.

Toda clearly asserts that there is a map
$$B:S^{2(p^2-1)}V(1)\longrightarrow V(1)\quad\mbox{for }p\ge 5$$
whose mapping cone is $V(2)$, and a map
$$C:S^{2(p^3-1)}V(2)\longrightarrow V(2)\quad\mbox{for }p\ge 7$$
whose mapping cone is $V(3)$ (Toda \cite[Corollary 4.3, p. 58]{T}).

By iterating $B$ we can construct
$$\xymatrix{
S^{2r(p^2-1)}V(1)\ar[r]&\dots\ar[r] &S^{4(p^2-1)}V(1)\ar[rr]^{S^{2(p^2-1)}B}
&&S^{2(p^2-1)}V(1)\ar[r]^{\quad\  B} &V(1)\ar[d]\\
S^{2r(p^2-1)}\ar[u]\ar[rrrrr]^{\beta_r}&&&&&S^{2p}
}$$
$\beta_r$ is of degree $2r(p^2-1)-2p=r(p+1)q-q-2$. By iterating $C$ we
can construct
$$\xymatrix{
S^{2r(p^3-1)}V(2)\ar[r]&\dots\ar[r] &S^{4(p^3-1)}V(2)\ar[rr]^{S^{2(p^3-1)}C}
&&S^{2(p^3-1)}V(2)\ar[r]^{\quad\  C} &V(2)\ar[d]\\
S^{2r(p^3-1)}\ar[u]\ar[rrrrr]^{\gamma_r}&&&&&S^{(p+2)q+3}
}$$
$\gamma_r$ has degree $2r(p^3-1)-(p+2)q-3=r(p^2 + p + 1)q-(p+2)q-3$. In particular,
$\gamma_1$ has degree $(p^2-1)q-3$.

Now, of course, if $C$ and $V(3)$ do not exist, then there is nothing to discuss, but it is generally supposed that they do exist.

Now we must discuss the $BP$-homology and cohomology of our spaces $V(n)$. It is necessary to recall that $\pi_*(BP)$ is a polynomial
algebra over $\mathbf{Z}_{(p)}$ on generators 
$v_1,v_2,v_3,\dots,v_i,\dots$, of dimension $2(p^i-1)$.
We will make the choice of generators $v_i$ more  precise later. From
the cofibering
$$\xymatrix{
S^0\ar[r]^p &S^0\ar[r] &V(0)\ar[r] &S^1\ar[r]^p &S^1}$$
we see that $BP_*(V(0))\cong\pi_*(BP)/p\pi_*(BP)$ on one generator of
dimension 0, coming from the injection $S^0\to V(0)$. Similarly,
$$BP^*(V(0))\cong\pi_*(BP)/p\pi_*(BP)$$
on one generator of dimension 1, coming from the projection $V(0)\to S^1$.

We now observe that
$$(H\mathbf{Z}_p)_*(V(n)\wedge BP)\longrightarrow (H\mathbf{Z}_p)_*(H\mathbf{Z}_p)$$
is mono, and epi in dimension $< 2p^{n+1}-1$. Therefore
$$\pi_i(V(n)\wedge BP) = 0\quad\mbox{for }0 <i<2p^{n+1}-2$$
i.e.,
$$BP_i(V(n))= 0\quad\mbox{for }0 <i<2p^{n+1}-2.$$
In particular, if $g_n\in BP_0(V(n))$ comes from the injection $S^0\to V(n)$,
then we have
$$v_ig_n=0\quad\mbox{for }i\le n$$

Consider now the cofibering
$$\xymatrix{
S^qV(0)\ar[r]^A &V(0)\ar[r] &V(1)\ar[r] &S^{q+1}V(0)\ar[r]^{SA} &SV(0)}.$$
We have $v_1g_0\mapsto 0$ in $BP_q(V(1))$, so we conclude 
$$A_*g_0 = c_1v_1g_0,\quad c_1\in\mathbf{Z}_p$$
with $c_1\ne 0$. Then the cofibering gives
$$BP_*(V(1))\cong\pi_*(BP)/(p,v_1)$$
generated by $g_1$. Similarly, we see that
$$B_*g_1 = c_2v_2g_1,\quad c_2\in \mathbf{Z}_p$$
with $c_2\ne 0$, that
$$BP_*(V(2))\cong\pi_*(BP)/(p,v_1,v_2)$$
generated by $g_2$, that
$$C_*g_2 = c_3v_3g_2,\quad c_3\in \mathbf{Z}_p$$
with $c_3\ne 0$, and that
$$BP_*(V(3))\cong\pi_*(BP)/(p,v_1,v_2,v_3)$$
generated by $g_3$.

Similar results are true for $BP$-cohomology. To deduce them, recall that
$$BP^*(X)\cong BP_*(DX)$$
where $DX$ is the Spanier-Whitehead dual of $X$. But the Spanier-Whitehead dual of a space $V(n)$ is again a space $V(n)$. Therefore we have
\begin{eqnarray*}
BP^*(V(0)) &\cong &\pi_*(BP)/(p)\\
BP^*(V(1)) &\cong &\pi_*(BP)/(p, v_1)\\
BP^*(V(2)) &\cong &\pi_*(BP)/(p, v_1 , v_2)\\
BP^*(V(3)) &\cong &\pi_*(BP)/(p, v_1 , v_z, v_3)
\end{eqnarray*}
on generators $g_0$, $g_1$, $g_2$, $g_3$ of degrees 1, $1+(2p-1)$,
$1+(2p-1)+(2p^2-1)$, $1+(2p-1)+(2p^2-1)+(2p^3-1)$ coming from the projections
$V(0)\to S^1$, $V(1)\to S^{2p}$, etc. Similarly, we have
\begin{eqnarray*}
A^*g_0 &= &c_1'v_1g_0\\
B^*g_1 &= &c_2'v_2g_1\\
C^*g_2 &= &c_3'v_3g_2
\end{eqnarray*}
where $c_1',c_2',c_3'\in\mathbf{Z}_p$ and are non-zero.

Without loss of generality we may replace $A$ by an integer multiple $mA$ where
$m\not\equiv 0$ mod $p$. Similarly for $B$ and $C$. Therefore we may suppose that
\begin{eqnarray*}
A^*g_0 &= &v_1g_0\\
B^*g_1 &= &v_2g_1\\
C^*g_2 &= &v_3g_2.
\end{eqnarray*}
This is the normalisation usually adopted by those who work with $BP$-cohomology.

Now I want to explain the nature of the proof that $\gamma_1\ne 0$.
First we form the complex
$$X_3=S^0\cup_{\gamma_1} e^{(p^2-1)q-2}.$$
Then we inspect the cofibering
$$\xymatrix{
S^{(p^2-1)q-3}\ar[r]^{\quad\gamma_1} &S^0\ar[r]^{i\qquad\qquad}
&X_3=S^0\cup_{\gamma_1} e^{(p^2-1)q-2}\ar[r] &S^{(p^2-1)q-2}\ar[r] &S^1.}$$
Here $BP^*(S^0)$ is non-zero only in dimensions $\equiv 0$ mod $q$, and
$BP^*(S^{(p^2-1)q-3})$ is non-zero only in dimensions $\equiv -3$ mod $q$, and
we may certainly suppose $q\ge 8$, so $BP^*(\gamma_1)=0$. There is a unique element 
$h\in BP^0(X_3)$ projecting to the generator of $BP^0(S^0)$ and we
take $k$ to be the image of the generator in 
$BP^{(p^2-1)q-2}(S^{(p^2-1)q-2})$.
We will produce a cohomology operation $\chi$ in $BP$-cohomology such that $\chi h$ is a non-zero multiple of $k$. (By ``multiple'' we mean a multiple by an element of
$\pi_*(BP)$.)	This will contradict the statement
$\gamma_1=0$, which implies $X_3\simeq S^0\vee S^{(p^2-1)q-2}$. For if so,
we would have
a map $j:X_3\to S^0$ such that $ji=1$; then $i^*j^*i^*h=i^*h$ and
$i^*:BP^0(X_3)\longrightarrow BP^0(S^0)$ is an isomorphism, so $h=j^*i^*h$. Then
$$\chi h\supseteq j^*i^*\chi h\supseteq j^*i^*(\mu k)\supseteq j^*0=0.$$
The only snag is that $\chi$ is a tertiary operation. However, when we see it we will be able to check that it is defined on $h$ and
therefore also on $i^*h$, and that its indeterminacy in $X_3$ is manageable
and its indeterminacy in $S^0$ is zero. I must now explain how to set up $\chi$ and calculate it.

For this I must explain about primary operations on $BP$-cohomology and the relations between them. In \cite{A2}, it is stated that
$$BP_*(BP)=\pi_*(BP)[t_1, t_2, t_3,\dots]$$
where the $t_i$ are well-defined elements such that $\mbox{deg }t_i=2(p^i-1)$
(Thm 16.1, p. 97). It follows that in $BP^*(BP)$ we have unique elements
$R_{i_1i_2\dots i_n}$ such that
$$\langle R_{i_1i_2\dots i_n},t_1^{j_1}t_2^{j_2}\cdots t_n^{j_n}\rangle
\ =\ \delta_{i_1j_1}\delta_{i_2j_2}\cdots\delta_{i_nj_n}.$$
In particular, we will need $R_1$, $R_p$ and $R_{01}$, which are the elements
of the ``dual base'' corresponding to $t_1,t_1^p,t_2$. They are of degree $q$, 
$pq$, $(p+1)q$.

\vspace{10pt}

\noindent
{\bf Lemma 7.1.} \textit{We have 
\begin{eqnarray*}
R_1R_p-R_pR_1 &= &R_{01}\\
R_1R_{01}-R_{01}R_1 &= &0\\
R_pR_{01}-R_{01}R_p &= &0.
\end{eqnarray*}}

\vspace{10pt}

\noindent
{\bf Proof.} By \cite[Prop. 3, p. 73]{A1}, the composition product 
in $BP^*(BP)$ is determined by the diagonal map in $BP_*(BP)$, according 
to the following formula. Suppose
$$\psi x=\sum_i e_i\otimes x_i$$
Then 
$$\langle ab,x\rangle\ =\ \sum_i (-1)^{|b|\,|e_i|} \langle a,e_i\langle b,x_i\rangle\rangle.$$
According to \cite[Th. 16.1, p. 98]{A2}, the diagonal map in $BP_*(BP)$
is given by the following formulae 
\begin{eqnarray*}
\psi t_1 &= &t_1\otimes 1+ 1\otimes t_2\\
\psi t_2 &= &t_2\otimes 1+t_1\otimes t_1^p+1\otimes t_2
-\sum_{\stackrel{i+j=p}{\scriptscriptstyle i>0,j>0}}\frac{(p-1)!}{i!j!}v_1t_1^i\otimes t_1^j.
\end{eqnarray*}
Here I need the fact that 
$$v_1=[CP^{p-1}]=p\,m_{p-1}$$
in the notation of the above. Moreover, I need to know that if I 
have a monomial $t^I=t_1^{i_1}t_2^{i_2}\cdots t_n^{i_n}$ of degree say $2d$,
then $\psi t^I$ has the form ${\sum c_{JK}t^J\otimes t^K}$, $c_{JK}\in\pi_*(BP)$, where
$0\le\mbox{deg }c_{JK}\le\frac{2d}{p}$, \,
$\frac{2d(p-1)}{p}<\mbox{deg}(t^J\otimes t^K)\le 2d$. 
So in evaluating a product $R^I\otimes R^J$ of 
degree $2e$, it will be sufficient to consider $\psi t^I$ with 
$\mbox{deg }t^I< \frac{2ep}{p-1}$. Hence in proving this lemma, we need not worry about monomials 
containing $t_3,t_4,\dots$ . The following tables give enough values of the 
pairing $\langle ab,t^I\rangle$ to prove the lemma.

\begin{center}
\begin{tabular}{l||c|c|c|c}
 &$t_1$ &$t_2$ &$t_1^{p+1}$ &$t_1t_2$\\
\hline\hline
$R_1R_p$ &$p+1$ &0 &1 &$-v_1$\\
$R_pR_1$ &$p+1$ &0 &0 &$-v_1$
\end{tabular}

\begin{tabular}{l||c|c|c}
 &$t_1^j$ &$t_1t_2$ &$t_2$\\
\hline\hline
$R_1R_{01}$ &0 &1 &0\\
$R_{01}R_1$ &0 &1 &0
\end{tabular}

\begin{tabular}{l||c|c|c|c|c|c}
 &$t_1^j$ &$t_1^pt_2$ &$t_1^{p+1}t_2$ &$t_1^{p+2}t_2$ &$t_2$ &$t_1t_2^2$\\
\hline\hline
$R_pR_{01}$ &0 &1 &0 &0 &0 &0\\
$R_{01}R_p$ &0 &1 &0 &0 &0 &0
\end{tabular}
\end{center}
\hfill$\Box$

\vspace{10pt}

Multiplying out these relations, I get:
\begin{eqnarray*}
R_pR_1^2-2R_1R_pR_1+R_1^2R_p &= &0\\
R_p^2R_1-2R_pR_1R_p+R_1R_p^2 &= &0.
\end{eqnarray*}
I now construct maps of spectra
$$\xymatrix{
S^{(2p+2)q}BP &S^{(p+2)q}BP\vee S^{(2p+1)q}BP\ar[l]_{d_2\qquad}
&S^qBP\vee S^{pq}\ar[l]_{\qquad\quad d_1} &BP\ar[l]_{\qquad d_0}}$$
with components given by the matrices
\begin{eqnarray*}
d_0 &=& \left[\begin{array}{c}R_1\\ R_p\end{array}\right]\\
d_1 &=& \left[\begin{array}{lc}R_pR_1-2R_1R_p &R_1\\ R_p^2 &R_1R_p\end{array}\right]\\
d_2 &= &\left[\begin{array}{lr}R_p &R_1\end{array}\right]
\end{eqnarray*}

Then by construction we have 
$$d_2d_1\simeq 0\qquad d_1d_0\simeq 0.$$
Consider also $h:X_3\to BP$. Then I claim we have $d_0h\simeq 0$. in fact, the maps
$$\xymatrix{
BP^q(X_3)\ar[rr]^{i^*} &&BP^q(S^0)\\
BP^{pq}(X_3)\ar[rr]^{i^*} &&BP^{pq}(S^0)
}$$
are isomorphisms and 
$$\xymatrix@R=5pt@C=1pt{
[S^0, S^qBP] &= &0\\
[S^0, S^{pq}BP ] &= &0\\
}$$
so $BP^q(X_3)=0$, $BP^{pq}(X_3)=0$. Hence the composite
$$\xymatrix{
S^0\ar[r] &X_3\ar[r]^h &BP\ar[r]^{d_0\qquad} &S^qBP\vee S^{pq}}$$
is zero, and so is $d_0h$. For the sake of convenience I will call the spectra
$$\xymatrix{
C_3 &C_2\ar[l]_{d_2} &C_1\ar[l]_{d_1} &C_0\ar[l]_{d_0}.}$$
Let me now nail my colours to the mast. From now on set $r=(p^2-1)q$.

\vspace{10pt}

\noindent
{\bf Theorem 7.2.} \textit{Let $X_3=S^0\cup_{\gamma_1}e^{r-2}$. Then
\begin{eqnarray*}
BP^*(X_3) &\cong &\pi_*(BP)\mbox{ on a generator $h$ of degree 0, restricting to the}\\
&&\mbox{\qquad\qquad generator of $S^0$}\\
&&\oplus\pi_*(BP)\mbox{ on a generator $\ell$ of degree $r-2$, coming from the}\\
&&\mbox{\qquad\qquad generator of $S^{r-2}$.}
\end{eqnarray*}
The operation $\{d_2,d_1,d_0,h\}$ is defined, and is a coset of maps
$$S^2X_3\longrightarrow C_3$$
and it is $-2v_2^{p-3}\ell \mod(p,v_1)\ell$.
}

\vspace{10pt}

At this time we only prove part of the theorem.
We have already noted that the facts about $BP^0(X_3)$ follow immediately from the cofibering 
$$\xymatrix{
S^{r-3}\ar[r]^{\gamma_1} &S^0\ar[r] &X_3\ar[r] &S^{r-2}\ar[r] &S^1.}$$

We first check that the operation is defined. As I stated above, $d_0h\simeq 0$.
In fact, we showed
$$\xymatrix@R=5pt{
[X_3,S^qBP]\ar[r]^{\cong\quad} &[S^0,S^qBP]=0\\
[X_3,S^{pq}BP]\ar[r]^{\cong\quad} &[S^0,S^{pq}BP]=0.}$$
So we will certainly be able to form the Toda brackets 
$$\{d_2,d_1, d_0\}, \quad\{d_1,d_0, h\}.$$
The first is a set of maps $SC_0\to C_3$ but all such maps are null-homotopic because $BP^*(BP)=0$ in dimensions $\equiv -1$ mod $q$. Similarly, the second is a set of maps
$SX_3\to C_2$; but all such maps are null-homotopic because $BP^*(X_3)=0$
in dimensions $\equiv -1$ mod $q$. So these brackets are zero mod zero. Therefore the quadruple Toda bracket
$$\{d_2,d_1,d_0,h\}$$
is defined, and it is a set of maps $S^2X_3\to C_3$.

Similarly, in the sphere $S^0$, $\{d_2, d_1, d_0, hi\}$ is defined,
and it  is a set of maps
$S^2\longrightarrow S^{(2p+2)q}BP$, i.e., it is zero 
mod zero.

Let us consider the indeterminacy of $\{d_2, d_1, d_0,h\}$.
We can vary the null-homotopies $d_2d_1\simeq 0$, $d_1d_0\simeq 0$,
$d_0h\simeq 0$ by maps 
$$SC_1\longrightarrow C_3, \quad SC_0\longrightarrow C_2,
\quad X_3\longrightarrow C_1.$$
All such maps are null-homotopic, the first two because $BP^*(BP)=0$
in dimensions $\equiv -1$ mod $q$, the third because $BP^*(X_3)=0$ in 
dimensions $\equiv -1$ mod $q$. We can vary the null-homotopies
$$\{d_2,d_1, d_0\}\simeq 0, \quad\{d_1,d_0, h\}\simeq 0$$
by maps 
$$S^2C_0\longrightarrow C_3,\quad S^2X_3\longrightarrow C_2.$$
Any map of the first sort is null-homotopic because $BP^*(BP)=0$ in dimensions
$\equiv -2$ mod $q$. Maps of the second sort may be non-zero, but
$$\xymatrix{
[S^r,C_2]\ar[rr]^{j^*} &&[S^2X_3,C_2]}$$
is an isomorphism. We conclude that the indeterminacy of our operation
is precisely
$$d_2[S^r,C_2]j.$$

I now seek to prove that this indeterminacy is contained in the set $(p,v_1)\ell$.
First I observe that the ideal $(p,v_1)\subseteq\pi_*(BP)$ is closed under all primary operations. This is clear, because it is the
kernel of an induced homomorphism
$$\pi_*(BP)=BP^*(S^{2p})\longrightarrow BP^*(V(1))=\frac{\pi_*(BP)}{(p,v_1)}.$$
Next we observe that any element of $[S^r, C_2]$ has the form
$$\left[\begin{array}{c}a\\ b\end{array}\right]$$
with $a, b\in\pi_*(BP)$ of degrees $(p^2 - p- 3)q$, $(p^2 - 2p - 2)q$. So they
must be of the form $a =v_1^pc$, $b=v_1d$ with $c$, $d$ of degree $(p-3)(p+1)q$.
So they lie in $(p,v_1)$, and so
$$d_2\left[\begin{array}{c}a\\ b\end{array}\right]=
\left[\begin{array}{lr}R_p &R_1\end{array}\right]
\left[\begin{array}{c}a\\ b\end{array}\right]$$
lies in $(p, v_1)$. It remains only to evaluate the operation, and this will take a bit more work.

First perhaps, you would like me to give you a bit more information on the action of operations $R_I$ on the coefficient ring
$\pi_*(BP)$. We have
$$R_I(xy)=\sum_{J+K=I}(R_Jx)(R_Ky)$$
so it is sufficient to give $R_I$ on the generators.

\vspace{10pt}

\noindent
{\bf Lemma 7.3.} \textit{Let the generators $v_i$ be defined by the formulae of Hazewinkel. Then
\begin{eqnarray*}
R_0v_1 &= &v_1\\
R_1v_1 &= &p\\
R_Iv_i &= &0 \mbox{ for } |I|> 1\\
R_0v_2 &= &v_2\\
R_1v_2 &= &-(p+1)v_1^p\\
R_iv_2 &\equiv &0\mod p^iv_i^{p+1-i}\mbox{ for }1<i<p\\
\end{eqnarray*}
\begin{eqnarray*}
R_pv_2 &\equiv &v_1\mod\,  p^{p-1}v_1\\
R_{p+1}v_2 &\equiv &0\mod\, p^p\\
R_{01}v_2 &= &p\\
R_Iv_2 &= &0\mbox{ for }|I|>p+1\\
R_1v_3 &= &-v_2^p\mod\, (p,v_1)\\
R_pv_3 &\equiv &0\mod\, (p,v_1)
\end{eqnarray*}
}

\vspace{10pt}

\noindent
{\bf Proof.} The formulae of Hazewinkel are as follows. We write
$m_i$ for what is called $m_{p^i-1}$ in \cite{A2}, i.e.,
$m_i=\frac{[CP^{p^i-1}]}{p^i}$. Then
\begin{itemize}
\item[(i)] $v_1=pm_1$
\item[(ii)] $v_2=pm_2 - v_1^pm_1$
\item[(iii)] $v_3=pm_3 -v_1^{p^2}m_2-v_2^pm_1$
\end{itemize}

Perhaps we should try to justify these formulae. First consider the homomorphism
$$Q_{2n}(\pi_*(BP))\longrightarrow Q_{2n}(H_*(BP))$$
where $Q_{2n}$ means the indecomposable quotient in dimension $2n$. It is known that its image is the subgroup of index $p$ if $n=p^i-1$ (\cite[Lemma 8.10, p. 58]{A2}). So we must have
$$v_i=\lambda pm_i\mbox{ mod decomposables}$$
where $\lambda=\frac{a}{b}$ with $a$ and $b$ prime to $p$, and we may as well normalise by taking $\lambda=1$. This disposes of $v_1$. For $v_2$ we must have
$$v_2=pm_2+\lambda v_1^{p+1}$$
where $\lambda$ is rational; and since $p^2m_2=[CP^{p^2-1}]$ is integral, we must
have $p\lambda\in\mathbf{Z}_{(p)}$. Of course, we get an equally good generator by changing
$\lambda$ mod $p$. So we may as well set
$$v_2 = pm_2 + \mu v_1^pm_1$$
where $\mu\in\mathbf{Z}_{(p)}$ and we want to determine $\mu$ mod $p$.

Consider now the Kronecker product $\langle R_p\otimes R_{p^2-p},\psi t_3\rangle$ .
 It must lie in $\pi_*(BP)$. But according to \cite[Thm 16.1, p. 98]{A2}, we have
$$\psi t_3+m_1(\psi t_2)^p+m_2(\psi t_1)^{p^2}\qquad\qquad\qquad\qquad\qquad\qquad\qquad\qquad\qquad$$
$$\qquad\qquad\qquad =t_3\otimes 1+t_2\otimes t_1^{p^2}+t_1\otimes t_2^p
+1\otimes t_3+m_1t_2^p\otimes 1+m_1t_1^p\otimes t_1^{p^2}$$
$$\qquad\qquad\qquad +m_11\otimes t_2^p+m_2t_1^{p^2}\otimes 1+m_21\otimes t_1^{p^2}.$$

Consider the pairing of these various terms with $R_p\otimes R_{p^2-p}$. On the
right-hand side we get 0. On the left-hand side,  the term $m_2(\psi t_1)^{p^2}$
gives
$$m_2\frac{p^2)!}{p!(p^2-p)!}\equiv pm_2\mbox{ mod }p^3m_2.$$
We have
$$\psi t_2=t_2\otimes 1+t_1\otimes t_1^p+1\otimes t_2-v_1(t_1\otimes t_1^{p-1}+\cdots
+t_1^{p-1}\otimes t_i).$$
Therefore
$$\langle(\psi t_2)^p, R_p\otimes R_{p^2-p}\rangle\ =\ (-1)^pv_1^p=-v_1^p.$$
We conclude that $pm_2-m_1v_1^p$ lies in $\pi_*(BP)$. So we may as well
choose it for $v_2$. 

We now conclude that the coefficient of $t_1^p\otimes t_1^{p^2-p}$ in
$\psi t_3$ is $-v_2$ mod $p^3m_2$. For $0\le i < p$ the coefficient of
$t_1^i\otimes t_1^{p^2-1}$ in $\psi t_3$ is 0 mod $p^2m_2$. 

Let us now proceed similarly with the Kronecker product 
$\langle R_{p^2}\otimes R_{p^3-p^2},\psi t_4\rangle$. The result is an element of $\pi_*(BP)$.
By the  same reference we have 
$$\psi t_4+m_1(\psi t_3)^p+m_2(\psi t_2)^{p^2}+m_3(\psi t_1)^{p^3}=\cdots,$$
where every term on the right-hand side yields 0 when paired with 
${R_{p^2}\otimes R_{p^3-p^2}}$. The term $m_3(\psi t_1)^{p^3}$ gives
$$m_3\frac{(p^3)!}{(p^2)!(p^3-p^2)!}\equiv pm_3\mbox{ mod } p^3m_3.$$
The term $m_2(\psi t_2)^{p^2}$ gives $m_2(-v_1)^{p^2}=-m_2v_1^{p^2}$.
The term $m_1(\psi t_3)^p$ gives
$$m_1\left((-v_2)^p\mbox{ mod }(p^2m_2)\right).$$ 
We conclude that there is an element of $\pi_*(BP)$ of the form 
$$pm_3 - m_2v_1^{p^2}- m_1a$$
where $a\equiv v_2^p\mbox{ mod}(p^2m_2)$.

By the same reference, the effect of the operation $R_I$ on the 
generators $m_i$ is
$$R_0m_1=m_1, R_1m_1=1, R_Im_1=0 \mbox{ otherwise}$$  
$$R_0m_2=m_2, R_pm_2=m_1, R_{01}m_2=1,  R_Im_2=0 \mbox{ otherwise}$$
$$R_0m_3=m_3, R_{p^2}m_3=m_2, R_{0p}m_3=m_1, R_{001}m_3=1, R_Im_3=0 \mbox{ otherwise.}$$
So we get 
$$R_1v_1=R_1(pm_1) = pR_1(m_1)=p$$
$$R_1v_2=R_1\left(pm_2-\frac{v_1^{p+1}}{p}\right)
=-\frac{p+1}{p}(R_1v_1)v_1^p=-(p+1)v_1^p$$
$$R_iv_2=R_i\left(pm_2-\frac{v_1^{p+1}}{p}\right)
=-\frac{(p+1)!}{i!(p+1-1)!}\frac{1}{p}p^iv_1^{p+1-i}\quad(i\ne p)$$
and for $1 <i < p$ the binomial coefficient contains a factor $p$.
$$R_pv_2=R_p\left(pm_2-\frac{v_1^{p+1}}{p}\right)
=pm_1-(p+1)\frac{1}{p}p^pv_1=v_1\mbox{ mod }p^{p-1}v_1$$
and
$$R_{01}v_2=R_{01}\left(pm_2-\frac{v_1^{p+1}}{p}\right)=p\cdot 1.$$
Now take
$$v_3=pm_3-m_2v_1^{p^2}-m_1a$$
where
$$a =v_2^p+(p^2m_2)b.$$
Apply $R_1$. We get
\begin{eqnarray*}
R_1v_3 &= &-m_2p^2pv_1^{p^2-1}-a-m_1R_1a\\
&= &-(m_2p^2)pv_1^{p^2-1}-v_2^p-(p^2m_2)b-m_1(R_1v_2^p)-m_1(R_1p^2m_2)b\\
&&-m_1(p^2m_2)R_1b\\
&=&-\left(pv_2+v_1^{p+1}\right)pv_1^{p^2-1}-v_2^p-\left(pv_2+v_1^{p+1}\right)b
+m_1p(p+1)v_1v_2^{p-1}\\
&&-0-m_1\left(pv_2+v_1^{p+1}\right)R_1b.
\end{eqnarray*}
Now $R_1v_3$ is certainly integral, and all the terms on the right are 
integrai except possibly $-m_1v_1^{p+1}R_1b$.
We conclude that $R_1b$ must be divisible by $p$, and the last term is
$v_1(pv_2+v_1^{p+1})c$. Thus
$$R_1v_3\equiv -v_2^p\mbox{ mod }(p,v_1^{p+1}).$$
For reasons of degree, we have
$$R_pv_3\in(v_1^2)$$
 This proves Lemma 7.3.\hfill$\Box$

\vspace{10pt}

Now I have to explain how I compute the fourfold bracket. Let me begin by explaining it in a different way: $\gamma_1$ is constructed as the following composite
(where $r = (p^2-1)q$).
$$\xymatrix{
S^{r-3}V(2)\ar[rrddd]^D\\
S^{r-3}V(1)\ar[u]^{i_2}\\
S^{r-3}V(0)\ar[u]^{i_1}\\
S^{r-3}\ar[u]^{i_0}\ar[rr]^{\gamma_1} &&S^0
}$$
Here $D$ induces a non-zero homomorphism of $BP^*$, that is, it carries the generator to $v_3g_2$. The other three maps induce the zero map of $BP^*$, 
but they give known elements in $\Ext^1_{BP^*(BP)}$. Therefore, the
composite $\gamma_1$ corresponds to a calculable element of $\Ext^3_{BP^*(BP)}$,
namely the Yoneda product.

Now I don't propose to do it exactly that way, but in a parallel way. That is, using a relation of the Peterson-Stein type, I will reduce
the calculation of our tertiary operation in $X_3=S^0\cup_{\gamma_1}CS^{r-3}$
to the calculation of a secondary operation in $X_2 = S^0\cup CS^{r-3}V(0)$. (The attaching map gives an element of $\Ext^2$.) Using another relation of the Peterson-Stein type, I will reduce the calculation of the secondary operation
in $X_2$ to the calculation of a primary operation in  $X_1 = S^0\cup CS^{r-3}V(1)$.
(The attaching map then gives an element of $\Ext^1$.) Finally, I calculate this primary operation using the known behaviour of $i_2$, $D$.

In practice, of course, I proceed from the known to the unknown, and calculate first in $X_1$, then in $X_2$, then in $X_3$.

In calculating with composite maps, we always use the following lemma.

\vspace{10pt}

\noindent	
{\bf Lemma 7.4.} (Verdier's Axiom) \textit{Given a composition
$$\xymatrix{X\ar[r]^f &Y\ar[r]^g &Z},$$
we can form the following diagram of cofiberings
$$\xymatrix{
Y\ar[dd]\ar[rr] &&Y\cup_f CX\ar[dr]\ar[dl]\\
&Z\cup_{gf}CX\ar[dd]\ar[rr] &&SX\ar[dd]\ar[dr]\\
Z\ar[ur]\ar[dr] &&&&SZ\\
&Z\cup_g CY\ar[rr]\ar[dr] &&SY\ar[dl]\ar[ur]\\
&&S(Y\cup_f CX)
}$$}

\vspace{10pt}

\noindent
{\bf Lemma 7.5.} \textit{Let $X_1=S^0\cup_{Di_2}CS^{r-3}V(1)$. Then
$$BP^*(X_1)\cong\pi_*(BP)\mbox{ on a unique generator $h_1$ of degree 0}$$
$$\qquad\qquad\qquad\qquad\mbox{restricting to the generator on $S^0$}$$
$$\qquad\qquad\qquad\qquad\oplus\frac{\pi_*(BP)}{(p,v_1)}\mbox{ on a generator $g_1$ of
degree $p^2q$ coming}$$
$$\qquad\qquad\qquad\qquad\qquad\mbox{from the generator on $S^{r-2}V(1)$.}$$
We have
$$d_0h_1=\left[\begin{array}{c}-v_2^{p-1}\\ 0\end{array}\right].$$
}

\vspace{10pt}

\noindent
{\bf Proof.} Take the diagram of Lemma 7.4 and make the following substitutions
$$\xymatrix{
X\ar[r]^f\ar@{=}[d] &Y\ar[r]^g\ar@{=}[d] &Z\ar@{=}[d]\\
S^{r-3}V(1)\ar[r]^{i_2} &S^{r-3}V(2)\ar[r]^{\quad D} &S^0
}$$
We obtain:
$$\xymatrix{
\qquad\parbox{100pt}{\quad$\frac{\pi_*(BP)}{(p,v_1,v_2)}$\\
on generator\\
$g_2$ of degree\\
$(p^2+p+1)q$
}\ar@{=}[d]
&&\qquad\parbox{100pt}{
\quad$\frac{\pi_*(BP)}{(p,v_1}$\\
on generator\\
 $\overline{g}_1$ of degree\\
$(p^2+p+1)q$
}\ar@{=}[d]\\
BP^*\left(S^{r-3}V(2)\right) &&BP^*\left(S^{(p^2+p)q-2}V(1)\right)\ar[ll]\\
&BP^*(X_1)\ar[ur]_{i^*}\ar[dl] &BP^*\left(S^{r-2}V(1)\right)\ar[u]\ar[l]\\
BP^*(S^0)\ar[uu]_{D^*}
&&\qquad\parbox{100pt}{\begin{center}$\frac{\pi_*(BP)}{(p,v_1)}$\\
on generator\\
$\tilde{g}_1$ of degree\\
$p^2q$\end{center}}\ar@{=}[u]\\
\parbox{100pt}{\begin{center}
$\pi_*(BP)$\\
on generator\\
of degree 0\end{center}}\ar@{=}[u]
&&BP^*\left(S^{r-2}V(2)\right)\ar[u]
}$$
Since $i_2$ induces the zero map of $BP^*(-)$, $Di_2$ must also do so, and the
exact sequence across the middle is short exact. Now it is clear that $BP^*(X_1)$
is as stated, except that we have to show the choice of $h_1$ is unique. Since
there is no monomial $v_2^iv_3^j\cdots$ in $\frac{\pi_*(BP)}{(p,v_1)}$	of the
relevant degree $p^2q$, this
is clear.
The generator $g_1$ in $BP^*(X_1)$ comes from 
$BP^*\left(S^{r-2}V(1)\right)$, and since the right-hand column is the cofibering which constructs $V(2)$ we have
$$g_1i=v_2\overline{g}_1.$$
Also $h_1$ maps to the generator in $BP^*(S^0)$ and then to $v_3$ in 
$BP^*\left(S^{r-3}V(2)\right)$; so
$$h_1i = v_3\overline{g}_1\mbox{ mod }(p, v_1, v_2).$$
But since $(p^2 +p+1)q$ is not divisible by $p+1$, no power $v_2^n\overline{g}_1$
lies in degree 0. Therefore
$$h_1i=v_3\overline{g}_1\mbox{ mod }(p, v_1 ),$$
i.e.,
$$h_1i=v_3\overline{g}_1.$$

Now in $S^0$ we have
$$\left[\begin{array}{c}R_1\\ R_p\end{array}\right]h_1=0.$$
So in $X_1$
$$\left[\begin{array}{c}R_1\\ R_p\end{array}\right]h_1$$
must be a multiple of $g_1$. But on the summand
$\frac{\pi_*(BP)}{(p,v_1)}g_1$, $i^*$ is a
 monomorphism; so it is sufficient to calculate in
$$S^{(p^2+p)q-2}V(1).$$
We find
\begin{eqnarray*}
\left[\begin{array}{c}R_1\\ R_p\end{array}\right]h_1i
&= &\left[\begin{array}{c}R_1\\ R_p\end{array}\right]v_3\overline{g}_1=
\left[\begin{array}{c}-v_2^p\\ 0\end{array}\right]\overline{g}_1
\mbox{ mod }(p, v_1)\mbox{ (by 7.3)}\\
&= &\left[\begin{array}{c}-v_2^p\\ 0\end{array}\right]g_1i.
\end{eqnarray*}
Therefore
$$\left[\begin{array}{c}R_1\\ R_p\end{array}\right]h_1=
\left[\begin{array}{c}-v_2^p\\ 0\end{array}\right]g_1.$$
\hfill$\Box$

To get any further we need Peterson-Stein relations. For the first
one, we suppose again we have maps
$$\xymatrix{
C_0\ar[r]^{d_0} &C_1\ar[r]^{d_1} &C_2}$$
with $d_1d_0\simeq 0$. Here $d_0$ and $d_1$ are in practice the maps specified above, but the precise spectra are not relevant. We interest ourselves in secondary operations of the form
$$\phi(\epsilon)=\left\{d_1,d_0,\epsilon\right\}$$
where $\epsilon:Y\to C_0$. We suppose we are given a cofibering
$$\xymatrix{
X\ar[r]^i &Y\ar[r]^j &Z\ar[r]^k\ar[d]^\zeta &SX\ar[r]^{Si}\ar[d]^\xi &SY\\
&&C_0\ar[r]^{d_0} &C_1
}$$
and maps $\zeta:Z\to C_0$, $\xi:SX\to C_1$ such that $\xi k=d_0\zeta$. Given $\zeta$,
two such $\xi$ differ by an element of
$$[SY,C_1]Si.$$
We have $d_0\zeta j\simeq\xi kj\simeq 0$ so
$$\{d_1,d_0,\zeta j\}\in \frac{[SY,C_2]}{d_1[SY,C_1]+[SC_0,C_2](S\zeta)(Sj)}$$
and
$$\{d_1,d_0,\zeta j\}Si\in\frac{[SX,C_2]}{d_1[SY,C_1]Si}.$$

\vspace{10pt}

\noindent
{\bf Lemma 7.6.} $\{d_1,d_0,\zeta j\}Si=-d_1\xi\mbox{ mod }d_1[SY,C_1]Si$.

\vspace{10pt}

\noindent
{\bf Proof.} This is simply the familiar identity
$$\{d_1,d_0,\zeta j\}Si = -d_1\{d_0,\zeta i j,i\}\supseteq -d_1\{d_0\zeta,j,i\}.$$
(The latter is actually an equality sign because the indeterminacies are the same.) In fact, the prescription for $\{d_0\zeta,j, i\}$ is that you take
$$\xymatrix{
SX\ar[r]^{\underline{i}} &Z\cup_j CY\ar@{=}[d]\ar[r]^{\overline{d_0\zeta}} &C_1\\
&(Y\cup_i CX)\cup CY\ar[d]^\simeq\\
&SX}$$
and the identity map qualities as $\underline{i}$.\hfill$\Box$

\vspace{10pt}

\noindent
{\bf Lemma 7.7.} \textit{Let $X_2=S^0\cup_{Di_2i_1}CS^{r-3}V(0)$. Then
$$BP^*(X_2)\cong\pi_*(BP)\mbox{ on a generator $h_2$ of degree 0 restricting}$$
$$\qquad\qquad\mbox{to the generator on $S^0$}$$
$$\quad\oplus\frac{\pi_*(BP)}{(p)}\mbox{ on a generator $g_0$ coming}$$
$$\qquad\qquad\qquad\qquad\mbox{from the generator in $S^{r-2}V(1)$}.$$
The operation $\{d_1, d_0, h_2\}$ is defined and it is a coset of maps
$SX_2\to C_2$ and it is
$$\left[\begin{array}{c}2v_1^pv_2^{p-3}g_0\\ 2v_1v_2^{p-3}g_0\end{array}\right]
\mod\, d_1[SX_2,C_1].$$
}

\vspace{10pt}

\noindent
{\bf N.B.} I want to leave the indeterminacy in this form because will disappear at the next step.

\vspace{10pt}

\noindent
{\bf Exercise.} Show directly that the value given does not lie in the indeterminacy.

\vspace{10pt}

\noindent
{\bf Proof of 7.7.} We consider $Di_2i_1$ as the composite $(Di_2)i_1$. We take the diagram of Lemma 7.4 and substitute
$$\xymatrix{
X\ar[r]^f\ar@{=}[d] &Y\ar[r]^g\ar@{=}[d] &Z\ar@{=}[d]\\
S^{r-3}V(0)\ar[r]^{i_1} &S^{r-3}V(1)\ar[r]^{\quad Di_2} &S^0
}$$
We get
$$\xymatrix@C=5pt@R=18pt{
&&\parbox{100pt}{\begin{center}
$\frac{\pi_*(BP)}{(p)}$\\
on generator $\overline{g}_0$\\
of degree $p^2q-1$
\end{center}}\ar@{=}[d]\\
&&BP^*\left(S^{p^2q-2}V(0)\right)\\
&BP^*(X_2)\ar[ur]^{i_*}&&BP^*\left(S^{r-2}V(0)\right)\ar[ul]\ar[ll]\\
&\parbox{105pt}{\begin{center}
$\pi_*(BP)\oplus\frac{\pi_*(BP)}{(p)}$\\
on generators $h_2$, $g_0$\\
of degrees 0, $r-1$\end{center}}\ar@{=}[u]\ar[dl]
&&\parbox{100pt}{\begin{center}
$\frac{\pi_*(BP)}{(p)}$\\
on generator $g_0$\\
of degree $r-1$\end{center}}\ar@{=}[u]\\
BP^*(S^0) &&&&BP^*(S^1)\ar[ul]\ar[dl]\\
&\parbox{105pt}{\begin{center}
$\pi_*(BP)\oplus\frac{\pi_*(BP)}{(p,v_1)}$\\
on generators $h_1$, $g_1$\\
of degrees 0, $p^2q$\end{center}}\ar[uu]_{j^*}\ar[ul]
&&\parbox{100pt}{\begin{center}
$\frac{\pi_*(BP)}{(p,v_1)}$\\
on generator $g_1$\\
of degree $p^2q$\end{center}}\ar[uu]\\
&BP^*(X_1)\ar@{=}[u] &&BP^*\left(S^{r-2}V(1)\right)\ar@{=}[u]\ar[ll]\\
&&BP^*\left(S^{p^2q-1}V(0)\right)\ar[ul]^{k^*}\ar[ur]\\
&&\parbox{100pt}{\begin{center}
$\frac{\pi_*(BP)}{(p)}$\\
on generator $\overline{g}_0$\\
of degree $p^2q$\end{center}}\ar@{=}[u]
}$$
Certainly the map of $BP^*$ induced by $Di_2i_1$ is zero, for those induced by
$i_2$ and $i_1$ are zero. So we have a short exact sequence for $BP^*(X_2)$.
We may take the generator $h_1$ in $BP^*(X_1)$ and map it into $X_2$; this gives an element restricting to the generator in $BP^*(S^0)$; we call this element $h_2$. Now it is clear that the structure of $BP^*(X_2)$ is as stated.

The generator $\overline{g}_0$ at the bottom of the diagram maps to $g_1$
in $S^{r-2}V(1)$ and to $g_1$ in $X_1$. That is, $\overline{g}_0h_1=g_1$

By 7.5, the map
$$X_1\stackrel{d_0h_1}{\longrightarrow}C_1$$
is
$$\left[\begin{array}{c}-v_2^{p-1}g_1\\ 0\end{array}\right],$$
and it is now clear that it factors through $k$; in fact, it is $\xi k$, where
$$\xi=\left[\begin{array}{c}-v_2^{p-1}\overline{g}_0\\ 0\end{array}\right].$$
So
$$d_0h_2=d_0h_1 j = \xi kj = 0,$$
and the bracket $\{d_1,d_0,h_2\}$ is defined. We must check its indeterminacy.
Since there is no way to change the homotopy $d_1d_0\simeq 0$ the indeterminacy
is $d_1[SX_2,C_1]$.

We now proceed to apply Lemma 7.6 taking the cofibering in that
lemma to be the one marked with $i^*,j^*,k^*$ in the previous diagram. We
take $\zeta=h_1$ and $\xi$ as above. The lemma gives
\begin{eqnarray*}
\{d_1,d_0,h_2\}Si &= &-d_1\xi\mbox{ mod }d_1[SX_2,C_1]Si\\
&= &\left[\begin{array}{ll}
R_pR_1-2R_1R_p &R_1^2\\
R_p^2 &-2R_pR_1+R_1R_p\end{array}\right]
\left[\begin{array}{c}-v_2^{p-1}\overline{g}_0\\ 0\end{array}\right]\\
&&\qquad\qquad\qquad\qquad\qquad\mbox{ mod }d_1[SX_2,C_1]Si.
\end{eqnarray*}
Now
$$R_pR_1-2R_1R_p=-R_{01}-R_1R_p,$$
and
$$R_{01}v_2^{p-1}\overline{g}_0= 0\mbox{ mod }p$$
$$R_pv_2^{p-1}\overline{g}_0=(p-1)v_1v_2^{p-2}\overline{g}_0
=-v_1v_2^{p-2}\overline{g}_0\mbox{ mod }p$$
\begin{eqnarray*}
-R_1R_pv_2^{p-1}\overline{g}_0&= &R_1(v_1v_2^{p-2}\overline{g}_0)\mbox{ mod }p\\
&= &v_1(p-2)(-1)v_1^pv_2^{p-3}\overline{g}_0\mbox{ mod }p\\
&= &2v_1^{p+1}v_2^{p-3}\overline{g}_0\mbox{ mod }p.
\end{eqnarray*}
Similarly
\begin{eqnarray*}
R_pv_2^{p-1}\overline{g}_0 &= &-v_1v_2^{p-2}\mbox{ mod }p\mbox{ (as above)}\\
R_p^2(v_2^{p-1}\overline{g}_0) &= &-v_1(p-2)v_1v_2^{p-3}\overline{g}_0\mbox{ mod }p\\
&= &2v_1^2v_2^{p-3}\overline{g}_0\qquad\mbox{ mod }p.
\end{eqnarray*}
Hence we obtain
$$\{d_1,d_0,h_2\}Si=\left[
\begin{array}{c}
2v_1^{p+1}v_2^{p-3}\overline{g}_0\\
2v_1^2v_2^{p-3}\overline{g}_0
\end{array}\right]\mbox{ mod }d_1[SX_2,C_1]Si.$$
Now from the commutativity of the upper triangle we see that $g_0i=\overline{g}_0$;
so we obtain
$$\{d_1,d_0,h_2\}Si=\left[
\begin{array}{c}
2v_1^{p}v_2^{p-3}g_0\\
2v_1v_2^{p-3}g_0
\end{array}\right]\mbox{ mod }d_1[SX_2,C_1]Si.$$
The injection $S^0\to X_2$
clearly annihilates $\{d_1,d_0,h_2\}$
(maps $S^1\to C_2$ are
zero), so it annihilates the indeterminacy $d_1[SX_2,C_1]$, and it clearly
annihilates
$$\left[
\begin{array}{c}
2v_1^{p}v_2^{p-3}g_0\\
2v_1v_2^{p-3}g_0
\end{array}\right].$$
But on the complementary summand $\frac{\pi_*(BP)}{(p)}$,\ \  $(Si)^*$ is mono, so
$$\{d_1,d_0,h_2\}=\left[
\begin{array}{c}
2v_1^{p}v_2^{p-3}g_0\\
2v_1v_2^{p-3}g_0
\end{array}\right]\mbox{ mod }d_1[SX_2,C_1].$$
This proves 7.7.\hfill$\Box$

\newpage

For the second Peterson-Stein relation, suppose we are given maps
$$\xymatrix{
C_0\ar[r]^{d_0} &C_1\ar[r]^{d_1} &C_2\ar[r]^{d_2} &C_3}$$
with $d_2d_1\simeq 0$, $d_1d_0\simeq 0$ (e.g. as above). We interest ourselves in the
tertiary operations of the form $\{d_2,d_1,d_0,\epsilon\}$ where $\epsilon:Y\to C_0$, 
so we will have to give ourselves the sort of data which ensure that such operations are defined.

We suppose given 0. cofibering
$$\xymatrix{
X\ar[r]^i &Y\ar[r]^j &Z\ar[r]^k &SX\ar[r]^{Si} &SY}$$
and two maps
$$\xymatrix{
S^2X\ar[d]^\xi &&Z\ar[d]^\zeta\\
C_1 &&C_0}$$
We suppose $d_0\zeta\simeq 0$. We also pick a homotopy $d_1d_0\simeq 0$;
in the applications there is only one way to pick it. We suppose that for this homotopy we have simultaneously
\begin{itemize}
\item[(i)] $\{d_2,d_1,d_0\}=0$ mod $[SC_1,C_3]Sd_0$
\item[(ii)] $\{d_1,d_0,\zeta\}=0$ in $[SZ,C_2]/d_1[SZ,C_1]$.
\end{itemize}
We pick a homotopy $d_2d_1\simeq 0$ which gives 0 as a representative map for
$\{d_2, d_1, d_0\}$. After this we are only willing to change it by a map
$\theta\in[SC_1, C_3]$ such that $\theta Sd_0\simeq 0$. We also pick a homotopy
$d_0\zeta\simeq 0$
which gives $\zeta Sk$ as the representative map for $\{d_1,d_0,\zeta\}$. After this we are willing to vary this homotopy by an element $\phi\in [SZ,C_1]$ such
that $d_1\phi\simeq 0$.

Then $\{d_1,d_0,\zeta j\}$ is represented by $\xi\cdot Sk\cdot Sj=0$, and we can form
\begin{eqnarray*}
\{d_2,d_1,d_0,\zeta\} &\in 
&[S^2Y,C_3]/\left(d_2[S^2C_0,C_2]+[S^2C_0,C_2]S^2\zeta\cdot S^2j\right.\\
&&\qquad\qquad+\{\{d_2,d_1,\phi Sj\}|\phi\in[SZ,C_1],d_1\phi\simeq 0\}\\
&&\qquad\qquad+\{\{\theta,Sd_0,S\zeta\cdot Sj\}|\theta\in[SC_1,C_3],\theta Sd_0\simeq 0\}\\
&&\qquad\qquad+\left.\{\theta'\cdot S\phi'\cdot S^2j|\theta'Sd_0\simeq 0,d_1\phi'\simeq 0\}\right)
\end{eqnarray*}
(here we must also include any indeterminacy implied by the Toda brackets which is not already accounted for). Then we may as well write $\{d_2,d_1,\phi\}S^2j$ for
$\{d_2, d_1, \phi\cdot Sj\}$ and $\{\theta, Sd_0, S\zeta\}S^2j$ for
$\{\theta, Sd_0, S\zeta\cdot Sj\}$ because these are defined, are contained in what we wrote before, and the larger indeterminacy is contained in what we have already written.

In the applications, of course, the indeterminacy reduces to
$d_2[S^2Y,SC_2]$, as we have seen.

In general $\{d_2,d_1,d_0,\zeta j\}S^2i$ lies in $[S^2X,C_3]/d_2[S^2Y,SC_2]S^2i$.

\vspace{10pt}

\noindent
{\bf Lemma 7.8.}
$$\{d_2,d_1,d_0,\zeta j\}S^2i=-d_2\xi:S^2X\to C_3\mbox{ mod }d_2[S^2Y,SC_2]S^2i.$$

\vspace{20pt}

\noindent
{\bf N.B.} We do not allow ourselves to vary $\xi$, because the indeterminacy of the quadruple bracket when we vary everything in sight may be greater than what we have written.

\vspace{20pt}

\noindent
{\bf Proof.} In Lemma 7.6 we substitute\\
\begin{tabular}{lllll}
&\qquad\qquad&$C_2\mapsto C_3$ &&$d_1\mapsto d_2$\\
&\qquad\qquad&$C_1\mapsto C_2$ &&$d_0\mapsto \overline{d}_1$\\
&\qquad\qquad&$C_0\mapsto C_1\cup_{d_0}CC_0$ &&(the choice of null-homotopy $d_1d_0\simeq 0$\\
&\qquad\qquad&&&\ gives a specific extension of $d_1$)
\end{tabular}\\
Now the old condition $d_1d_0\simeq 0$ becomes the two new conditions
$d_2d_1\simeq 0$ and $\{d_2,d_1,d_0\}=0$ mod $[SC_1,C_3]Sd_0$.
For the cofibering in 7.6, we substitute
$$\xymatrix{
SX\ar[r]^{Si} &SY\ar[r]^{Sj} &SZ\ar[r]^{Sk}\ar[d]^{\underline{\zeta}}
&S^2X\ar[r]^{S^2i}\ar[d]^{\xi} &S^2Y\\
&&C_1\cup_{d_0}CC_0\ar[r]^{\qquad d_1} &C_2
}$$
and for the maps, $\xi$ and a coextension $\underline{\zeta}$ of 
$\zeta:Z\to C_0$. Here the composite
$$\xymatrix{
SZ\ar[r]^{\underline{\zeta}\qquad} &C_1\cup_{d_0}CC_0\ar[r]^{\qquad d_1} &C_1}$$
of course represents $\{d_1,d_0,\zeta\}$, and
we choose the coextension $\underline{\zeta}$ so that we actually get the element
$\xi\cdot Sk$ in $[SZ,C_2]$.
Then the diagram is commutative. The triple bracket
$\{d_1,d_0,\underline{\zeta}\cdot Sj\}$  is now an alternative construction
for $\{d_2,d_1,d_0,\zeta j\}$ and
if you don't like quadruple brackets, take this as the definition. Applying
Lemma 7.6, we come out with
$$\{d_2,d_1,d_0,\zeta j\}=-d_2\xi\mbox{ mod }d_2[S^2Y,C_2]S^2i.$$
\hfill$\Box$

\vspace{10pt}

\noindent
{\bf Proof of Theorem 7. 2 completed.} We continue the process exhibited in 7.5, 7.7 applying Lemma 7.8. We consider	$\gamma_1=Di_2i_1i_0$ as the composite $(Di_2i_1)i_0$.
We take the diagram of Lemma 7.4, and substitute
$$\xymatrix{
X\ar[r]^f\ar@{=}[d] &Y\ar[r]^g\ar@{=}[d] &Z\ar@{=}[d]\\
S^{r-3}\ar[r]^{i_0\quad} &S^{r-3}V(0)\ar[r]^{\quad Di_2i_1} &S^0}$$
We obtain
$$
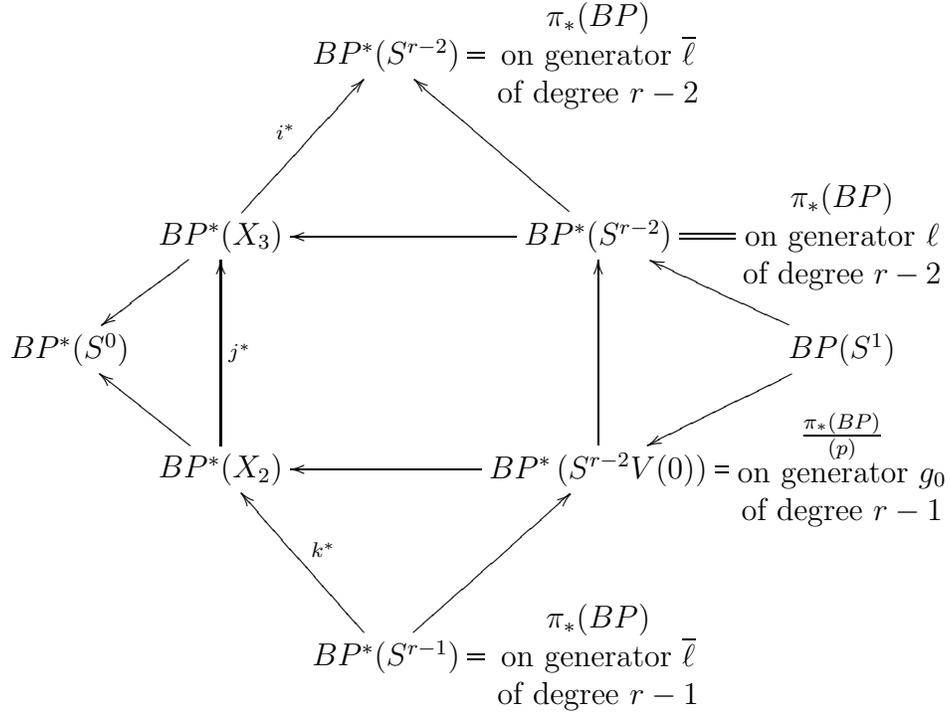
\xymatrix@C=5pt@R=0pt{
&&BP^*(S^{r-2})\ar@{=}[r] 
&\parbox{80pt}{\begin{center}
$\pi_*(BP)$\\
on generator $\overline{\ell}$\\
of degree $r-2$
\end{center}}\\
&BP^*(X_3)\ar[dl]\ar[ur]^{i^*} &&BP^*(S^{r-2})\ar[ul]\ar[ll]
&\parbox{75pt}{\begin{center}
$\pi_*(BP)$\\
on generator $\ell$\\
of degree $r-2$\end{center}}\ar@{=}[l]\\
BP^*(S^0) &&&&BP(S^1)\ar[ul]\ar[dl]\\
&BP^*(X_2)\ar[uu]_{j^*}\ar[ul] &&BP^*\left(S^{r-2}V(0)\right)\ar[uu]\ar[ll]
&\parbox{80pt}{\begin{center}
$\frac{\pi_*(BP)}{(p)}$\\
on generator $g_0$\\
of degree $r-1$\end{center}}\ar@{=}[l]\\
&&BP^*(S^{r-1})\ar[ul]_{k^*}\ar[ur]
&\parbox{80pt}{\begin{center}
$\pi_*(BP)$\\
on generator $\overline{\ell}$\\
of degree $r-1$
\end{center}}\ar@{=}[l]
}$$
We need some induced maps. The right-hand vertical column is induced by the
cofibering
$$S^0\stackrel{p}{\longrightarrow} S^0\longrightarrow V(0),$$
so $\overline{\ell}$ at bottom maps to $g_0$, and $\ell$ maps
to $p\overline{\ell}$. We now proceed to apply Lemma 7.8 to the cofibering marked
$i^*,j^*,k^*$.
We take	$\zeta h_0$ to be $h_2:X_2\to C_0$. That $d_0h_2=0$ was already proved in Lemma 7.7. We now wish to produce $\xi$.

Now Lemma 7.6 contains an indeterminacy $d_1[SX_2 C_1]$. But every
element of $[SX_2, C_1]$ lifts to $[S^{r-1}V(0, C_1]$ by 7.7;
and every element of $[S^{r-1}V(0, C_1]$ clearly lifts to $[S^r,C_1]$.
Thus every element of
$d_1[SX_2,C_1]$ lies in $d_1[S^r,C_1]Sk$. We can also lift the representative
$$\left[\begin{array}{c}
2v_1^pv_2^{p-3}g_0\\
2v_1v_2^{p-3}g_0
\end{array}\right]$$
by taking
$$\xi=\left[\begin{array}{c}
2v_1^pv_2^{p-3}\overline{\ell}\\
2v_1v_2^{p-3}\overline{\ell}
\end{array}\right].$$
So we see that we can in fact lift any element of $\{d_1,d_0,h_2\}$ by a lift
congruent to
$$\left[\begin{array}{c}
2v_1^pv_2^{p-3}\overline{\ell}\\
2v_1v_2^{p-3}\overline{\ell}
\end{array}\right]\mbox{ mod } d_1[S^R,C_1].$$
Now Lemma 7.8 states that 
$$\{d_2, d_1, d_0, h_3 \}S^2i = -d_2\xi:S^2X_3\to C_3\mbox{ mod }d_2[S^2X_3,C_2]S^2i.$$
Here $d_2[S^2X_3,C_2]$
is in our case exactly the indeterminacy of $\{d_2, d_1, d_0, h_3 \}$
and we have already shown that it is contained in $(p,v_1)\ell$. Therefore
$$d_2[S^2X_3,C_2]S^2i\subseteq (p^2,pv_1)\overline{\ell}.$$
As for $-d_2\xi$ consider first the term
$$-d_2d_1[S^r, C_1].$$
It is 0, because $d_2d_1\simeq 0$. Consider next the term
$$-\left[\begin{array}{ll}R_p &R_1\end{array}\right]\left[
\begin{array}{c}
2v_1^pv_2^{p-3}\overline{\ell}\\
2v_1v_2^{p-3}\overline{\ell}
\end{array}\right].$$
We have
$$R_p\left(v_1^pv_2^{p-3}\right)=\sum_{i+j=p}(R_iv_1^p)(R_jv_2^{p-3}).$$
Here the first term 
$$R_iv_1^p=\frac{p!}{i!(p-i)!}p^iv_i^{p-i}$$
lies in $(p^2,pv_1)$ unless $i=0$, in which case
$$R_p(v_2^{p-3})=(p-3)v_1v_2^{p-4}\mbox{ mod }p^2$$
and we get
$$R_p\left(v_1^pv_2^{p-3}\right)=(p-3)v_1^{p+1}v_2^{p-4}\mbox{ mod }(p^2,pv_1).$$
We have by 7.3
\begin{eqnarray*}
R_1(v_1v_2^{p-3}) &= &pv_2^{p-3}+v_1(p-3)(-(p+1))v_1^pv_2\\
&= &pv_2^{p-3}-(p-3)v_1^{p+1}v_2^{p-4}\mbox{ mod }pv_1^{p+1}.
\end{eqnarray*}
Totalling we get 
$$-\left[\begin{array}{ll}R_p &R_1\end{array}\right]\left[
\begin{array}{c}
2v_1^pv_2^{p-3}\overline{\ell}\\
2v_1v_2^{p-3}\overline{\ell}
\end{array}\right]=-2pv_2^{p-3}\overline{\ell}\mbox{ mod }(p^2,pv_1).$$
Since $\ell(S^2i) = p\overline{\ell}$, we conclude that 
$$\{d_2,d_1,d_0,h_3\}=-2v_2^{p-3}\ell\mbox{ mod }(p,v_1).$$
This completes the proof, and shows that $\gamma_1\ne 0$.\hfill$\Box$

\vspace{20pt}

Oka and Toda say they have not checked the assertion of Larry 
Smith that $\beta_p\ne 0$ and if this is false their proof collapses. So perhaps 
it will be as well to examine $\beta_p$. It turns out that by using an indirect
method, we can avoid even mentioning secondary operations; we need some
preliminary information first.

The symbol $r$ is now freed for use, $q=2(p-1)$.

\newpage

\noindent
{\bf Lemma 7.9.} \textit{Let $X$ be a complex $X=S^0\cup_f e^{rq}$ where 
$p^2<r<p^2+p$. Then
\begin{eqnarray*}
BP^*(X) &\cong &\pi_*(BP)\mbox{ on one generator $h$ of degree 0 restricting}\\
&&\qquad\qquad\mbox{ to a generator on }S^0\\
&&\oplus\pi_*(BP)\mbox{ on one generator $\ell$ of degree $rq$ coming from $S^{rq}$}\\
\end{eqnarray*}
and
$$R_{p^2}h\equiv 0\mod\, p^p$$
(mod p suffices for what follows).}

\vspace{10pt}

\noindent
{\bf Proof.} From the $BP^*$ exact sequence of the cofibering, it is clear that we have a short exact sequence
$$0\longrightarrow BP^*(S^{rq})\longrightarrow BP^*(X)\longrightarrow BP^*(S^0)
\longrightarrow 0$$
Everything that follows would be equally valid for a short exact sequence of modules over $BP^*(BP)$
$$\xymatrix{
0\ar[r] &BP^*(S^{rq})\ar@{=}[d]\ar[r] &M\ar[r] &BP^*(S^0)\ar@{=}[d]\ar[r] &0\\
&N&&L}$$
Such a module defines an element of $\Ext^i_{BP^*(BP)}(L,N)$. It is clear that the structure of $M$ as a module over $\pi_*(BP)$ is as stated in the lemma; if we choose a generator $h$, then the operations in $M$ are given by a vector
$$\left[\begin{array}{l}R_1\\ R_p\\ R_{p^2}\end{array}\right]h$$
In $N^3=N\oplus N\oplus N$. We may alter $h$ to $h+a\ell$ ($a\in\pi_{rq}(BP)$); this alters
$$\left[\begin{array}{l}R_1\\ R_p\\ R_{p^2}\end{array}\right]h
\mbox{\qquad to\qquad }
\left[\begin{array}{l}R_1\\ R_p\\ R_{p^2}\end{array}\right]h +
\left[\begin{array}{l}R_1\\ R_p\\ R_{p^2}\end{array}\right]a\ell
$$
We thus get a description of
$\Ext^i$ as ``cocycles'' modulo ``coboundaries'' ; the cocycles are the vectors
$$\left[\begin{array}{l}R_1\\ R_p\\ R_{p^2}\end{array}\right]h;$$
the coboundaries are the vectors
$$\left[\begin{array}{l}R_1\\ R_p\\ R_{p^2}\end{array}\right]a\ell.$$

We first check that the result claimed does not depend on the choice of cocycle mod coboundaries. In fact, in the dimension considered $a$ is a sum of monomials
$$v_1^iv_2^j\mbox{ with }j\le p-1.$$
Consider $R_{p^2}(v_1^iv_2^j)\ell$ and expand it by the Cartan formula. We have
$R_kv_2= 0$ for $k> p+1$ and $R_{p+1}v_2=0$ mod $p^p$, which only leaves for consideration terms containing a factor $R_k(v_1^i)$ with $k\ge p$. Such a factor is $\equiv 0$ mod $p^k$, i.e., mod $p^p$.

Now we have a choice of arguments. First, we assume it known
that $\Ext^i_{BP^*(BP)}\left(\pi_*(BP),\pi_*(BP)\right)$ is in this dimension
 $\mathbf{Z}/p\mathbf{Z}$ generated by $\alpha_r$. Then we can easily work out the operations in the corresponding complex. We apply 7.4 to the composite
$$\xymatrix{
X\ar@{=}[d]\ar[r]^f &Y\ar@{=}[d]\ar[r]^g &Z\ar@{=}[d]\\
S^{rq-1}\ar[r] &S^{rq-1}V(0)\ar[r] &S^0
}$$

We obtain the following diagram
$$
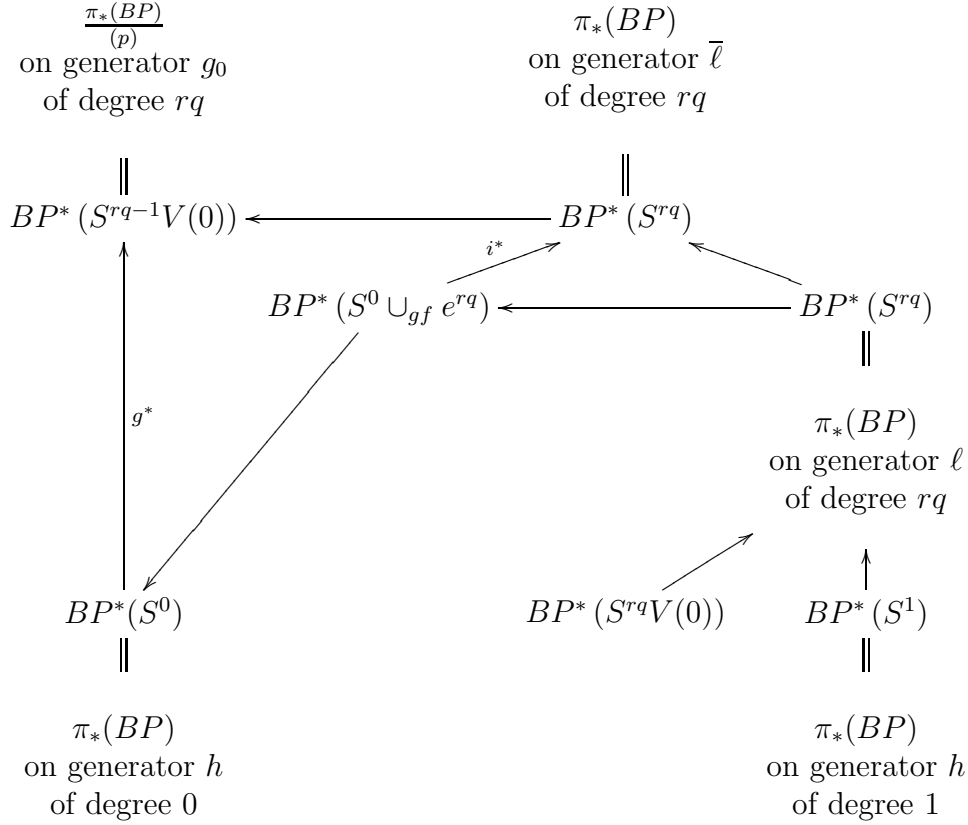
\xymatrix@C=5pt@R=15pt{
\parbox{80pt}{\begin{center}
$\frac{\pi_*(BP)}{(p)}$\\
on generator $g_0$\\
of degree $rq$
\end{center}}\ar@{=}[d]
&&\parbox{80pt}{\begin{center}
$\pi_*(BP)$\\
on generator $\overline{\ell}$\\
of degree $rq$
\end{center}}\ar@{=}[d]\\
BP^*\left(S^{rq-1}V(0)\right) &&BP^*\left(S^{rq}\right)\ar[ll]\\
&BP^*\left(S^0\cup_{gf}e^{rq}\right)\ar[ddl]\ar[ur]^{i^*}
&&BP^*\left(S^{rq}\right)\ar[ul]\ar[ll]\\
&&&\parbox{80pt}{\begin{center}
$\pi_*(BP)$\\
on generator $\ell$\\
of degree $rq$
\end{center}}\ar@{=}[u]\\
BP^*(S^0)\ar[uuu]_{g^*}&&BP^*\left(S^{rq}V(0)\right)\ar[ur]
&BP^*\left(S^1\right)\ar[u]\\
\parbox{80pt}{\begin{center}
$\pi_*(BP)$\\
on generator $h$\\
of degree 0
\end{center}}\ar@{=}[u]
&&
&\parbox{80pt}{\begin{center}
$\pi_*(BP)$\\
on generator $h$\\
of degree 1
\end{center}}\ar@{=}[u]\\
}$$

We see that
\begin{eqnarray*}
BP^*\left(S^0\cup_{gf}e^{rq}\right)
&= &\pi_*(BP)\mbox{ on one generator $h$ of degree 0 restricting}\\
&&\qquad\qquad\mbox{ to a generator on }S^0\\
&&\oplus\pi_*(BP)\mbox{ on one generator $\ell$ of degree $rq$ coming}\\
&&\qquad\qquad\quad\mbox{from $S^{rq}$}
\end{eqnarray*}
as in the statement of the lemma. The map $S^{rq}\to S^{rq}$ is part of the cofibering defining $V(0)$, and has degree $p$, so $\ell$ maps to $p\overline{\ell}$.

Also $g^*$ maps $h$ to $v_1^rg_0$ (by construction); so $i^*$ maps any
choice of $h$ to $v_1^r\overline{\ell}$ mod $p$. But $i^*$ maps $\ell$ to 
$p\overline{\ell}$, so we can change
the choice of $h$ to ensure that $hi=v_1^r\overline{\ell}$. Then
$$R_{p^2}hi=R_{p^2}v_1^r\overline{\ell}=cp^{p^2}v_1^{p-1}\overline{\ell}
=cp^{p^2}v_1^{p-1}\ell i,\quad c\in\mathbf{Z}.$$
Now $R_{p^2}h$ clearly maps to 0 on $S^0$ and on the complementary summand
$\pi_*(BP)\ell$ , $i^*$ is mono; so with this choice of $h$ we have
$$R_{p^2}h=cp^{p^2-1}v_1^{p-1}\ell.$$

Secondly, we can merely assume it known that $\Ext^i_{BP^*(BP)}(L,N)$
is finite. Let us denote
$$d_0 = \left[\begin{array}{l}R_1\\ R_p\\ R_{p^2}\end{array}\right]$$
Then in this case, for each cocycle $d_0h$ there is an integer $m$ such that
$md_0h$ is a coboundary $d_0(a\ell)$ where $a\in\pi_{qr}(BP)$. Now with 
$0 < j \le p-1$
the coboundary $d_0(v_1^iv_2^j\ell)$ is
$$\left[\begin{array}{cl}
-jv_1^{i+p}v_2^{j-1} &\mbox{mod }p\\
jv_1^{i+1}v_2^{j-1} &\mbox{mod }p^{p-1}\\
0  &\mbox{mod }p^p
\end{array}\right]\ell$$
while $d_0(v_1^r\ell)$ is
$$\left[\begin{array}{cl}
prv_1^{r-1}\\
0 &\mbox{mod }p^p\\
0 &\mbox{mod }p^{p^2}
\end{array}\right]\ell$$
Note $r\not\equiv 0$ mod $p$.

It is now clear that if $c,c_{ij}$ are rational coefficients, then the 
coboundary
$$d_0\left(cv_1^r\ell+\sum_{0<j\le p-1}c_{ij}v_1^iv_2^j\ell\right)$$
is integral only if the $c_{ij}$ are integral (i.e., lie in $\mathbf{Z}_{(p)}$) and $pc$ is integral, say $c=\frac{m}{p}$. Therefore the most general cocycle is
as obtained above. 

Finally, you can explicitly obtain the cocycles as the kernel of $d_1$, and if you think this method is more reliable than the two already presented you are welcome to go ahead.\hfill$\Box$

\vspace{10pt}

Next I undertake to define an invariant defined on elements of order $p$
in $\pi_{rq-2}(S^0)$, $p^2 < r < p^2 +p$. We know that $\beta_p$ is an element of order $p$, because it extends over $S^{rq-2}\cup_p e^{rq-1}=S^{rq-2}V(0)$
and even over $S^{rq-2}V(1)$
by construction. Toda's elements $\varepsilon_i$ of degree $(p^2+i)q-2$
 ($1\le i\le p-1$) are also asserted to be of order $p$. So suppose given a map 
$f:S^{rq-2}\to S^0$ of order $p$, and extend it to 
$\overline{f}:S^{rq-2}\cup_p e^{rq-1}\longrightarrow S^0$. Form the mapping cone
$$X=S^0\cup_{\overline{f}} C\left(S^{rq-2}\cup_p e^{rq-1}\right).$$
Then evidently we have a short exact sequence
$$\xymatrix{
0\ar[r] &BP^*\left(S^{rq-2}V(0)\right)\ar[r]\ar@{=}[d] &BP^*(X)\ar[r]
&BP^*(S^0)\ar@{=}[d]\ar[r] &0\\
&\parbox{100pt}{\begin{center}$\pi_*(BP)$ on\\ generator $g_0$\\ of degree $qr$
\end{center}}
&&\parbox{100pt}{\begin{center}$\frac{\pi_*(BP)}{(p)}$ on\\ generator $h$\\ of degree 0
\end{center}}
}$$
We take an element $h\in BP^0(X)$ projecting to the generator in $BP^*(S^0)$
and we consider $R_{p^2}h=ag_0$. Now changing $h$ to $h+ bg_0$ will of course
change $R_{p^2}h$ to $R_{p^2}h+R_{p^2}bg_0$, but we have already checked that
$R_{p^2}bg_0$ is zero mod $p^p$, and here we are working in $\frac{\pi_*(BP)}{(p)}$,
so the resulting change is zero. We must also not forget that we can change the extension $\overline{f}$
of $f$ by an element of $\pi_{qr-1}(S^0)$ but this changes our extension precisely by an element of the group $\Ext_{BP^*(BP)}^i(L, N)$ which I so carefully discussed in Lemma 7.9, and by Lemma 7.9, the resulting change in $R_{p^2}h$ is zero
mod $p$. This sets up our invariant.

It remains to calculate this invariant for $\beta_p$.

\vspace{10pt}

\noindent
{\bf Theorem 7.10.} \textit{If $f=\beta_p$, we have
$$R_{p^2}(h)=v_1^{p-1}g_0$$
and therefore $\beta_p\ne 0$.}

\vspace{10pt}

\noindent
{\bf Proof.} Factor $\beta_p$ in the form
$$\xymatrix{
X\ar[r]^f\ar@{=}[d] &Y\ar[r]^g\ar@{=}[d] &Z\ar@{=}[d]\\
S^{rq-2}V(0) &S^{rq-2}V(1) &S^0
}$$
where $g^*(\mbox{generator})= v_2^pg_1$. Apply 7.4. We obtain the following diagram
$$\xymatrix@R=12pt@C=15pt{
\parbox{80pt}{\begin{center}
$\frac{\pi_*(BP)}{(p,v_1)}$\\
on generator $g_1$\\
of degree $(r+1)q$
\end{center}}\ar@{=}[d]
&&\parbox{80pt}{\begin{center}
$\frac{\pi_*(BP)}{(p)}$\\
on generator $\overline{g}_0$\\
of degree $(r+1)q$
\end{center}}\ar@{=}[d]\\
BP^*\left(S^{rq-2}V(1)\right) &&BP^*\left(S^{(r+1)q-1}V(0)\right)\ar[ll]\\
\\
&BP^*\left(S^0\cup_{gf} C\left(S^{rq-2}V(0)\right)\right)\ar[uur]^{i_*}\ar[dl]
&BP^*\left(S^{rq-1}V(0)\right)\ar[uu]\ar[l]\\
BP^*(S^0)\ar[uuu]&&
\parbox{80pt}{\begin{center}
$\frac{\pi_*(BP)}{(p)}$\\
on generator $\overline{g}_0$\\
of degree $(r+1)q$
\end{center}}\ar@{=}[u]
}$$
The map $S^{rq-1}V(0)\longrightarrow S^{(r+1)q-1}V(0)$ comes from the cofibering defining $V(1)$
and therefore carries $g_0$ to $v_1\overline{g}_0$. So in our complex
$S^0\cup_{gf} C\left(S^{rq-2}V(0)\right)$ we have $g_0i=v_1\overline{g}_0$.

Let us choose a generator $h$ in 
$BP^*\left(S^0\cup_{gf} C\left(S^{rq-2}V(0)\right)\right)$ so that
it maps to a generator in $S^0$ and therefore to $v_2^pg_1$ in $BP^*(S^{rq-2}V(1))$. So
we must have
$$hi = v_2^p\overline{g}_0\mbox{ (mod $v_1$)}.$$
But changing $h$ by $ag_0$ changes $hi$ by $v_1a\overline{g}_0$ so we can choose the
generator $h$ so that it maps to $v_2^p\overline{g}_0$ exactly.

Now $R_{p^2}h$ clearly maps to zero on $S^0$. The complementary summand
$\frac{\pi_*(BP)}{(p)}g_0$ maps monomorphically under $i^*$, so it is sufficient to calculate $R_{p^2}hi$. We find
\begin{eqnarray*}
R_{p^2}hi &= &R_{p^2}v_2^p\overline{g}_0\\
&= &v_1^p\overline{g}_0\mbox{ mod $p^p$ at least}\\
&= &v_1^{p-1}g_0i.
\end{eqnarray*}
We conclude that
$$R_{p^2}h=v_1^{p-1}g_0$$
in $\frac{\pi_*(BP)}{(p)}g_0$.
\hfill$\Box$

\newpage

\section{Epilogue}\label{sec8}

Thomas and Zahler published their proof that $\gamma_1\ne 0$ in 1974, cf. \cite{TZ}. Their proof uses a tertiary operation
in Brown-Peterson cohomology to detect the nontriviality of the element. Oka and Toda tracked down the mistake in their
original calculation and published an alternative proof of this result in 1975, cf. \cite{OT}. The methods they employ involve
secondary compositions and extended powers of complexes. Later that year Miller, Ravenel and Wilson \cite{MRW} announced a
proof that $\gamma_t\ne 0$ for all $t>0$ and all $p\ge 7$, using the Adams-Novikov spectral sequence.

In 1975 Bousfield published a proof of Theorem 4.6, cf. \cite{B2}. His proof uses simplicial methods and proves a stronger
result. Bousfield constructs a strict functor on the category of simplicial sets which induces Adams' conjectured localisation
upon passage to homotopy. However it should be noted that there is a very simple way of repairing Adams' original approach
to localisation given in these notes, which we explain below.

\vspace{10pt}

We first introduce a new axiom on the class $S$ of morphisms we are inverting.

\vspace{10pt}

\noindent
{\bf Axiom 8.1.} For each pair of objects $X$, $Y$ in $C$ there  is a set of diagrams
$$\left\{\xymatrix{X\ar[r]^{f_\alpha} &Z_\alpha &Y\ar[l]_{s_\alpha}}\right\}$$
with $s_\alpha$ in $S$ such that for any diagram
$$\xymatrix{X\ar[r]^f &Z &Y\ar[l]_s}$$
with $s\in S$ we can find a morphism $t:Z_\alpha\longrightarrow Z$ in $S$ so that the
following diagram commutes
$$\xymatrix{
&Z_\alpha\ar[d]^t\\
X\ar[ur]^{f_\alpha}\ar[r]^f &Z &Y\ar[l]_s\ar[ul]_{s_\alpha}
}.$$

\vspace{10pt}

Axiom 8.1 serves as a substitute for Axiom 3.4.  Its virtue is that it is easily verified for the examples
of interest in these notes, in contrast to 3.4.  We next prove the relevant revision of Theorem 3.8.

\vspace{10pt}

\noindent
{\bf Theorem 8.2.} \textit{Let $\C$ be the category in which the objects are connected CW-complexes
with basepoint and the maps are homotopy classes. Let $S$ be a subclass of the morphisms of $\C$,
satisfying 3.1, 3.2, 3.3, 3.5, 3.6 and 8.1 (or alternatively 3.1, 3.5, 3.6, 4.1 and 8.1). Then $S$ arises by 2.6 from a pair $(E,\eta)$ satisfying 2.1 and 2.2.}

\vspace{10pt}

\noindent
{\bf Proof.} We first need to show that the quotient category $S^{-1}\C$ is well defined, i.e. 
$[QX,QY]_{S^{-1}\C}$ is a set for any pair of objects $X$, $Y$ in $\C$. We see this by showing
that any morphism in $[QX,QY]_{S^{-1}\C}$, represented by a diagram
$$\xymatrix{X\ar[r]^f &Z &Y\ar[l]_s}$$
in $\C$, has an equivalent representative
$$\xymatrix{X\ar[r]^{f_\alpha} &Z_\alpha &Y\ar[l]_{s_\alpha}}.$$
This is demonstrated by the following commutative diagram in $\C$:
$$\xymatrix{
&&Z\ar@{=}[drr]\\
X\ar[urr]^f\ar[drr]_{f_\alpha} &&Y\ar[u]^s\ar[d]_{s_\alpha}\ar[rr]^s &&Z\\
&&Z_\alpha\ar[urr]_{t}
}$$
The rest of the proof is identical to that of Theorem 3.8.  Briefly we use the Brown Representability
Theorem to construct a right adjoint to ${Q:C\to S^{-1}\C}$:
$$[QX, Y]_{S^{-1}\C}\longleftrightarrow[X,RY]_\C$$
and define $E$ to be the composite
$$\C\stackrel{Q}{\longrightarrow} S^{-1}\C\stackrel{R}{\longrightarrow}\C.$$
\hfill$\Box$

\vspace{10pt}

Finally we need to check that the class $S$ of our main example satisfies 8.1. This is demonstrated by the following lemma, which is essentially Lemma 11.2 of \cite{B2}.

\vspace{10pt}

\noindent
{\bf Lemma 8.3.} {\it Suppose that
$$X\stackrel{f}{\longrightarrow}Z\stackrel{s}{\longleftarrow} Y$$
is a diagram in $\C$, the homotopy category of based connected CW complexes. Suppose that $s$ induces an isomorphism
with respect to a generalized homology theory $K_*$. Then there is a commutative diagram
$$\xymatrix{
&W\ar[d]^{s''}\\
X\ar[ur]^{f'}\ar[r]_f &Z &Y\ar[l]^s\ar[ul]_{s'}
}$$
such that $s'$ and $s''$ induce isomorphisms with respect to $K_*$ and such that
$$\# W\le\max\{\# X, \# Y, \# K_*,\aleph_0\}.$$
Here $\# W$, $\# X$, $\# Y$ denote the cardinalities of the sets of cells in these CW complexes and $\# K_*$ denotes the
cardinality of $\bigoplus_{m\in\mathbf{Z}}K_m(\mathrm{pt})$.
}

\vspace{10pt}

\noindent
{\bf Proof.} Without loss of generality we may take $f$ and $s$ to be inclusions of complexes. We then
construct a sequence of subcomplexes 
$$W_0\subset W_1\subset W_2\subset\cdots\subset Z$$
with the following properties:
\begin{itemize}
\item[(i)] The inclusions $W_n\to Z$ induce epimorphisms in $K_*$.
\item[(ii)] $\mbox{ker}\left(K_*(W_n)\to K_*(Z)\right)=\mbox{ker}\left(K_*(W_n)\to K_*(W_{n+1})\right)$.
\item[(iii)] $\# W_n\le\max\{\# X, \# Y, \# K_*,\aleph_0\}$.
\end{itemize}
We proceed by induction. We take $W_0=X\vee Y$. Then (iii) is immediate and (i) follows from the fact that $Y\subset Z$
induces an isomorphism in $K_*$.

Having constructed $W_n$ we note $K_*(Z)=\colim_{\alpha}K_*(Z_\alpha)$, where $Z_\alpha$ runs over all subcomplexes
of $Z$ obtained by attaching finitely many cells to $W_n$.  Hence for each element $x\in\mbox{ker}\left(K_*(W_n)\to K_*(Z)\right)$
we can find such a subcomplex $Z_{\alpha_x}$ so that $x\in\mbox{ker}\left(K_*(W_n)\to K_*(Z_{\alpha_x})\right)$. We define
$W_{n+1}=\bigcup_{x}Z_{\alpha_x}$.  Then (ii) follows from construction and (i) holds for $W_{n+1}$, since it contains $W_0$.
To check that (iii) holds for $W_{n+1}$, we first observe that it follows from the induction hypothesis and the Atiyah-Hirzebruch
spectral sequence that $K_*(W_n)$, and hence also  $\mbox{ker}\left(K_*(W_n)\to K_*(Z)\right)$,  has cardinality bounded by
$\max\{\# X, \# Y, \# K_*,\aleph_0\}$, which implies the same bound on $\# W_{n+1}$.

Finally we define $W=\bigcup_{n=1}^\infty W_n$.  Since $K_*(W)=\colim_n K_*(W_n)$, it follows from (ii) that $K_*(W)\to K_*(Z)$
is injective. Since $W$ contains $W_0$, this map is also surjective, and hence an isomorphism.  It follows
that $Y\subset W$ is also an isomorphism in $K_*$.  The bounds on $\# W_n$ imply the same bound on $\# W$.  Thus the desired
commutative diagram is provided by the inclusions $X\subset W\subset Z$ and $Y\subset W\subset Z$.\hfill $\Box$  

\vspace{15pt}

If we take $S$ to be the class of morphisms in $\C$ which induce isomorphisms in $K_*$, as in Theorem 4.6, then  Lemma 8.3 shows that $S$ satisfies Axiom 8.1 (we may take the underlying sets of the $Z_\alpha$ in the axiom to be subsets of some fixed set of large enough cardinality).  Thus Theorem 8.2 implies Theorem 4.6. Note that we do not need to prove that $S$ satisfies Axiom 3.4. Indeed attempting to derive the proof of Theorem 4.6 from Axiom 3.4 appears to be a blind alley.

\vspace{20pt}

A natural follow-up question, which Adams did not address in these notes, is whether there are
localisations with respect to generalised cohomology theories.  It is clear that if one takes $S$ to be
the class of morphisms which are inverted by a generalised cohomology theory, then $S$ evidently satisfies
Axioms 3.1, 3.5, 3.6, and 4.1.  However the argument of Lemma 8.3 breaks down completely and
there does not appear to be any alternative argument to show that $[QX, QY]_{S^{-1}\C}$ is a set, and
thus no apparent way to show the existence of localisation with respect to $S$.
In \cite{B1} Bousfield shows that for ordinary cohomology theories (or more generally for anti-connective
cohomology theories) such localisations exist because the class $S$ is the same as the class of morphisms
inverted by $H_*(-;A)$ for $A$ a subring of $\mathbf{Q}$ or $A$ a subgroup of 
$\bigoplus_{\mbox{\scriptsize $p$ prime}}\mathbf{Z}/p\mathbf{Z}$.  In \cite{H},  Hovey similarly demonstrates the
existence of localisations with respect to a larger class of generalised cohomology theories by 
constructing corresponding generalised homology theories with the same classes $S$ of inverted morphisms. He moreover conjectures that for any generalised cohomology theory there is a corresponding
generalised homology theory such that the class of morphisms $S$ that they invert are the same.

However Casacuberta, Scevenels and J. H. Smith \cite{CSS} have shown that the existence of arbitrary
cohomological localisation functors is a logical consequence of a large-cardinal
axiom, Vop\v{e}nka's principle, and that a closely related conjecture (Farjoun's question)
can not be proved using the ordinary ZFC axioms of set theory.  A recent preprint of Bagaria, Casacuberta,
Mathias and Rosick\'y \cite{BCMR} proves the existence of cohomological localisations under a
much weaker large-cardinal hypothesis, namely the existence of a proper class of supercompact cardinals.  
They also show that the classes of maps inverted by a generalized homology theory may be defined by a
logical formula which has a lower formal complexity than a formula for defining classes inverted by a
generalized cohomology theory. This work seems to suggest that the existence of cohomological localisations
may be undecidable in ZFC, whereas Hovey's conjecture may be disprovable in ZFC.

The reader is also referred to \cite{B3} for a survey of results about localisation
functors in homotopy theory.

\end{document}